\documentclass[11pt,leqno]{article}
\usepackage{amsthm,amsfonts,amssymb,epsfig,graphics,amsmath,eufrak,oldgerm}
 \usepackage[latin1]{inputenc}\relax

\newcommand{\CalO}{\mathcal{O}}

\newcommand{\trans}{\text{\rm tr}}

\newcommand{\RR}{{\mathbb R}}

\newcommand{\CC}{{\mathbb C}}
\newcommand{\EE }{{\mathbb E}}

\newcommand{\FF}{{\mathbb F}}

\newcommand\cA{{\cal  A}}

\newcommand\cE{{\cal  E}}

\newcommand\cM{{\mathcal M}}

\newcommand{\rv}{{\mathrm{v } }}

\newcommand{\bdiag}{{\text{\rm block-diag} }}

\newcommand{\cz}{\check \zeta}
\newcommand{\cg}{\check \gamma}
\newcommand{\ch}{\check \eta}
\newcommand{\ct}{\check \tau}

\newcommand{\cx}{\check \xi}
\newcommand{\ccG}{\check G}
\newcommand{\ccA}{\check A}

\newcommand{\ccEE}{{\check {\mathbb E}}}
\newcommand{\ccSi}{\check \Sigma}

\newcommand{\up}{\underline p} 
\newcommand{\uq}{\underline q} 
\newcommand{\utau}{\underline \tau}
\newcommand{\uxi}{\underline \xi}
\newcommand{\ueta}{\underline \eta}
\newcommand{\uga}{\underline \gamma}
\newcommand{\uzeta}{\underline \zeta}

\newcommand{\ls}{{\ \lesssim \ }}

\def\eps{\varepsilon }

\def\D{\partial }
\newcommand\adots{\mathinner{\mkern2mu\raise1pt\hbox{.}
\mkern3mu\raise4pt\hbox{.}\mkern1mu\raise7pt\hbox{.}}}

\newcommand{\Id}{{\rm Id }}
\newcommand{\im}{{\rm Im }\, }
\newcommand{\re}{{\rm Re }\, }

\renewcommand{\div}{{\rm div}}
\newcommand{\rot}{{\rm curl}}
\newcommand{\curl}{{\rm curl}}
\newcommand{\na}{{\nabla}}

\newcommand{\la}{\langle }
\newcommand{\ra}{\rangle }
\newcommand{\pr}{{\parallel }}

\newtheorem{theo}{Theorem}[section]
\newtheorem{prop}[theo]{Proposition}
\newtheorem{cor}[theo]{Corollary}
\newtheorem{lem}[theo]{Lemma}
\newtheorem{defi}[theo]{Definition}
\newtheorem{ass}[theo]{Assumption}

\newtheorem{exam}[theo]{Example}
\newtheorem{rem}[theo]{Remark}
\newtheorem{rems}[theo]{Remarks}
\newtheorem{exams}[theo]{Examples}

\numberwithin{equation}{section}

\title{Hyperbolic Boundary Value Problems for Symmetric Systems with 
Variable Multiplicities}

\author{\sc \small Guy M\'etivier\thanks{
MAB Universit\'e de Bordeaux I,
33405 Talence Cedex , France; metivier@math.u-bordeaux.fr., 
partially supported by European network HYKE,  HPRN-CT-2002-00282.
},
Kevin Zumbrun\thanks{Indiana University, Bloomington, IN 47405;
kzumbrun@indiana.edu:
K.Z. thanks the University of Bordeaux I
for its hospitality during the visit in which
this work was carried out.
Research of K.Z. was partially supported
under NSF grants number DMS-0070765 and DMS-0300487.
 }}
\begin{document}

\maketitle

\begin{abstract}
We extend the Kreiss--Majda theory of stability of 
hyperbolic initial--boundary-value and shock
problems to a class of systems, notably including the
equations of magnetohydrodynamics (MHD), for which
Majda's block structure condition does not hold:
namely, simultaneously symmetrizable systems with
characteristics of variable multiplicity, satisfying
at points of variable multiplicity either a ``totally
nonglancing'' or a ``nonglancing and linearly splitting''
condition.
At the same time, we give a simple characterization of
the block structure condition as ``geometric regularity'' 
of characteristics, defined as analyticity of associated
eigenprojections
The totally nonglancing or nonglancing and linearly splitting conditions
are generically satisfied in the simplest case
of crossings of two characteristics, and likewise for
our main physical examples of MHD or Maxwell equations for a crystal.
Together with previous analyses of spectral stability
carried out by Gardner--Kruskal and Blokhin--Trakhinin,
this yields immediately a number of new results of nonlinear 
inviscid stability
of shock waves in MHD
in the cases of parallel or transverse magnetic field,
and recovers the sole previous nonlinear result, obtained 
by Blokhin--Trakhinin by direct ``dissipative
integral'' methods, of stability in the zero-magnetic field limit.
Our methods apply also 
to the viscous case.
\end{abstract}

\clearpage
\tableofcontents

\clearpage
\bigbreak

\section{Introduction}

There are two cases where the analysis of hyperbolic boundary value problems is well
developed: first, when the system is symmetric and the boundary conditions are dissipative
\cite{F.1, F.2, FL};
second, when the system is hyperbolic with constant multiplicities and the boundary conditions 
satisfy a Lopatinski conditions
\cite{Sa.1, Sa.2, Kr, MaOs, MetLMS}.
In both cases, the main energy estimate is proved using symmetrizers and 
integrations by parts. 
In the first case,  the $L^2$ estimate is immediate: the  symmetrizer is given 
and  the dissipativity  property is assumed.   In the second case, 
the construction of  symmetrizers is   delicate. 
It was achieved first by Kreiss (see also Chazarain-Piriou \cite{Ch-P}) 
in the strictly hyperbolic case.  
Kreiss' construction was next  extended by Majda-Osher and  
Majda to systems satisfying a ``block structure condition'',
 which  holds in particular 
as soon as the eigenvalues of the system are semi-simple and have constant 
multiplicity (see \cite{MetLMS}).
Kreiss's symmetrizers  are tangential  pseudo-differential operators. 
For symmetrizable systems with dissipative boundary conditions, a trivial,
constant symmetrizer is available in the form $SA_d$, where $S$ is the symmetrizer 
for the initial value problem, and $A_d$ is the normal coefficient matrix; 
see discussion above Theorem \ref{2main}.
 
The constant multiplicity assumption is satisfied in many applications
such as the linearized Euler equations or isotropic Maxwell' equation, but it is 
not  satisfied in other interesting example such as MHD or Maxwell's equation in crystals. 
On the other hand, the symmetry condition is  commonly satisfied as a consequence
of the existence of a conserved energy. For instance, the last two examples are symmetric. 
But  in the  symmetric case, the boundary conditions are not always dissipative, while the 
 sharp criterion of stability is given by a  Lopatinski condition
(see examples in Majda-Osher).  In particular, for shock waves 
considered as transmission problems, the Rankine-Hugoniot transmission conditions 
are not dissipative in general,  and Majda's  analysis of shock waves 
relies on the use of Kreiss symmetrizers.  
Therefore, with the goal of studying shock waves for   MHD, it is natural to 
extend the construction of Kreiss symmetrizers to   cases 
where the eigenvalues may have variable multiplicity. 

Such an  extension is  given below: for symmetric systems, we extend the 
construction of symmetrizers to a framework which contains MHD.  
First, we   classify the multiple eigenvalues 
into categories:  algebraically regular, geometrically  regular and nonregular. 
We show that Majda's block structure condition is equivalent to 
geometric regularity.  
This indicates clearly where Kreiss' construction stops.
Next, we extend the notion of being nonglancing 
to multiple eigenvalues.    Finally, using  the symmetry of the system,  
we extend the construction of   Kreiss' symetrizers   to 
the case of   eigenvalues that are not
geometrically regular provided that they are ``totally  nonglancing''
or else ``linearly splitting''
(defined in Section 3); see Theorem \ref{theo52}.

For applications to nonconstant coefficient or quasilinear systems,
it is important that symmetrizers vary smoothly with respect
to frequency, i.e., that the choice of symmetrizer be robust
under changes in the various parameters.
The basic construction for the totally nonglancing case is smooth.
However, the one in the linearly splitting case is not; thus, it
is useful mainly in the constant coefficient case.
A second possibility, under sufficiently strong structural assumptions,
would be to perform a second microlocalization about the locus 
of variable multiplicity; however, we do not pursue this direction, 
lacking a motivating physical example.

Instead, we pursue a different course, investigating in Section 6
the conditions under which one may obtain a smooth symmetrizer
in the original frequency variables, without further microlocalization.
This turns out to have a simple and pleasing answer, which also
clarifies the relation between Kreiss' theory and the earlier
Friedrichs theory of symmetrizable systems with 
dissipative boundary conditions.
Namely, one may define at points of variable multiplicity an appropriate
reduced system, with reduced boundary condition, involving only
the modes that change multiplicity.  Then, necessary and sufficient
conditions for existence of a smooth symmetrizer are that the full 
system satisfy the uniform Lopatinski condition, and the reduced 
problem admit a {\it constant symmetrizer},
or equivalently be symmetrizable with maximally dissipative boundary condition;
see Theorems \ref{2main}--\ref{3main}.
In the simplest case of a double root, for which the reduced system is 
$2\times2$,  it turns out that the Lopatinski conditions
is equivalent to symmetrizability/maximal dissipativity (Appendix D;  
see also \cite{MaOs}).  Thus, we find that
existence of a smooth symmetrizer requires only the Lopatinski
condition and nothing further; see Theorem \ref{4main}. 

Finally, in Section 7, we apply the tools we have developed
to the motivating problem of shock stability in MHD.
In this case, it turns out that all characteristics are at least
algebraically regular.  More, they are either geometrically
regular or totally nonglancing; see calculations, Appendix A.
Thus, our results yield existence of smooth symmetrizers,
and consequent linearized and nonlinear stability, provided
the uniform Lopatinski condition is satisfied, thus reducing
the stability problem as in the constant-multiplicity case
to a (nontrivial!) linear algebraic calculation.
This immediately gives several new results and recovers
the sole existing result, obtained by Blokhin and Trakhinin 
\cite{BT.1,BT.4}
by direct methods, of stability in the zero-magnetic field limit;
see Corollary \ref{pert}.

Though we concentrate here on inviscid, i.e., first-order problems,
the methods introduced extend in straightforward fashion to the
viscous case considered in 
\cite{ZKochel,ZCime, GMWZ.1, GMWZ.2, GMWZ.3, GMWZ.4}; see Appendix E.
In particular, this extends the results obtained in 
\cite{Zhandbook, ZCime} for compressible
Navier--Stokes equations to the case of (viscous) MHD.

\section{Multiple eigenvalues of hyperbolic systems}

In this section, we recall several basic definitions and start a classification of 
multiple eigenvalues. 

\subsection{Basic definitions}

Consider an $N \times N$ first order system with symbol
\begin{equation}
\label{defL} 
L(p, \tau , \xi) = \tau \Id + A(p, \xi) = \tau \Id + \sum_{j=1}^d \xi_j A_j (p). 
\end{equation}
The characteristic polynomial is 
$\Delta(p, \tau, \xi) = \det L(p, \tau, \xi)$.
Recall the following definition.

\begin{defi}
 (i)  The homogenous polynomial $ \pi  (\eta) $  is  hyperbolic in the real direction 
$\nu$  if and only if $\pi (\nu ) \ne 0$ and for all real 
$\eta' \notin \RR \nu$, all the roots in $z\in \CC $ of $\pi ( z \nu + \eta') = 0 $ are real. 

(ii)  The system $\eqref{defL}$ is hyperbolic in the direction $\nu$ if
its characteristic polynomial is. 

(iii) $L$ is symmetric hyperbolic in the direction $\nu$ in the sense of 
Friedrichs if there  is a smooth self adjoint matrix 
$S(p)$, called a Friedrichs symmetrizer,  such that all the matrices $S(p) A_j(p)$ are self adjoint
and $S(p) L(p, \nu)$ is definite positive. This implies hyperbolicity. 
  
\end{defi}

Below, we always assume that $L$ is hyperbolic in the direction 
$dt = (1, 0, \ldots, 0)$. 
This means that all the eigenvalues of $A(p, \xi)$ are real when 
$\xi $ is real.   
   
Suppose that $\Delta (\up, \utau, \uxi) = 0$. Then $- \utau$ is an eigenvalue of 
$A(\up, \uxi)$.   
Its   algebraic multiplicity $m$ is   the multiplicity of $-\utau$ as a root of the polynomial 
$\Delta(\up, \cdot, \uxi)$. 
Its  geometric multiplicty is  $ m_g =  \dim \ker L(\up, \utau, \uxi)$. 
The eigenvalue is semi-simple when $m_g = m$. 
Recall that for  symmetric hyperbolic  systems, the eigenvalues are always 
semi-simple.

The simplest cases of multiple roots are described in the following definition. 

\begin{defi}
\label{def22}
Consider a root $(\up, \utau, \uxi)$ of $\Delta(\up, \utau, \uxi)= 0$, of 
algebraic multiplicity $m$ in $\tau$. 

i)  $(\up, \utau, \uxi)$ is algebraically regular, if on a neighborhood $\omega$ of 
$(\up, \uxi)$ there are $m$ smooth real functions  
$\lambda_j(p, \xi)$, analytic in $\xi$,  such that for 
$(p, \xi) \in \omega$: 
\begin{equation}
\label{22}
\Delta (p, \tau,  \xi)  = e(p, \tau, \xi) \prod_{j=1}^m \big(\tau + \lambda_j(p, \xi) \big)
\end{equation}
where $e$ is a polynomial in $\tau$ with smooth coefficients such that 
$e(\up, \utau, \uxi) \ne 0$. 

ii) $(\up, \utau, \uxi)$ is geometrically regular if in addition 
there are $m$ smooth functions $e_j(p, \xi)$ on $\omega$ with values in 
$\CC^N$,  analytic in $\xi$, such that 
\begin{equation}
\label{23}
A(p, \xi) e_j(p, \xi) = \lambda_j(p, \xi) e_j(p, \xi), 
\end{equation}
and the $e_1, \ldots, e_m$ are linearly independent. 
\end{defi}
 
  \begin{rems}
{\bf a) }\textup{Simple roots  (i.e.  $m= 1$),  are  geometrically regular.}  

{\bf b) }\textup{Roots of constant multiplicity are algebraically regular.  In this case 
all the $\lambda_j$ are equal.  They are geometrically regular if and only if in addition  they are semi-simple : the $m= m_g$ vectors $e_j$  form a smooth basis of the eigenspace.  } 
\end{rems}

\begin{exams}
{\bf a) }\textup{For MHD, the multiple eigenvalues are algebraically regular, but some are not 
geometrically regular; see Appendix A. }

{\bf b) }\textup{For Maxwell's equation in a biaxial crystal, the multiple eigenvalues are not algebraically regular; see Appendix B. }
\end{exams}

 Systems with only geometrically regular multiple eigenvalues play an important role, since  
 they provide exactly the class of boundary value problems satisfying 
 the block structure condition (see Theorem \ref{theo31} below).  
To treat examples such as MHD or non isotropic Maxwell Equations,
we have to go beyond this class. 

\subsection{Block reduction}

Consider a root $(\up, \utau, \uxi)$ of $\Delta(\up, \utau, \uxi)= 0$, of 
algebraic multiplicity $m$ in $\tau$. 
When $(p, \xi)$ is close to $(\up, \uxi)$, $\Delta(p, \cdot, \xi)$ has exactly 
$m$ roots, counted with their multiplicity, close to $\utau$ and there is a smooth
block reduction of $A(p, \xi)$: 
\begin{equation}
\label{nn24}
 U^{-1} (p, \xi)      A(p, \xi) U(p, \xi) = \begin{pmatrix} 
    A^\flat (p, \xi)    &   0 \\
  0    &   \tilde A(p, \xi)
\end{pmatrix}  
\end{equation}
with $A^\flat$ of dimension $m$ and   $\tilde A$ has no eigenvalue in a neighborhood 
of $\- \utau$. Moreover, 
\begin{equation}
\label{n27}
A^\flat(\up, \uxi) = - \utau \Id \qquad \mathrm{ if } \ \utau \ 
\mathrm{is \ semi simple}. 
\end{equation}
Denoting by $V$ [resp. $W$ ] the $m$ first columns of $U$  
[resp. $m$ first rows of $U^{-1}$],  
 there holds, 
\begin{equation}
\label{n28}
 A(p, \xi) V(p, \xi)  = V(p, \xi) A^\flat(p, \xi) , 
 \quad   W (p, \xi) A(p, \xi)    =  A^\flat(p, \xi) W(p, \xi) 
 \end{equation}
\begin{equation}
\label{n28bb}
A^\flat(p, \xi)  = W (p, \xi )  A(p, \xi) V(p, \xi)    \qquad
W(p, \xi)   V(p, \xi) = \Id. 
\end{equation}
Note that when $L$ is symmetric hyperbolic, then 
$A^\flat$ is symmetric since we we can choose 
$U(p, \xi)$ such that $U^{-1}(p, \xi) = U^*(p, \xi)  S(p)$ implying that
\begin{equation}
\label{n28b}
W(p, \xi) = V^*(p, \xi) S(p) , \qquad 
A^\flat(p, \xi) = V^*(p, \xi) S(p) A(p, \xi) V(p, \xi). 
\end{equation}
There holds
\begin{equation}
\label{n29}
\Delta(p, \tau, \xi) = e(p, \tau, \xi) \det \big(\tau \Id + A^\flat(\tau, \xi) \big)
\end{equation}
with $e(\up, \utau, \uxi) \ne 0$.

\subsection{Examples  and linearly splitting eigenvalues} 

We  consider first the  example of double roots for 
symmetric systems. 
 In this case, the reduced matrix $A^\flat$
has the form
\begin{equation}
\label{n213}
A^\flat(p, \xi) = \lambda(p, \xi) \Id  + \begin{pmatrix}
  a(p, \xi)    &   b(p, \xi) \\
  b(p, \xi)     &   - a(p, \xi) 
\end{pmatrix}
\end{equation}
with $\lambda$, $a$ and $b$ smooth and homogeneous of degree $1$ in $\xi$.
Double roots occur on the set 
\begin{equation}
\label{n215}
\cM = \{ (p, \xi) : a(p, \xi) = b(p, \xi) = 0 \}
\end{equation} 
and  $(\up, \uxi) \in \cM$. 
Several  cases can be considered. 

\medbreak

{\bf 1) }  $a = b \equiv 0$; this is the constant multiplicity case. 

\medbreak

{\bf 2) }   $\cM $ is a smooth manifold of codimension 1, 
of equation $\{ \phi (p, \xi) = 0 \} $.  
  This means that $a$ and $b$ vanish on $\cM$, and one can factor out 
$\phi$, or powers of $\phi$ in $a$ and $b$. Assume for instance that 
\begin{equation}
\label{nn211}
a = \phi^k \tilde a , \quad  b = \phi^k \tilde b , \quad  \mathrm{with} \  ( \tilde a, \tilde b) \ne (0, 0) 
\ \mathrm{at} \ (\up, \uxi). 
\end{equation} 
This occurs with $k =1$ for some eigenvalues of MHD. 
Then the eigenvalues are $\lambda + \phi^k \tilde \lambda_j$ 
where the $\tilde \lambda_j$  are the smooth and distinct eigenvalues 
of 
$$
\tilde A =  \begin{pmatrix}
  \tilde a     &   \tilde b  \\
  \tilde b     &   - \tilde a  
\end{pmatrix}. 
$$
The eigenvectors of $A^\flat$ are those of $\tilde A$. 
We are in the geometrically regular case. 
 
\medbreak

{\bf 3) } 
$\cM$ is a smooth manifold of codimension 2 and more precisely
\begin{equation}
\label{n214}
d_\xi a(\up, \uxi) \wedge d_\xi b(\up, \uxi) \ne 0. 
\end{equation}  
 Then $a$ and $b$ can be taken as independent coordinates near
 $(\up, \uxi)$, transversal to $\cM$. The eigenvalues 
 are $\lambda \pm \sqrt{ a^2 + b^2}$.  
 $(\up, \uxi)$ is not an algebraically regular point. 
 This situation occurs for Maxwell's equations in  bi-axial crystals. 
 
 \medbreak
 {\bf 4) } We give now an example of a more degenerate 
 situation, which occurs for MHD. 
 $\cM$ is a manifold of codimension 2, 
 given by the equations
 \begin{equation}
\label{nn214}
\cM= \{  \phi = \psi = 0 \} , \quad   d \phi \wedge d \psi \ne 0 . 
\end{equation}
Moreover, $a$ and $b$  vanish at the second order on $\cM$, 
and there is a smooth function $c$ such that 
$a^2 + b^2 = c^2$. For instance, this holds when 
\begin{equation}
\label{nn215}
a = \phi^2 - \psi^2 , \quad b = 2 \phi \psi ,  \quad c = \phi^2 + \psi^2. 
\end{equation}
In this case the eigenvalues are smooth, equal to $\lambda \pm c$. 
 $(\up, \uxi)$ is  an algebraically regular point. 
 But it is not geometrically regular, since the eigenvectors are 
 \begin{equation*}
 \begin{pmatrix}
     \phi  \\
        \psi 
\end{pmatrix}, \quad \begin{pmatrix}
    -  \psi  \\
     \phi 
\end{pmatrix}
\end{equation*}
which have no limit as $(\phi, \psi) \to (0, 0)$. 

\medbreak 

These examples can be generalized to higher order roots. 
Suppose that we are given  a smooth conic manifold 
$\cM$ of codimension $\nu$ and a smooth function 
$\lambda$ on $\cM$  such that  for all $(p, \xi ) \in \cM$, 
$\lambda(p, \xi)$ is a semi-simple eigenvalue of $A(p, \xi)$ of multiplicity 
$m$. Suppose that $(\up, \uxi) \in \cM$ and $  \utau  + \lambda(\up,   \uxi) = 0$. 

\medbreak 

{\bf 1) } The trivial case $\nu = 0$ corresponds to constant multiplicity. 

\medbreak 

{\bf 2 a ) } 
Suppose  that there are no parameter $p$ and  $d = 2$. By homogeneity, 
$\nu$ is at most one. When $\nu = 1$, on the unit sphere 
$\vert \xi \vert = 1$, $\cM$ consists of isolated points,  and $A^\flat (\xi)$ 
is an analytic family depending on one parameter. If $L$ is symmetric, then so does 
$A^\flat$, and therefore $A^\flat$ is analytically diagonalizable: 
$(\utau, \uxi)$ is a geometrically regular characteristic point.  

\medbreak 

{\bf 2 b ) } Suppose that $\nu = 1$. Then $\cM$ is given
by an equation $\phi = 0$, and, extending $\lambda$ outside $\cM$, 
$A^\flat -  \lambda  \Id $ vanishes when $\phi = 0$.  
Suppose that 
\begin{equation}
\label{nn212e}
A^\flat (p, \xi) = \lambda(p, \xi) \Id + \phi^k \tilde A (p, \xi)
\end{equation}
where  $\tilde A (\up, \uxi) $ has distinct eigenvalues. 
This framework extends \eqref{nn211}; again 
$(\utau, \uxi)$ is a geometrically regular characteristic point.

\medbreak

{\bf 3) } 
Suppose that $\nu \ge 1$ and 
$\cM$ is given by equations 
$\phi_1 = \ldots = \phi_\nu = 0$ with 
$d_\xi \phi_1, \ldots, d_\xi \phi_\nu$ linearly independent. 
Then 
\begin{equation}
\label{nn213e}
A^\flat (p, \xi) = \lambda(p, \xi) \Id + \sum_{j=1}^\nu 
\phi_j \tilde A_j(p, \xi). 
\end{equation}
Extending  \eqref{n214}, consider the case where 
\begin{equation}
\label{nn214e}
\check A (\theta ) := \sum_j \theta_j \tilde A_j(\up, \uxi) \ \ \mathrm{is \ 
strictly \ hyperbolic }, 
\end{equation}
meaning that for all   $\theta \in \RR^{\nu} \backslash \{  0\} $, $\check A (\theta ) $ has only real and simple eigenvalues.  
Then, for $\vert \theta  \vert = 1$ and $(p, \xi)$ close to 
$(\up, \uxi)$,  $\sum \theta_j \tilde A_j (p, \xi)$ has distinct  simple eigenvalues 
$ \check \lambda _j ( p, \xi, \theta)$, with eigenvectors $ e_j (p, \xi, \theta)$. 
The eigenvalues of $A^\flat$ are 
\begin{equation}
\label{nn215e}
\lambda_j (p, \xi) = \lambda(p, \xi) +  \vert \phi \vert \check \lambda_j \Big( p, \xi, \frac{\phi}{\vert \phi \vert} \Big) . 
\end{equation}
The eigenvalues are thus simple away from $\cM$ and, in general, 
$(\up, \utau, \uxi)$ is not algebraically regular. 
 
In the next section, we give an intrinsic formulation of 
\eqref{nn214e}. We will refer to this example as  the \textit{linearly splitting} case. 

\subsection{The tangent systems at semi-simple multiple eigenvalues}

First, we recall   the following result about  multiple roots of hyperbolic polynomials. 

\begin{prop}
\label{propn22}
Assume that $\pi(\tau, \xi)$ is hyperbolic in the direction $dt = (1, 0, \ldots, 0)$. 
Consider  $(\utau, \uxi) \ne 0$   and assume that 
$\utau $ is a root of multiplicity $m\ge 1 $ of  $\pi (\cdot, \uxi) = 0$. 
Then 
\begin{equation}
\label{nnn215}
\pi (\utau + \tau,\uxi+  \xi)  =   \underline \pi^{(m)}  (\tau, \xi) + 
O( \vert \tau , \xi \vert^{m+1})
\end{equation}
where  $\underline \pi^{(m)}$ is homogeneous of degree $m$ in $(\tau, \xi)$ 
and  hyperbolic in the direction $dt$. 
\end{prop}

This means that all  terms of degree smaller that $m$
in  the Taylor  expansion of $\pi$  with respect to all variables, vanish. 
 
Consider a system \eqref{defL} and a root $(\up, \utau, \uxi)$ of $\Delta$, 
of multiplicity $m$. Denote by $A^\flat$ the associated $m \times m$ block 
as in \eqref{nn24} and by $W$ and $V$ the matrices defined in \eqref{n28}
\eqref{n28bb}. 
 Following Proposition \ref{propn22}, we also introduce the $m$-order term in the Taylor expansion of $\Delta(\up, \cdot) $  at $(\utau, \uxi)$, that we denote by $\underline \Delta$. 

\begin{prop}
\label{propnn27}
Assume that $-\utau$ is a semi-simple eigenvalue of 
$A(\up, \uxi)$ of order $m$. 
Then 
\begin{equation}
\label{nn216}
A^\flat (\up, \uxi + \xi) = - \utau \Id + \underline A' (\xi) + O( \vert \xi \vert^2) 
\end{equation}
with
\begin{equation}
\label{n211}
  \underline A'(\xi) = \sum_j
\xi_j \underline W  A_j (\up) \underline V , \quad 
\underline W= W(\up, \uxi), \quad \underline V= V(\up, \uxi)
\end{equation}
Moreover, there is $\underline e \ne 0$, such that 
\begin{equation}
\label{nn218}
 \det \big( \tau \Id + \underline A'(\xi) \big)
 = \underline e  \underline \Delta( \tau, \xi).  
 \end{equation}
\end{prop}

\begin{proof}
Differentiating the relations  \eqref{n28} \eqref{n28bb}  and using 
that $A(\up, \uxi) = - \utau \Id$, yields \eqref{nn216}. 

The determinant $\Delta^\flat (p, \tau, \xi)$ is a factor of 
$\Delta (p, \tau,  \xi)$, as in \eqref{n29}. 
By \eqref{nn216}, 
$$
\det\Big( (\utau + \tau) \Id + A^\flat ( \uxi + \xi) \Big) = 
\det \Big( \tau \Id + \underline A' (\xi) \Big) + O(\vert \tau, \xi \vert^{m+1}). 
$$
Comparing with \eqref{nnn215}, implies \eqref{nn218}. 
\end{proof}
  
  We will refer to   $\underline L'  = \D_t + \underline A' (\D_x)$ as the 
  \textit{tangent system to $L$ at $(\up, \utau, \uxi)$}.

\begin{rem}
\label{groupeq}
 \textup{ The operator 
 $\underline L' $ 
 is well known in the analysis of the propagation singularities or of   wave packets 
 \begin{equation*}
u (t, x) = e^{ i (t \utau + x \uxi)/ \eps}  \big( a_0(t, x) + \eps a_1(t, x) + \ldots \big)
\end{equation*}
solutions to $L(\up, \D_t, \D_x) u = 0$. The propagation  is given by the laws of geometric optics. 
The principal term $a_0$ satisfies  
\begin{equation}
\label{GO}
  a_0 =  \underline V  a^\flat_0 \,, \quad   \  \underline L'  (\D_t, \D_x)  a^\flat _0 = 0
\end{equation}
(see \cite{Lax, JMR} ;  see also \cite{Texierthesis, Texieradvances} 
for \eqref{nn218}). 
This means that  the high frequency oscillations  are transported 
by $\underline L'  (\D_t, \D_x)$.  }
 \end{rem}

 \begin{rem}
  \label{remnn29e}
  \textup{By homogeneity, $\underline A' (\uxi) = - \utau \Id$ and 
  $(-\utau, \uxi)$ is always a characteristic direction  for  $\underline L'$.}
  \end{rem}

\begin{exam}
\label{ex210}
 \textup{If  $\utau + \lambda(\uxi) = 0 $ 
with $\lambda$  an eigenvalue of $A(\xi)$ of constant multiplicity $m$ 
near $\uxi$, then
\begin{equation}
\label{defv}
\begin{aligned}
 & \underline L' = ( \D_t +   \rv \cdot \D_x ) \Id  , 
\qquad \mathrm{with } \   \rv =  \D_{\xi } \lambda (\uxi) ,  
\\
& \underline \Delta (\tau, \xi) = \beta (\tau +  \rv \cdot \xi)^m , \quad
\beta = \frac{1}{m!} \D_\tau^m \Delta(\up, \utau, \uxi). 
\end{aligned}
\end{equation}
$\underline L'$ is the transport field at the group velocity  $\rv$.  
 In this case, \eqref{GO} means that the amplitudes $a_0$ are transported along
the rays of geometric optics, which are the integral curve of the transport field
$\D_t + \rv \cdot \D_x$. }
\end{exam}

\begin{exam}
\label{ex211}
\textup{Suppose that  $(\utau, \uxi)$ is geometrically regular. With $\lambda_j$ and 
$e_j$ as in \eqref{22} \eqref{23},  in the basis $e_j(\uxi)$, 
\begin{equation}
\label{212}
\begin{aligned}
& \underline L'  = \mathrm{diag} (\D_t + \rv_j \cdot \D_x) , 
\qquad  \rv_j =  \D_\xi \lambda_j(\uxi), 
\\
& \underline \Delta (\tau, \xi) = \beta  \prod (\tau +  \rv_j  \cdot \xi)^m , \quad
\beta = \frac{1}{m!} \D_\tau^m \Delta(\up, \utau, \uxi). 
\end{aligned}
\end{equation} 
 }
\end{exam}

\begin{exam}
\label{ex212}
\textup{ Suppose that $A^\flat$ has the form   
  \eqref{nn213e}, with $d_\xi \phi_1, \ldots, d_\xi \phi_\nu$ linearly independent. 
Using the notation $\check A$ as in \eqref{nn214e}, the symbol of 
$\underline L'$ is 
\begin{equation}
\label{nn226e}
\underline L' (\tau, \xi) =  (\tau + \rv \cdot \xi) \Id  +  \check A ( \Phi \xi) 
\end{equation} 
with $\rv = \D_\xi \lambda(\up, \uxi)$ and $\Phi = \D_\xi \phi (\up, \uxi)$. 
$\Phi \xi = 0$ if and only if $\xi$ is tangent to 
$\cM_{\up} := \{ \phi_1(\up, \cdot)  = \ldots = \phi_\nu(\up, \cdot)  = 0\} $. 
Thus, the condition \eqref{nn214e} is equivalent to
\begin{equation}
\label{nn228e}
 \mathrm{ for \ all \ }   \xi  \notin T_{\uxi}  \cM_{\up}  , 
 \    \underline A' (\xi)  \  \mathrm{ has \ only  \ real \ and \ simple \ eigenvalues}.  
\end{equation} 
This motivates the following definition. }
\end{exam}

\begin{defi}
\label{def213}
A multiple root $(\up, \utau, \uxi)$ is called linearly splitting transversally to 
  a smooth manifold $\cM = \{ \phi_1 = \ldots = \phi_\nu = 0 \} $ with
  $  \phi_1, \ldots,  \phi_\nu$ analytic in $\xi$ and 
  $d_\xi  \phi_1, \ldots, d_\xi \phi_\nu$ linearly independent, if $\cM$ contains
   $(\up,   \uxi)$ and 
   
   \quad i)  there is a smooth real valued function $\lambda(p, \xi)$, analytic in $\xi$,
   such that for all $(p, \xi) \in \cM$, $\lambda(p, \xi)$ is a semi-simple eigenvalue
   of $A(p, \xi)$ of constant multiplicity $m$, 
   
   \quad ii) the condition $\eqref{nn228e}$  is satisfied. 
    \end{defi} 
    
\begin{rem}
\label{rem213}
\textup{In this case, the manifold 
 \begin{equation}
 \label{nn229e} 
 \widetilde \cM = \{ (p, \tau, \xi) \ : \ (p, \xi) \in \cM,  \ \tau = - \lambda(p, \xi) \} 
 \end{equation}
 is a smooth submanifold of the characteristic variety of $L$ which contains
 $(\up, \utau, \uxi)$. The corresponding  block   
 \begin{equation} 
\label{nn230e}
L^\flat (p, \tau, \xi ) := \tau \Id + A^\flat (p, \xi) 
\end{equation} 
 satisfies $ L^\flat (p, \tau, \xi ) = 0$ on $\widetilde \cM$. Therefore, 
  \begin{equation} 
\label{nn231e}
\underline L' ( \tau, \xi ) = 0 , \quad \mathrm{for} \  (\tau, \xi) \in T_{\utau, \uxi}  \widetilde \cM_{\up}.
\end{equation} 
 This means that  $\underline L' $ only depends on frequency variables which are transversal 
 to $\widetilde \cM$. Introduce the quotient 
 $ \check E = \RR^{1+d} / T_{\utau, \uxi}  \widetilde \cM_{\up}$, 
 and denote by $\varpi$ the projection from $\RR^{1+d}$ onto $\check E$. By \eqref{nn231e}, 
 there is 
 $\check L$ on $\check E$ such that  
 \begin{equation}
 \label{nn232e}
 \underline L' (\tau, \xi) = \check L \big(\varpi(\tau, \xi)\big) 
 \end{equation}
 Because $dt$ is transversal to $\tilde \cM$, $\varpi dt \ne 0$, and the condition 
 \eqref{nn228e} is equivalent to 
 \begin{equation}
\label{nn214ee}
\check L (\theta), \ \theta \in \check E,   \ \mathrm{is \ strictly \ hyperbolic \ in \ the \ direction \ } \varpi dt.  
\end{equation}
This approach gives a completely intrinsic definition of linear splitting. 
In particular, it shows that one can replace $dt$  by  any hyperbolic direction. }

\textup{
Alternatively, consider  $E$ such that 
$\RR^{1+d} = \EE \oplus T_{\utau, \uxi}  \widetilde \cM_{\up}$ and a basis in $E$. 
Identifying $E$ and $\check E$, 
 the operator $\check L$ has the form 
$  \check L =  \sum_{j} \check A_j \D_{y_j}$. 
One can choose $E$ and the basis such that $dt = dy_0$; in this case   
    \eqref{nn214ee} or \eqref{nn228e} means that  $  \check L$ is strictly 
hyperbolic in the direction $dy_0$. 
For example, for Maxwell's equation in a bi-axial crystal, the multiplicity  is two and 
$ \underline L'  $ is   a  $ 2\times 2$,   strictly hyperbolic system in space dimension 2, equivalent to a  2-D wave equation (see Appendix B).  That one space dimension (at least) is lost when 
passing from $L$  to $\underline L'$ is a general  fact that  follows from Remark~\ref{remnn29e}. 
}
 \end{rem}

     \begin{rem}
     \label{rem214y}
     \textup{If $\cM$ is of codimension one, then \eqref{nn212e} holds with $k= 1$, and 
     \eqref{nn228e} implies that $\tilde A(\up, \uxi) $  has distinct eigenvalues. This shows that 
   if   $(\up, \utau, \uxi)$ is  linearly splitting transversally to 
  a smooth manifold $\cM$ of codimension 1, then $(\up, \utau, \uxi)$ is geometrically regular. }
     \end{rem}



\section{The boundary block analysis}

In this section we start the investigation of noncharacteristic boundary value problems
for hyperbolic systems. In particular, we propose an extension of  the definition of
 glancing  to multiple eigenvalues with variable multiplicity. 

We consider a planar boundary. Changing notations and calling 
$(y,x) $ the spatial components,  the boundary is  $\{ x    = 0\} $.  
We assume that the boundary is not characteristic, that is 
\begin{equation}
\label{n31}
\det A_d \ne 0. 
\end{equation} 
The spatial Fourier frequency variables are denoted by $(\eta, \xi)$. 
and  $\tau - i \gamma $ is the complex time Fourier-Laplace  frequency. 
Given a hyperbolic system $L$ \eqref{defL}, we consider
 \begin{equation}
\label{nn32}
G(p, \zeta) = A_{d}(p)^{-1}\Big( (\tau- i \gamma) \Id + \sum_{j=1}^{d-1} 
\eta_j A_j(p) \Big)
\end{equation}
with  $\zeta = (\tau, \eta, \gamma) \in \RR^{d+1}$. We denote by 
$\RR^{d+1}_+ $ the half space $\{ \gamma > 0 \}$. By homogeneity, we can 
restrict attention to $\zeta $ in the unit sphere $S^d$; we denote by 
$S^d_+$ the open half sphere $S^d \cap \{ \gamma > 0 \}$ and by 
$\overline S^d_+$ its closure $S^d \cap \{ \gamma \ge 0 \}$.

Denoting by $\Delta$ the   characteristic polynomial  of $L$, there holds
\begin{equation}
\label{31}
\Delta(p, \tau - i \gamma, \eta, \xi) = \det A_d(p) \ 
\det \big( \xi \Id + G(p, \zeta) \big).  
\end{equation}  
The consequence of hyperbolicity, is that for $\gamma \ne 0$, 
$G(p, \zeta)$ has no real eigenvalue.


\subsection{Block reduction}

The goal is to construct symmetrizers for $G(p, \zeta)$. 
This is done locally, near each point 
$(\up, \uzeta)$ with $\uzeta \in \overline S^d_+$ using  a diagonal block reduction 
of  $G(p, \zeta)$. 

Accordingly, consider an eigenvalue $- \uxi$  of
$G(\up, \uzeta)$,   $\uzeta \in \overline S^d_+$. 
Denote by $m'$ its algebraic multiplicity.  There are 
$N \times m'$, $ m' \times m'$ and  $m' \times N$ matrices,  $V_b(p, \zeta)$, 
$G^\flat(p, \zeta)$ and 
$W_b(p, \zeta)$ respectively, depending smoothly on $(p, \zeta) $ near $(\up, \uzeta)$, 
and holomorphically on $\tau - i \gamma$, such that 
\begin{eqnarray}
\label{nn34}
 G   V_b   = V_b  G^\flat,  
 &\quad  & W_b    G    =  G^\flat  W_b , 
\\
\label{nn35}
G^\flat   = W_b  G V_b ,    &\quad & 
W_b    V_b = \Id. 
\end{eqnarray}
  Moreover, the spectrum of 
  $G(\up, \uzeta)$ is $\{ - \uxi \}$.  
  
  Hyperbolicity of $L$ implies the  following. 
  
 \begin{lem}
 \label{lem31}
 With notations as above, $G^\flat$ has no real eigenvalues when $\gamma > 0$.  
  \end{lem}
 
Consider $\uzeta = (\utau, \ueta, \uga) \ne 0  $ with   $ \uga = 0$, and a \textit{real} eigenvalue 
$- \uxi$ of $G^\flat(\up, \uzeta)$. Thus, $(\up,  \utau, \ueta, \uxi)$ is a real root of $\Delta$. 
Denote by $m$ its order in $\tau$ and by $A^\flat(p, \eta, \xi) $ the $m \times m$ matrix 
corresponding to this root in the block reduction of $A(p, \eta, \xi)$. 
There is no explicit relation between $A^\flat$ and $G^\flat$. However, there holds: 

\begin{lem}
\label{lem32}
If  $(\up,  \utau, \ueta, \uxi)$ is a real root of $\Delta$, then the polynomial in $\xi$
$\det (\xi + G^\flat(p, \zeta)$ has real coefficients when $\gamma = 0$. 

Moreover, for  
$(p, \tau - i \gamma , \eta, \xi)$ in a neighborhood of 
  $(\up,  \utau, \ueta, \uxi)$, complex with respect to the second and fourth argument,  there holds
\begin{eqnarray}
\label{nn36}
 \det \big(   \xi  \Id +   G^\flat (p,   \zeta )  \big) &
= &  e (p,   \zeta,  \xi)\det \big( (  \tau - i   \gamma) \Id +   A^\flat (p,   \eta,  \xi)  \big) , 
\\
\label{nn37}
\ker \big(\xi \Id  + G^\flat(p, \zeta)\big) & = & E(p, \zeta, \xi) 
\ker \big( (\tau - i \gamma )\Id + A^\flat(p, \eta, \xi) \big) , 
\end{eqnarray} 
with $e (\up, \uzeta, \uxi) \ne 0$ and  $E (p, \zeta, \xi) = W_b(p, \zeta) V(p, \eta, \xi) $.

If in  addition $- \tau$ is a semi-simple eigenvalue of $A(\up, \ueta, \uxi)$, then 
$m' \ge m$ and   $E (\up, \uzeta, \uxi) $ is of maximal  rank $m$. 
\end{lem}

 \begin{proof}
By hyperbolicity $\det A_d $ is real and $\det (\xi \Id + G)$ has real coefficients
 when $\gamma = 0$. If  $\uxi $ is real,  then for $(p, \tau, \eta, 0)$ close to $(\up, \uzeta)$, 
 the roots close to  $\uxi$ go by conjugated pair, implying 
 that $\det (\xi \Id + G^\flat  (p, \tau, \eta, 0)) $ has real coefficients. 
 
 By definition of $A^\flat $ and $G^\flat$,  
 both determinant in \eqref{nn37} are equal to $\Delta$   up to a nonvanishing factor near 
$(\up, \utau, \ueta, \uxi)$. 

Moreover, 
  \begin{equation*}
\label{412c}
\begin{aligned}
V(p, \eta, \xi) \ker \Big(  (  \tau - i  \gamma) \Id +  A^\flat (p,  \eta,   \xi)  \Big)
& = \ker \Big(  (  \tau - i  \gamma) \Id +  A(p,  \eta,   \xi)  \Big)
\\
= \, \ker \Big( \xi  \Id +  G(p,  \zeta)  \Big)  & = 
V_b(p, \zeta) \ker \Big(  \xi  \Id +  G^\flat (p,  \zeta)  \Big). 
\end{aligned}
\end{equation*}
This implies \eqref{nn37}. If $- \utau$ is semi-simple, then 
$ \underline \EE := \ker \big( \utau \Id + A(\up, \ueta, \uxi) \big)$ has dimension 
$m$ and is the image of $V(\up, \ueta, \uxi)$. It is equal to $\ker (\uxi \Id + G^(\up, \uzeta))$
thus contained in the invariant space which is the image of $V_b(\up, \uzeta)$. 
This implies that $m' \ge m$ and that  $W_b(\up, \ueta) $ is injective on $\underline \EE$.
 \end{proof}

 
 \subsection{The block structure condition}
 
 It was pointed out by A.Majda, \cite{Maj} (see also \cite{MaOs}) that Kreiss' construction 
 of symmetrizers can be carried out if $G$ satisfies the following 
 \textit{block structure condition}.

\begin{defi}
\label{defbl2} A matrix 
$G(p, \zeta) $ has the block structure property near $(\up, \uzeta)$ if there exists a smooth matrix 
$V$ on a neighborhood of that point such that 
$V^{-1} G V = \mathrm{diag} (G_k) $ is block diagonal, with blocks $G_k$
of size $\nu_k \times \nu_k$,    having one 
the following properties:

 i) 	  the spectrum of $G_k(p, \zeta) $ is contained in $\{ \im \mu \ne 0 \}$. 

ii)  $\nu_k = 1$,  $G_k(p, \zeta)$  is real when when 
$\gamma = 0$, and  
$ \partial_\gamma G_k(\up, \uzeta)  \ne 0$, 

iii)   $\nu_k > 1$, $G_k(p, \zeta)$  has real coefficients   
  when $\gamma = 0$,   there is
 $ \underline \mu_k  \in \RR$ such that
 \begin{equation}
\label{bl2}
 G_k(\up, \uzeta)    = \underline \mu_k \Id +  
\left[\begin{array}{cccc}
0  & 1 & 0&   
\\
0  &0  & \ddots  &  0
\\
  & \ddots &      \ddots & 1 \\
 &  &  \cdots & 0
\end{array}\right]\, , 
\end{equation}
 and  the lower left hand corner of 
${\partial_\gamma  G_k} (\up, \underline \zeta)$ 
does not vanish. 

\end{defi}

Consider a block $G^\flat $ associated to an eigenvalue $-\uxi$ of 
$G(\up, \uzeta)$. 
If the block is \textit{elliptic},  that is if  $\uxi $,  then $G^\flat$ satisfies 
property $i)$ near $(\up, \uzeta)$, and this finishes the analysis when $\uga > 0$.  
 
 \begin{theo}
 \label{theo31}
 Consider   $\uzeta = (\utau, \ueta, 0) $ and 
 $- \uxi $  a real eigenvalue of $G(\up, \uxi)$. Then, the associated  block 
 $G^\flat (p, \zeta)$ satisfies Majda's block structure condition on a neighborhood
 of $(\up, \uzeta)$ if and only if $(\utau, \ueta, \uxi)$ is geometrically regular for 
 $L$ is the sense of Definition $\ref{def22}$. 
 \end{theo}
 
 The if part was proved by Kreiss 
for strictly hyperbolic systems. It was next extended 
to hyperbolic systems
 with semi-simple eigenvalues of constant multiplicity 
  in \cite{MetLMS}. This  proof, which  only uses  the factorization     
    \eqref{22} and the existence of smooth eigenvectors \eqref{23},
 is easily extended to the geometrically regular case.
 Surprisingly, the condition of geometric regularity is also necessary.   
 We postpone the independent proof of this theorem to  Appendix C. 
  
 
 \subsection{Glancing modes}
 
 We want to extend the definition of \textit{glancing}  to general eigenvalues. 
 Recall   first  the case of  a simple eigenvalue  of $A(p, \ueta, \uxi)$ or of an  eigenvalue
  of constant  multiplicity $m$. Suppose that the real eigenvalue $- \uxi$ of $G(\up, \zeta)$ 
  satisfies 
 \begin{equation} 
 \label{33}
   \utau + \lambda(\up, \ueta, \uxi) = 0. 
 \end{equation}
 Following \cite{Kr} and \cite{MetLMS}, the  analysis of $G^\flat $ depends on the 
multiplicity  of $\uxi$ as root of \eqref{33}. 
\textit{Hyperbolic modes} correspond to  simple roots 
and \textit{glancing modes} to multiple roots. Only the later  yield non trivial Jordan blocks in \eqref{bl2}. 
The glancing condition  reads
\begin{equation}
\label{34}
\D_\xi \lambda(\up, \ueta, \uxi) =   0. 
\end{equation}
Introducing the transport field  \eqref{defv} and  the speed 
$\rv = (\rv_1, \ldots, \rv_d) = \D_{\eta, \xi} \lambda(\up, \ueta, \uxi) $, \eqref{34}  is equivalent to 
\begin{equation}
\label{35}
\rv_d  =   0. 
\end{equation}
 In this case, the rays  which are the integral curves of the transport fields are tangent 
 to the boundary.  
 According to Example  \ref{ex210}, the conditions 
 \eqref{34} or \eqref{35} are satisfied if and only if the  
 boundary is characteristic for 
 $ \underline \Delta^{(m)} $,  the principal term of $\Delta$ at 
 $({\up, \utau, \ueta, \uxi})$. 
 This motivates  the  following definition. 

 \begin{defi}
 \label{def32}
A root $( \up, \utau, \ueta, \uxi)$ of $\Delta$, of multiplicity 
 $m$ in $\tau$  
  is nonglancing, if and only if the $m$-th order term in the Taylor 
  expansion of $\Delta (\up, \cdot) $ at $(\utau, \ueta, \uxi)$ satisfies  
 \begin{equation}
 \label{nn310}
 \underline \Delta^{(m)}   (dx  ) \ne 0
 \end{equation} 
  where $ dx = (0, \ldots, 0, 1) $ is the (space-time) conormal to the boundary. 
 \end{defi} 
 
 Also recall that in the constant multiplicity analysis, the two cases
 $\rv_d > 0$ and $\rv_d < 0$    are different. 
  In the first case, the rays launched at  the boundary $\{x = 0\}$ at time 
 $t = 0$ enter the domain $\{ x > 0 \} $ for  $t > 0$;  
 the mode is said  \textit{incoming}. In the second case, the rays go off the domain, 
 and the mode is   \textit{outgoing}. 
 This has a  natural extension to general eigenvalues. 
 By Proposition \ref{propn22}, $\underline \Delta^{(m)}$ is hyperbolic in the direction 
 $dt$. Introduce the component of $dt$ in 
 $\underline \Delta ^{(m)} > 0$. It is an open  convex cone 
 $\underline \Gamma_+ $  and $\underline \Delta^{(m)} $ is hyperbolic
 in any direction in $\underline \Gamma_+$. Its dual cone 
 $\underline{\hat \Gamma}_+ = \{  (t, y, x) : \forall (\tau; \eta, \xi) \in 
 \underline \Gamma_+ , \  
 t \tau + y \eta + x \xi \ge 0 \}$  is the   forward cone of propagation
 (see e.g. \cite{Hormander}). 
 In the constant multiplicity case, 
 $\underline \Gamma_+   = \{ \tau + \rv' \eta + \rv_d \xi > 0\}$  and 
 $ \underline {\hat \Gamma}_+  = \overline \RR_+ (1, \rv) $. 
 Thus,    $  \rv_d >  0$ if and only if   to one of 
  the following two equivalent conditions hold: 
 \begin{eqnarray}
 \label{36}
 & &   dx   \in \underline \Gamma_+   , 
 \\
\label{37}
 &  &  \underline {\hat \Gamma} _+  \backslash\{0\}   \subset \{ x > 0\}  . 
\end{eqnarray}
 Similarly,  the outgoing nonglancing condition $  \rv_d <   0$ is equivalent to one of the following conditions. 
 \begin{eqnarray}
 \label{38}
 & &  -  dx    \in \underline \Gamma_+  , 
 \\
\label{39}
 &  &  \underline {\hat \Gamma} _+ \backslash\{0\}   \subset \{ x <  0\}  . 
\end{eqnarray}
 
 This leads to the following extension. 
 
 \begin{defi}
 \label{def33}
A root  $( \up, \utau, \ueta, \uxi)$ 
  of $\Delta$, of multiplicity $m$ in $\tau$, is said to be 
 totally incoming  if one of the two equivalent conditions $\eqref{36}$ 
 $\eqref{37}$ is satisfied. It is totally outgoing if the equivalent conditions $\eqref{36}$ 
 $\eqref{37}$ hold.  It is totally nonglancing if it is either totally incoming or outgoing. 
 \end{defi} 
 
 \begin{exam}
 \textup{
 We have already seen the example of constant multiplicity \eqref{33}. 
 Consider  now the case of an algebraically regular root
 as in Example \ref{ex211}. Then 
 \begin{equation*}
\underline  \Delta   = \beta 
 \prod_j \big(\tau + \rv_j' \cdot \eta  + \rv_{j, d} \xi\big) , 
 \end{equation*}
see \eqref{212}. 
 Each mode  can be glancing, incoming or outgoing 
 depending on the sign of  $\rv_{j, d}$. 
 According to Definition \ref{def32}, the mode is nonglancing  when all the 
 $\rv_{j, d} $ are different from $0$. 
It is totally incoming [resp., outgoing] if  all  the  $\rv_{j, d} $ are positive
[resp., negative].  This explains    the terminology 
 totally incoming or totally outgoing.} 
 \end{exam}
 
 \bigbreak
When  $- \utau$ is  a semi-simple eigenvalue of 
  of $A(\up, \ueta, \uxi)$ of multiplicity $m$, we have introduced the tangent 
  system $\underline L' (\tau, \eta, \xi) $ at $(\up, \utau, \ueta, \uxi)$. 
 With notations as in \eqref{n211}, it reads
  \begin{equation}
\label{n312}
\underline L' (\tau, \eta, \xi) = \tau \Id + 
\sum_{j= 1}^{d-1} \underline A_j \eta_j + \underline A_d \xi \,, 
\qquad \underline A_j = \underline W A_j (\up) \underline V. 
\end{equation}
 
 \begin{lem}
 \label{lem35}
 Assume that  $- \utau$ is  a semi-simple eigenvalue of $A(\up, \ueta, \uxi)$. 
 Then
 
 i)  $ {\up, \utau, \ueta, \uxi}$ is nonglancing if and only if 
 the boundary is noncharacteristic for 
   $\underline L'$, which means that
   $\underline A_d$ is invertible,

 ii) $ {\up, \utau, \ueta, \uxi}$ is totally incoming [resp., outgoing]  if and only 
   $\underline L'$, is hyperbolic in the direction $dx$ normal to the boundary
   and $dx$ [resp. $-dx$ ] is in the same component as $dt $; this holds if and only if 
 the spectrum of  $\underline A_d$ is positive [resp., negative].   
 
 \end{lem}

 \begin{proof}
By  \eqref{nn218}, $\underline \Delta$ is the characteristic polynomial of 
$\underline L'$, up to a nonvanishing factor. This implies $i)$. 
  
  Since  $\underline \Delta $ is hyperbolic in the time direction, 
 the eigenvalues of $\underline A_d$ are  real. 
  In addition,  the spectrum of $\underline A_d $ is positive if and only if 
  the matrix $\alpha \Id + (1- \alpha) \underline A_d $ is invertible 
  for all $\alpha \in [0, 1]$, thus, by definition and convexity of $\underline \Gamma_+ $,  
   if and only if  $dx $ belongs to this cone.  In the totally outgoing case, the proof  is similar. 
\end{proof}

 When $-\utau$ is semi-simple and nonglancing, there is a simple link between the first order 
 Taylor expansion of the reduced 
 symbols $\tau \Id + A^\flat(p, \eta, \xi)$ and 
 $ G^\flat(p, \zeta) + \xi \Id$ defined at \eqref{n28bb} and \eqref{nn35} respectively. 
 Similarly  to \eqref{nn216}, we introduce 
 $\underline G' (\zeta)$  the first order variation of 
 $G^\flat$ at the base point: 
 \begin{equation}
\label{nn313}
G^\flat(\up, \uzeta + \zeta) = G^\flat(\up, \uzeta) + 
\underline G' (\zeta) + O( \vert \zeta \vert^2). 
\end{equation}

 \begin{prop}
\label{prop34}
 Suppose that $-\utau$ is a semi-simple eigenvalue of $A(\up, \ueta, \uxi)$ of multiplicity 
 $m$ and 
$( \up, \utau, \ueta, \uxi)$  is nonglancing.  Then 
$- \uxi$ is a semi-simple eigenvalue of $G(\up, \uzeta)$  of multiplicity 
$m  $.
 Moreover, one can choose bases such that 
 \begin{equation}
\label{nn314}
\underline G'  (\zeta) + \xi \Id = (\underline A_d)^{-1} 
\Big( ( \tau - i \gamma) \Id + \underline A'  (\eta, \xi) \Big). 
\end{equation}

 \end{prop}

 \begin{proof}
{\bf a) }   We fix $p = \up$ and we omit the parameter $\up$ in the notations below. 
  There holds
 \begin{equation}
 \label{310}
 \Delta ( \utau + \tau, \ueta + \eta, \uxi + \xi ) = 
 e (\tau, \eta, \xi) \Delta' (\tau, \eta, \xi)
 \end{equation}
 with $e ( \utau, \ueta, \uxi) \ne 0$  and 
 $$
 \Delta'(\tau, \eta, \xi) = \tau^m + \sum_{j= 1}^m a_{m-j}(\eta, \xi)  \tau^j , \qquad
 a_{m-j} (\eta, \xi) = O\big( \vert \eta, \xi \vert^{m-j} \big)\,. 
 $$ 
Thus, 
\begin{equation}
\label{310b}
\underline  \Delta (\tau, \eta, \xi) = e(0,0,0) 
 \Big( \tau^m + \sum_{j= 1}^m   a^{(m-j)}_{m-j}(\eta, \xi)  \tau^j \Big) 
 \end{equation}
 where $a^{(m-j)}_{m-j}(\eta, \tau)$ denotes the homogeneous term of degree 
 $m-j$ in the Taylor expansion of $a_{m-j}$. 
 The nonglancing condition implies that  
 \begin{equation}
 \label{310c}
 a^{(m)}_m (0, \xi) = \underline a_{m} \xi^m  \quad \mathrm{with} 
 \quad \underline a_m \ne 0. 
 \end{equation}
 Thus, 
 $$
 \Delta(\utau, \ueta, \uxi+ \xi) = e_0(0,0,0) \underline a_m^{(m)} \xi^m 
 + O\big(\xi^{m+1}\big)
 $$
 showing that 
 $\uxi $  is a root of exact  order  $m$ of  $\Delta(\utau, \ueta,  \cdot ) = 0$.  
 Thus, $-\uxi$ is an eigenvalue of $G(\uzeta)$ of algebraic multiplicity $m$.

 Since $- \utau $ is semi-simple, the dimension of  $\ker (\utau \Id + A(\ueta, \xi) $ 
 is equal to $m$. 
 Since $\utau \Id + A(\ueta, \xi) = A_d^{-1} (\uxi \Id + G(\uzeta) $, 
 this space is equal to the kernel of   $(\uxi \Id + G(\uzeta)$
 showing that the geometric multiplicity of $-\uxi$ is  also equal to $m$. 
 Thus,    $-\uxi$ is   semi-simple.

 \medbreak
 {\bf b) }    Introduce the splittings 
\begin{equation}
\label{319}
\CC^N = \EE ( \zeta) \oplus \FF  ( \zeta) , \quad 
\CC^N = \EE_b (  \zeta) \oplus \FF_b (  \zeta)
\end{equation}
where $\EE  $ and $\FF $ [resp. $\EE_b$ and $\FF_b$] 
denote  the invariant spaces of $A$ [resp. $G$ ] associated 
to eigenvalues near and away from $-\utau$ [resp.  $- \xi$]. 
 We  have   shown  that 
 $\EE (\ueta, \uxi) = \EE_b (\uzeta)$. 
 Thus, one can choose bases such that the matrices $V$ and $V_b$ occurring 
 in \eqref{n28} and \eqref{nn34} respectively, satisfy
 \begin{equation}
\label{nn318}
\underline V = V(\ueta, \uxi) = V_b (\uzeta) . 
\end{equation}
The matrix   $W$ [resp. $W_b$ ]   vanishes  on $\FF$  [resp. $W_b$] 
and is equal to the inverse of $V$ [resp. $V_b $  on 
$\EE  $ [resp. $\EE_b$].  Since $-\utau$ and $-\uxi$ are semi-simple, 
\begin{equation}
\label{320}
\begin{aligned}
\FF_b ( \uzeta)  & = \mathrm{Range}  \Big(G( \uzeta) + \uxi \Id \Big) 
\\ 
& =  A^{-1}_d\mathrm{Range}  \Big(\utau \Id + A( \ueta,  \uxi )  \Big)
=  A^{-1}_d  \FF(  \ueta, \uxi) 
\end{aligned}
\end{equation}
  Since $\underline A_d = \underline W A_d(\up) \underline V$ is invertible
by Lemma \ref{lem35}, this implies that 
\begin{equation}
\label{322}
\underline W_b := W_b(\uzeta) = (\underline A_d)^{-1} \underline W A_d 
\quad \mathrm{and} \quad  \underline W_b \underline V = \Id . 
\end{equation}
Indeed, the matrix in the right hand side vanishes on $\FF_b(\uzeta)$ and 
multiplying it  by $\underline V$ on the right  gives $\Id$.   
 
 Differentiating \eqref{nn34} \eqref{nn35} and using that $G(\uzeta) = - \uxi \Id$, implies that 
\begin{equation}
\label{323}
\underline G' (\zeta) = \underline W_b G(\up, \zeta) \underline V . 
\end{equation}
With \eqref{nn32}, \eqref{n211} and \eqref{322} the identity \eqref{nn314} follows. 
\end{proof}


\subsection{Linearly splitting modes}

We now study the structure of $G^\flat$ when $(\up, \utau, \ueta, \uxi)$ is a real   
root of $\Delta$, of multiplicity $m \ge 2$, linearly splitting transversally to $\cM$ and nonglancing. 
We use the notations of Definition \ref{def213}. 
 We first assume that 
 \begin{equation}
\label{nn326}
dx \notin T_{\ueta, \uxi} \cM_{\up}. 
\end{equation}
With \eqref{n211} and \eqref{nn228e}, this is equivalent to the
condition that $\underline A_d$ has distinct eigenvalues, or, since $m \ge 2$,  
   \begin{equation}
\label{n43}
\underline A_d \notin \RR \Id. 
\end{equation}
In this case, there 
  are coordinates $(\eta_1, \ldots, \eta_{d-1}) $ in $\RR^{d-1}$ such that 
  $\cM$ is parametrized by $q = (p, \eta')$, with $\eta' = (\eta_1, \ldots, \eta_{d- \nu} )$,  close to $\uq = (\up, \ueta')$: 
  $\cM$ is locally given by 
  \begin{equation}
\label{n44}
  \eta''     = \phi (q ),  
\quad 
\xi    =  - \mu (q ) , 
\end{equation}
with $\eta'' = (\eta_{d- \nu +1} , \ldots, \eta_{d-1})$ and $\phi$, $\mu$ smooth. 
{}From now on, we work in such  coordinates. 
  With some abuse of notations, we write
  \begin{equation}
\label{n45}
\lambda(q)  := \lambda \big(p, \eta', \phi(q)  ,  
- \mu(q) \big).  
\end{equation}

Introduce the following notations  
    \begin{equation}
\label{n47}
\begin{array}{lll}
q = (p, \eta') ,  & 
\tilde \eta   = \eta''- \phi  (q), &  \tilde \xi = \xi + \mu(q),
\\
\tilde \tau = \tau +\lambda(q) ,  & 
\tilde \gamma = \gamma, &  \tilde \zeta = (\tilde \tau, \tilde \eta,  \tilde \gamma) ,  
\end{array}
\end{equation}
By  assumption,   $\lambda(q)$  is a semi-simple  eigenvalue   of $A\big(p, \eta',\phi(q)  ,  
- \mu(q)\big)$.  Therefore, 
\begin{equation}
\label{nn333}
A^\flat(p, \eta, \xi)    =  \psi(q) \Id + \tilde A (q, \tilde \eta, \tilde \xi),  
\quad \mathrm{with} \quad  \tilde A(q, 0, 0) = 0.  
\end{equation}
Moreover, Proposition~\ref{prop34} implies that 
   $\mu(q) $ is a semi-simple eigenvalue with  multiplicity $m$ of 
  $G \big( p, \lambda(q), \eta',   \phi(q') ,  0) \big)$, hence
  \begin{equation}
\label{nn334}
G^\flat(p, \zeta)   =  \mu (q) \Id  + \tilde G (q, \tilde \zeta) 
\quad \mathrm{with} \quad  
\tilde G(q, 0) = 0.    
\end{equation}
To blow up the singularity near $\cM$, we introduce  the matrices 
$\ccG(q, \rho, \cz) $ and $\ccA(q, \rho, \ct, \cx)$ such that 
for $\rho \cz$ and $\rho \cx$ small enough: 
\begin{equation}
\label{nn335}
\tilde G (q, \rho \cz  )  = \rho  \ccG (q, \rho, \cz ), 
\quad 
 \tilde A (q, \rho \ch , \rho \cx )  = \rho \ccA( q, \rho, \ch, \cx) . 
\end{equation}

\begin{prop}
\label{prop310} 
 Suppose that $(\up, \utau, \ueta, \uxi)$ is a real root of $\Delta$, of multiplicity 
 $m$, linearly splitting transversally to the manifold $\cM$, nonglancing and satisfying
$ \eqref{nn326}$. Then, with the notations as above, 
for $q$ in a neighborhood of $\uq$ and $\rho \cz$ and $\rho \cx$ small, 
\begin{equation}
\label{nn336}
  \det \big( \cx  \Id + \ccG(q,\rho, \cz )  \big) 
 = e(q, \rho \cz,  \rho\cx)  \det \big( (  \ct - i \cg ) \Id + \ccA(q, \rho, \ch, \cx)  \big) , 
\end{equation}
\begin{equation}
\label{nn337}
 \ker \big(  (\ct  u - i \cg ) \Id + \ccA(q,\rho , \ch, \cx )  \big) 
=  E(q, \rho \cz, \rho \cx )  \ker \big( \cx  \Id + \ccG(q, \rho, \cz )  \big) ,  
 \end{equation}
with $e (\uq, 0,0) \ne 0$, and 
 $E(q, 0, 0) = Id  $. Moreover,  
 the polynomial in $\xi$,  $\det \big( \cx  \Id + \ccG(q,\rho, \cz )  \big)$, 
 has real coefficients when $\cg = 0$, 
 and  $ \ccA(q, \rho, \ch, \cx)$ has only simple real eigenvalues when
 $\cx $ is real and $(\ch, \cx) \ne 0$. 
\end{prop}

\begin{proof}
The identities \eqref{nn336} and \eqref{nn337} follow directly from Lemma~\ref{lem32}
and \eqref{322}, through the change of variables \eqref{n47}.  
Similarly, that  $\det \big( \cx  \Id + \ccG(q,\rho, \cz )  \big)$ is real when $\cg = 0$ follows 
from Lemma~\ref{lem32}.  

Since the mode is  linearly splitting, the condition \eqref{nn228e} implies that 
$ \ccA(\uq, 0, \ch, \cx) $ has only real and simple eigenvalues when
$\vert (\ch , \cx)\vert = 1$ and $\cx$ is real. By continuity this extends 
to $ \ccA(q, \rho, \ch, \cx) $ for $q$ close to $\uq$ and $\rho$ small and $\vert (\ch , \cx)\vert = 1$.  
Since 
$$
\ccA(q, \rho, \ch, \cx)  = \frac{1}{\rho} \tilde A (q, \rho \ch, \rho \cx) = 
\vert (\ch, \cx) \vert  \ccA(q, \rho \vert (\ch, \cx) \vert  ,  (\ch, \cx) / \vert (\ch, \cx) \vert ), 
$$
$ \ccA(q, \rho, \ch, \cx)$ has only simple real eigenvalues when
 $\cx $ is real and $\rho \vert (\ch, \cx) \vert $ is small enough. 
      \end{proof}

  Next, we consider the case where 
  \begin{equation}
\label{nn338}
dx \in T_{\ueta, \uxi} \cM_{\up}. 
\end{equation}
Equivalently, this means that \begin{equation}
\label{nn339}
\underline A_d = a \Id. 
\end{equation}
  The nonglancing hypothesis implies that $a \ne 0$. In particular, the mode is totally incoming 
  when $a > 0$ and totally outgoing when $a < 0$.  
  
  As in Remark~\ref{rem213}, introduce the manifold $\widetilde \cM$
  of  points $(p, \tau, \eta, \xi)$ with $(p, \eta, \xi) \in \cM$ and 
  $\tau = \lambda(p, \eta, \xi)$. Since 
  $\underline L' (\dot \tau, \dot \eta, \dot \xi)$ vanishes when 
  $ (\dot \tau, \dot \eta, \dot \xi) \in T_{\utau, \ueta, \uxi} \widetilde \cM_{\up}$, the
  nonglancing hypothesis implies that 
  $\underline A_d = \underline L' (dx) \ne 0$, thus 
    \begin{equation}
\label{nn338b}
dx \notin T_{\utau, \ueta, \uxi} \widetilde \cM_{\up}. 
\end{equation}
  where $dx$ is considered in the now considered as the space-time conormal 
  to the boundary. However, \eqref{nn339} implies that 
  $\underline L'( dt + a dx) = 0$. By  \eqref{nn214e}, if 
  $\theta$ is not tangent to $\widetilde \cM_{\up}$, the dimension of k 
    $\ker \underline L'(\theta)$ is at most 1, thus 
  $ dt + a dx$ must be tangent to $\widetilde \cM_{\up}$. 
  This shows that, in contrast with the previous case, we cannot take
    $\tau$ and $\xi$ as independent variables transversal to $\widetilde  \cM$. 
  However, \eqref{nn338b} implies that there are coordinates  
  $(\eta_1, \ldots, \eta_{d-1}) $ in $\RR^{d-1}$ such that 
  $\cM$ is parametrized by $(q, \tau)$, with $q = (p,  \eta_1, \ldots, \eta_{d- \nu-1} )$,  
   and  locally given by 
  \begin{equation}
\label{n44b}
  \eta''     = \phi (q, \tau  ),  
\quad 
\xi    =  - \mu (q, \tau ) , 
\end{equation}
with $\eta'' = (\eta_{d- \nu +1} , \ldots, \eta_{d-1})$ and $\phi$, $\mu$ smooth, 
analytic in $\tau$. Denote by $\cM_b$ the manifold 
$\{ \eta '' = \phi (q, \tau) \}$.  

  Proposition~\ref{prop34} implies that $- \xi$ is a semi-simple eigenvalue with multiplicity 
  $m$ of $G(p, \tau, \eta, 0)$ if $(p, \tau, \eta, \xi) \in \widetilde \cM$. 
Thus, 
\begin{equation}
\label{nn342}
G^\flat(p, \tau, \eta, 0)   =  \mu (q, \tau   ) \Id  + \tilde G (q, \tau , \tilde \eta) 
\quad \mathrm{with} \quad      \tilde \eta   = \eta''- \phi  (q, \tau)     
\end{equation}
and  $ \tilde G(q, \tau , 0) = 0$.    
Thus, there are  matrices 
$\ccG(q, \rho, \tau, \ch) $   such that 
for $\rho \ch$   small enough: 
\begin{equation}
\label{nn343}
\tilde G (q, \tau, \rho \ch  )  = \rho  \ccG (q, \rho, \tau, \ch ). 
\end{equation}

\begin{prop}
\label{prop311}
 Suppose that $(\up, \utau, \ueta, \uxi)$ is a real root of $\Delta$, of multiplicity 
 $m$, linearly splitting transversally to the manifold $\cM$, nonglancing and satisfying
$ \eqref{nn326}$. Then, with the notations as above, 
for $(q, \tau) $ in a neighborhood of $(\uq, \utau)$ and $\rho \ch$  small, 
$\ccG(q, \rho , \tau, \ch)$ has only real and simple eigenvalues when $\ch \ne 0$. 
Moreover, 
$\D_\gamma G^\flat (\up, \uzeta) =  -  i a^{-1}  \Id $ 

\end{prop}
  
\begin{proof}
By \eqref{nn339}, we are in a   totally nonglancing  framework,  implying that 
$\cx + \underline G' (\up, \tau, \eta, 0)$ is  hyperbolic in the direction $dx$. 
Taking $(\tilde \eta, \xi)$ as variables transversal to $\widetilde \cM$, 
\eqref{nn214},   implies  that 
$\cx + \underline G' (\up, \utau, \ueta', \ch, 0)
= \cx + \ccG (\uq, 0, \utau, \ch) $ is strictly hyperbolic in the direction 
$dx$. This shows that $ \ccG (\uq, 0, \utau, \ch) $ has only real and simple eigenvalues 
when $\ch \ne 0$. This is preserved for $(q, \tau)$ close to $(\uq, \utau)$, 
$\ch \ne0$ and $\rho$ such that $\vert \rho \ch \vert$ is small enough. 

By \eqref{nn314}, 
$$
\D_\gamma G^\flat (\up, \uzeta) = - i \D_\tau G^\flat(\up, \uzeta) = 
-  i  \underline A_d^{-1}  =  - i a^{-1} \Id. 
$$
\end{proof}

   
 
 \section{The Lopatinski condition and symmetrizers}

 \subsection{Maximal estimates}
 
 Consider a system $L(p, \D_t, \D_y, \D_x) = A_d (p) 
 ( \D_x+ G(p, \D_t, \D_y))  $, hyperbolic in the direction $dt$ such  that 
 the boundary $\{ x = 0 \}$ is not characteristic. 
 The classical plane wave analysis yields the equations
 \begin{equation}
 \label{eq51}
 \D_x u +  i G (p, \zeta ) u = f, 
 \end{equation} 
 with $\zeta = (\tau - i \gamma, \eta)$, identified with $(\tau, \eta, \gamma)$. 
 We supplement the operator $L$ with  boundary conditions, which read, 
 after Laplace-Fourier transform, 
 \begin{equation}
 \label{eq52}
  M (p, \zeta) u _{\vert x = 0 } = g  . 
 \end{equation}
  We   assume that $M$ depends smoothly on 
  $p$ and $\zeta \in \RR^{d+1}\backslash\{0 \} $ and is homogeneous of degree 
  $0$ in $\zeta$. We allow $M$ to depend on $\zeta$ for two reasons: 
  first, it is so in the shock problem and second,   
  the block reduction process leads us to consider such cases. 
   By homogeneity, we can restrict attention to 
 $\zeta \in S^d_+ = \{ \vert \zeta \vert = 1 ; \gamma > 0 \}$.

  Recall that the hyperbolicity condition implies that for 
$\zeta \in S^d_+$, the eigenvalues of $G(p, \zeta)$ are nonreal. 
The $L^2$ solutions of   $ (\D_x + i G )u = 0 $
consist of all $u(x) = e^{ - i x G} h$ with $h \in \EE_-(p, \zeta)$,
where $\EE_-(p, \zeta)$ denotes
the invariant space of $G(p, \zeta)$ 
associated to the spectrum in $\{ \im \mu <  0 \}$.  Recall that  
\begin{equation}
\label{eq56}
\dim \EE_-(p, \zeta) = N_+  
\end{equation}
 is the number of positive eigenvalues of $A_d$. 
  We   assume that we have the correct number of boundary conditions, that is 
\begin{equation}
\label{eq53}
  \dim \ker M(p, \zeta) = N_- = N- N_+ ,  
  \end{equation}
    the number of negative eigenvalues of $A_d$.  
 For simplicity, we also assume that $M$ is surjective, that is that 
 $M$ is an $N_+ \times N$ matrix. 
 
 \smallbreak
 Following Kreiss and Majda, we are interested in the maximal estimates for solutions 
 of \eqref{eq51} \eqref{eq52}:
 \begin{equation}
\label{eq54}
\gamma \Vert u \Vert^2_{L^2(\RR_+) } + 
\vert u(0) \vert^2 \ls \frac{1}{\gamma} \Vert f \Vert^2_{L^2(\RR_+) } + 
\vert g  \vert^2, 
\end{equation}
 where $\ls$ means that the left hand side is smaller than $C$ times the right hand side
 where $C$ is independent of $u, f,g$ and $\zeta \in S^d_+$.  
A necessary condition is that
\begin{equation*}
\vert u(0) \vert \ls \vert M (p, \zeta) u(0) \vert
\end{equation*}
for all $u \in L^2$ such that  $\D_x u + i G u = 0$.
 Thus, a necessary condition for  \eqref{eq54} is that 
for all $\zeta \in S^d_+$  there holds 
\begin{equation}
\label{eq55}
\forall h \in \EE_-(p, \zeta) , \quad 
\vert  h \vert \le  C \vert M (p, \zeta) h \vert .   
\end{equation}

\begin{defi}\label{lop}
The  Lopatinski determinant  associated with $L$, $M$
is defined for $gamma >0$ as
  \begin{equation}
  \label{eq57}
  D(p, \zeta) = \det \big( \EE_-(p, \zeta) , \ker M(p, \zeta) \big).
  \end{equation}
We say that $L$, $M$ satisfy the
 uniform Lopatinski condition on   $\omega$    if there is a constant
 $c > 0$  such that 
 \begin{equation}
 \label{eq58}
 \forall (p, \zeta) \in \omega \times S^d_+  : \qquad \vert D(p, \zeta) \vert \ge c 
 \end{equation}
\end{defi}

In \eqref{eq57}, the determinant is obtained by taking orthonormal bases in 
each space $\EE_-$ and $\ker M$ (note that 
$\dim\EE_- + \dim \ker M = N$) and  independent of this choice. 
It depends on the scalar product used in $\CC^N$ to form orthonormal bases, but the condition \eqref{eq58} is independent of this scalar product. 
Denoting by $\pi_{M^\perp}$ the orthogonal projection 
on $\ker M^\perp$, then 
$D(p, \zeta)$ is the determinant of the 
restriction of $ \pi_{M^\perp} $ to $   \EE_-$ (in orthonormal bases).  
Thus, using that $M$ is an isomorphism from 
$\ker M^\perp$ to $\CC^{N_+}$, uniformly bounded 
with uniformly bounded inverse
for $\zeta \in \overline S^d_+$, we see that the condition 
\begin{equation}
\label{eq59}
\vert D (p, \zeta ) \vert \ge c,
\end{equation}
implies    \eqref{eq55}  with $C$ depending only on $c$. Conversely, if 
\eqref{eq55} holds, then \eqref{eq59} is satisfied for some $c > 0$ which depends only on 
$C$. 
 
 This shows that the uniform Lopatinski condition is necessary for \eqref{eq54}, see \cite{Kr}. 
  Conversely,  there holds: 
  
  \begin{theo}[First Main Theorem]
  \label{theo52}
  Suppose that  $L$ is symmetric and that 
   the roots of the characteristic equations are either geometrically regular or  
  totally nonglancing, or nonglancing and linearly splitting transversally to 
  smooth manifolds. 
  Then,  the uniform Lopatinski condition at 
  $\up$ implies the maximal estimate $\eqref{eq54}$ for $p$
  in a neighborhood of $\up$. 
  \end{theo}

 This was established by Kreiss (\cite{Kr}; see also \cite{Ral, Ch-P}) 
in the strictly hyperbolic case by an ingenious construction of   symmetrizers.
It was extended to systems with block structure by Majda and Osher
\cite{Maj, MaOs}.
By Theorem \ref{theo31}  geometric regularity implies
block structure. 
In the totally nonglancing case, we extend Kreiss's construction of smooth symmetrizers
using the symmetry of $L$ (Theorem~\ref{2main}). 
In the presence of  nonglancing linearly splitting modes, 
we extend further the construction
of smooth symmetrizers when there is at most    one linearly splitting mode
(for each $ \zeta$) corresponding to a double root (Theorem~\ref{4main}).  
In the general case  of 
nonglancing linearly splitting modes, we prove Theorem~\ref{theo52} using nonsmooth symmetrizers.  
 For applications to nonconstant coefficients and nonlinear equations, the 
 use of nonsmooth symmetrizers seems very delicate. 
 In this respect, Theorems~\ref{2main} and  \ref{4main} below are 
more important than Theorem~\ref{theo52}.
  
    \begin{rem}
    \label{rem43x}
    \textup{ When the vector bundle $\EE_-(p, \zeta)$ has a continuous extension 
    to  $\omega \times \overline S^d_+$, for some  neighborhood $\omega$ of $\up$, 
    then $D$ extends as a continuous function on $\omega \times \overline S^d_+$,
    and the uniform Lopatinski condition holds, on a possibly smaller neigborhood
    of $\up$, if and only if 
    \begin{equation}
    \label{eq410x}
    \forall \zeta \in \overline S^d_+, \quad  D(\up, \zeta) \ne 0, 
    \end{equation}
    or, equivalently, if and only if 
      \begin{equation}
    \label{eq411x}
    \forall \zeta \in \overline S^d_+, \quad  \EE_-(\up, \zeta) \cap \ker M(\up, \zeta) = \{ 0 \}. 
    \end{equation}
    The continuous extendability of $\EE_-$ to boundary values $\gamma = 0$, 
    is true for strictly hyperbolic systems
    (\cite{Kr}), or when the block structure condition is satisfied (\cite{Maj}), or,  more generally,
    when the 
    nongeometrically regular modes are totally nonglancing (see Theorem~\ref{2main}). 
    In the other cases, one can introduce the set of limits of sequences of vectors  
    $\EE_-(p_n, \zeta_n)$ as $(p_n, \zeta_n)$ tends to $(p, \zeta)$ and $\gamma_n > 0$:  
    \begin{equation}
    \label{eq412x}
    \hat \EE_{-} (p, \zeta) = \limsup_{\substack{(p', \zeta') \to (p, \zeta)\\ \gamma' > 0}}
    \EE_-(p'; \zeta') . 
    \end{equation}
    This is a closed cone in $\CC^N$, equal to $\EE_-(p, \zeta)$ if 
    $\gamma > 0$ and equal to the continuous extension of $\EE_-$  when 
    such an extension exists. If the uniform estimate \eqref{eq55} holds 
    for $p$ close to $\up$ and $\zeta \in S^d_+$, then 
    \begin{equation}
    \label{eq413x}
  \forall \zeta \in \overline S^d_+, \   \forall h \in \hat \EE_-(\up, \zeta) , \quad  \vert h \vert \le C \vert M(\up, \zeta) h \vert 
    \end{equation}
    Conversely, this implies \eqref{eq55} on a neighborhood of $\up$, with a possibly larger 
    constant $C$. By homogeneity, it is sufficient to check \eqref{eq413x} for 
    $h$ in the unit sphere, and by compactness this shows that this condition is equivalent to
    the following analogue of \eqref{eq411x}: 
      \begin{equation}
    \label{eq411xx}
    \forall \zeta \in \overline S^d_+, \quad  \hat \EE_-(\up, \zeta) \cap \ker M(\up, \zeta) = \{ 0 \}. 
    \end{equation}
    }
    \end{rem}

 
 \subsection{Symmetrizers}

 We recall first the method of symmetrizers. 
 We fix a neighborhood $\omega$ of $\up$. 

\begin{defi}
\label{def53}
A   bounded    symmetrizer  for  $G(p, \zeta) $ on 
$ \Omega =  \omega \times U $, $U \subset  S^d_+$, 
is a smooth matrix $\Sigma (p, \zeta)$ 
on  $ \Omega$, 
such that, for some  $C$, $c>0$, there holds for 
all $(p, \zeta) \in \Omega $, 
 \begin{eqnarray}
\label{global41}
&&\Sigma (p, \zeta) = \Sigma ^*(p, \zeta)  ,   
\\
\label{global42}
&& \vert \Sigma (p, \zeta) \vert \le C,    
\\
\label{global43}
&& \im \Sigma (p, \zeta) G(p, \zeta)  \ge c \gamma \Id,  
\end{eqnarray}
It is a Kreiss symmetrizer for $(G, M)$ if in addition 
\begin{equation}
\label{global44}
 \Sigma (\up, \uzeta) + CM^*(p, \zeta)  M(p, \zeta)  \ge c \Id.
\end{equation}
The symmetrizer is smooth, if it extends smoothly to $\omega \times \overline U 
\subset    \omega \times \overline S^d_+$. 
\end{defi}

Taking the scalar product of the equation \eqref{eq51} with $\Sigma u$ and integrating by 
parts, immediately implies the following.

\begin{lem}
\label{lem54}
If there exists a  bounded  Kreiss symmetrizer  for $(L, M)$ on 
$\omega \times S^d_+$ 
then the maximal estimate \eqref{eq54} holds for all 
$(p, \zeta) \in \omega \times S^d_+$. 
\end{lem}

\begin{rem}
\label{rem55}
\textup{\eqref{global44} implies that
\begin{equation}
\label{eq514}
\Sigma(p, \zeta)   \ge c \Id  \quad   \mathrm{ on } \  \ker M (p, \zeta).
\end{equation} 
Conversely, if this inequality holds, then 
\begin{equation*}
( \Sigma h , h) +  C_1 \vert \pi_{\ker M^{\perp} } h \vert^2  \ge \frac{c}{2} \vert h \vert^2 
\end{equation*} 
if $C_1 \ge C + c/2 + 2 C^2/ c$  and $C \ge \vert \Sigma \vert$ as in \eqref{global42}. 
With $C_2$ such that 
\begin{equation}
\label{eq515}
\vert \pi_{\ker M^{\perp} } h \vert \le C_2 \vert M h \vert, 
\end{equation} 
this shows that \eqref{eq514} implies \eqref{global44} with $c$ replaced by $c/2$ 
and $C $ by $ C_1 C_2^2$. This shows that one can replace the condition \eqref{global44} by \eqref{eq514}. }
\end{rem}

\begin{rem}
\label{lem56}
\textup{ If $S(p) $ is a symmetrizer for $L$,  that is 
$S(p)= S^*(p)$ is positive definite and all the $S(p)A_j(p)$ are self adjoint, then 
$\Sigma (p) = - S (p)A_d(p)$ is a symmetrizer for $G(p, \zeta)$ :  the properties \eqref{global41} \eqref{global42} 
and \eqref{global43} are satisfied. When the  property \eqref{global44} or equivalently 
\eqref{eq514} holds, the boundary conditions are said to be 
\textit{strictly dissipative}. In this case, the maximal estimates \eqref{eq54}
are satisfied. For the theory of dissipative boundary value problems, 
we refer to  \cite{F.1, F.2,  FL, Rau}. 
}
\end{rem} 
    
    We recall the following general fact about symmetrizers. 
    \begin{lem}\label{sign}
Let $\Sigma$ be a symmetrizer for $G$.  Then, for $\gamma > 0$,
$\Sigma(p, \zeta) $ is negative definite on the stable subspace $\EE_-(p, \zeta) $ of $G(p, \zeta)$. 
\end{lem}

\begin{proof}
Let $v$ satisfy   $\D_x v +  i  Gv= 0 $, $v(0)=u \in \EE_- $, $u \ne 0$.
Then, $v \to 0$ exponentially fast  as $x \to +\infty$.  On the other hand, using
symmetry of $\Sigma $, we have 
\begin{equation}
\begin{aligned}
(d/dx )\langle \Sigma v, v \rangle&= 
2 \re  \langle \Sigma  \D_x v, v \rangle 
= - 2 \re \langle i \Sigma  G v,  v\rangle \\
& = 2  \langle (\im \Sigma G) v , v \rangle \ge  c \gamma  |v|^2 >0.
\end{aligned}
\end{equation}
Since $\langle v,\Sigma v\rangle \to 0  $ and  is integrable  as $x \to+\infty$, we
must have, therefore 
\begin{equation}
\langle u,\Sigma u\rangle \le -  c \gamma \Vert v \Vert^2_{L^2(\RR_+)}  < 0 . 
\end{equation}
as claimed.
\end{proof}

   In the analysis of  \cite{Kr, Maj}, there is an intermediate step in the construction of   symmetrizers: 
    one construct first a family of symmetrizers $\Sigma^\kappa$, which satisfy 
    \eqref{global41} \eqref{global42} and \eqref{global43}, 
     such that the negative cone of $S^\kappa$ is  an arbitrarily  
    small conic neighborhood of $\EE_-$. 
    Next, one uses the uniform Lopatinski condition : because $\ker M$
    does not intersect $\EE_-$, it is contained in the positive cone of
    $\Sigma^\kappa$ for f $\kappa$  large enough
implying  
    \eqref{global44}.    
    
    \begin{defi}
    \label{def57}
Consider a family $\Sigma^\kappa$  of bounded [resp. smooth ] symmetrizers for $G$   on 
$\Omega^\kappa = \omega^\kappa \times U^\kappa$. 
It is a K-family of bounded [resp. smooth] symmetrizers,    if  
 there are $C$ and  projectors 
$\Pi_-(p, \zeta)  $ on $\EE_-(p, \zeta)$ such that 
   \begin{eqnarray}
 \label{517b}
 &&\forall \kappa \ge 1, \  \forall (p, \zeta) \in \Omega^\kappa, 
 \quad \vert \Pi_-(p, \zeta) \vert \le C , 
  \\
  \label{518b}
&&  
\forall \kappa \ge 1, \ \forall (p, \zeta) \in \Omega^\kappa, 
\quad   \Sigma^\kappa  (\up, \uzeta) \ge  m(\kappa)  \Pi_+^*   \Pi_ +  - 
 \Pi_-^*   \Pi_-. 
\end{eqnarray}
  where $  \Pi_+ = \Id - \Pi_- $ and $m(\kappa) \to + \infty $ as $\kappa \to + \infty$. 
 
   \end{defi}
 
 \begin{prop}
 \label{prop58}
 Suppose that   $\Sigma^\kappa$ 
 is a  K-family of symmetrizers for $G$   on $\Omega^\kappa$.  
 Suppose that $M$ satisfies  the uniform Lopatinski condition. Then, 
  for $\kappa$ large enough, $  \Sigma^\kappa$ is a Kreiss symmetrizer 
 on $\Omega^\kappa$ for $(G, M)$.  
 \end{prop}
 
 \begin{proof}
The Lopatinski condition implies that there is a constant $C$ such that 
\eqref{eq55} holds. Therefore,
$$
\vert \Pi_- h \vert  \le C  \vert  M h \vert +  C \vert M \vert \,  \vert \Pi_+ h \vert. 
$$ 
Thus,   there are $C_1$ and $C_2$ such that 
$$ 
\vert h \vert^2 \le  2 \vert \Pi_- h \vert^2 +  2 \vert \Pi_+ h \vert^2 
\le C_1 \vert M h \vert^2 +  C_2  \vert \Pi_+ h \vert^2 -  \vert \Pi_- h \vert^2  . 
$$
Therefore, if $m(\kappa) \ge C_2$, \eqref{global44} follows. 
\end{proof}

\begin{rem}
\label{rem59}
\textup{The choice of the projector $\Pi_-$ is irrelevant, as long as the uniform bound
\eqref{517b} holds. 
If $\tilde \Pi_-$ is another projector on 
$\EE_-$ satisfying \eqref{517b} with a constant $\tilde C$, then $\tilde \Pi_+   \Pi_- =  0 $ and 
\begin{equation*}
\begin{aligned}
& \vert \tilde \Pi_+ h \vert = \vert \tilde \Pi_+ \Pi_+ h \vert  \le  \tilde C \vert   \Pi_+ h \vert , 
\\
&  \vert   \Pi_- h \vert \le   C  ( \vert  \tilde  \Pi_+ h \vert  +  \vert  \tilde  \Pi_- h \vert).
\end{aligned}
\end{equation*}
Thus, 
\begin{equation*}
m(\kappa) \vert \Pi_+ h \vert^2 - \vert \Pi_- h \vert^2 
\ge (  m(\kappa)/ \tilde C^2  - 2 C^2)  \vert \tilde \Pi_+ h \vert^2  - 2 C^2 \vert \tilde \Pi_- h \vert^2. 
\end{equation*}
Therefore, changing $\Sigma^\kappa$ to $C^{-2} \Sigma^\kappa$ 
 we see that 
\eqref{518b} for $\Sigma^\kappa$  and $\Pi_\pm$ implies  similar estimates
for $\tilde \Sigma^{\tilde \kappa}$  and $\tilde \Pi_\pm$, with 
$ \tilde m(\kappa) = m(\kappa) / C^2 \tilde C^2 - 2  $. 
In particular, we can always choose $\Pi_- (p, \zeta)$ to be the orthogonal projector
on $\EE_-(p, \zeta)$.  
 }
\end{rem}

 
 \subsection{Localization and block reduction}
 
 We collect here several remarks concerning the construction of symmetrizers. 
First, the construction is local. 

\begin{lem}
\label{lem510}
Suppose that for all $(\up, \uzeta ) \in \overline \omega \times \overline S^d_+$, there are neigborhoods
$\Omega^\kappa$ of $(\up, \uzeta)$ and self adjoint  matrices 
$\Sigma^\kappa(p, \zeta)$  for $(p, \zeta) \in \Omega^\kappa_+ = 
\Omega^\kappa \cap \{ \gamma > 0 \}$ such that 
the $\Sigma^\kappa$ form a $K$-family of symmetrizers for $G$ on 
 $\Omega^\kappa_+$. 
 Then there exists neigborhoods $\omega^\kappa$ of $\overline \omega $ and 
$K$- families of  symmetrizers $\tilde \Sigma^\kappa $ on 
 $\omega^\kappa \times S^d_+$. 
 
 If the local $\Sigma^\kappa$ are smooth, 
then the global  $\tilde \Sigma^\kappa$ can be chosen smooth. 
 
\end{lem}

\begin{proof}
 By Remark \ref{rem59}, we can assume that the local symmetrizers 
 $\Sigma^\kappa(p, \zeta)$  on $\Omega^\kappa_+$  satisfy \eqref{518b}  
 with $\Pi_-$ equal to the orthogonal projector on $\EE_-$. Relabeling the families, 
 we can also assume that they satisfy \eqref{518b} with  $m(\kappa) = \kappa$. 
 
For all $\kappa$, we can find 
a finite covering of $ \overline \omega \times \overline S^d_+$ by open sets 
 $\Omega^\kappa_j $, such that 
 a $K$-family of symmetrizers $\Sigma^\kappa_j$  exists on   $\Omega^\kappa_{j , +}$. 
 Consider a 
a partition of unity
$1 = \sum \chi^\kappa_j$ 
on $\overline \omega \times \overline S^d_+$, with $\chi^\kappa_j$ supported in 
$\Omega^\kappa_j$. Let
\begin{equation}
\label{ne419}
\Sigma^\kappa (p, \zeta) = \sum_j \chi_j(\zeta) \Sigma^\kappa_j (p, \zeta). 
\end{equation}
Because the covering is finite, we can choose uniform constants $C$ and $c$ 
in  \eqref{global41}, \eqref{global42}, \eqref{global42} and  \eqref{518b} and the lemma follows. 
\end{proof}

There is an analogous result for  Kreiss symmetrizers. 

\begin{lem}
\label{lem511}
Suppose that for all $(\up, \uzeta) \in \overline \omega \times \overline S^d_+$, there is a  neigborhood
$\Omega$ of $(\up, \uzeta)$ and a  bounded [resp. smooth] Kreiss symmetriser for $(G, M)$ on  
$ \Omega _+ =  \Omega^\kappa \cap \{ \gamma > 0 \}$.  Then there is a neigborhood $\omega' $ of $\overline \omega $ and a bounded [resp. smooth]  symmetrizer   on 
 $\omega' \times S^d_+$. 
 
\end{lem}

Next, we consider a smooth diagonal block reduction of $G$ on a  neighborhood
$\Omega$ of $(\up, \uzeta)$:  
\begin{equation}
\label{eq517}
U^{-1}(p, \zeta	) G( p, \zeta) U(p, \zeta) = \mathrm{block\, diag} (G_k(z, \zeta)). 
\end{equation} 
For instance, we  can consider the distinct eigenvalues $\underline \mu_k $, 
$k \in \{1, \ldots, \underline k \}$,  of 
$G (\up, \uzeta)$, small discs $D_k$ centered at   
$\underline \mu_k $ that do not intersect each other, and for 
$(p, \zeta)$ close to $(\up, \uzeta)$, the diagonal block reduction where 
the spectrum of $G_k$ is contained in $D_k$.  

Equivalently, taking appropriate blocks in $U$ and $U^{-1}$, one can 
introduce smooth matrices $V_k$ and $W_k$ such that 
\begin{equation}
\label{ne423}
G = \sum V_k G_k W_k , \quad   W_k    V_j  = \delta_{j,k} \Id.
\end{equation}

 For $\gamma > 0$, we denote by $\EE_{k, -}(p, \zeta)$ the invariant subspace of $G_k $  
  associated to eigenvalues in $\{ \im \mu < 0 \}$. Thus, 
 \begin{equation}
\label{eq516b}
\EE_-(p, \zeta)  = \bigoplus_k V_k(p, \zeta) \EE_{k,-} (p, \zeta). 
\end{equation}

Given  symmetrizers    $\Sigma _k$ for $G_k$ on 
$\Omega  $, 
\begin{equation}
\label{522c}
\Sigma  =  U^{-1\, *}  \mathrm{diag} (\Sigma _k ) U^{-1}  = 
\sum     W_k^* \Sigma_k  W_k  .
\end{equation} 
 is a symmetrizer for $G$. 
 
 \begin{lem}
 \label{lem512}
 Suppose that for all $k$, 
$ \Sigma^\kappa_k$  is a $K$-family of symmetrizers for $G_k$ on 
$\Omega^\kappa$, then there are  $K$-families  of symmetrizers 
$\Sigma^\kappa$
for $G$ on $\Omega^\kappa$.  

If the $\Sigma^\kappa_k$ are smooth, 
then    $ \Sigma^\kappa$ can be chosen smooth. 
 \end{lem}    
     
   \begin{proof}
   Consider the symmetrizers  $\Sigma^\kappa$ associated  by \eqref{522c} 
   to the $\Sigma^\kappa_k$. 
  By Remark \ref{rem59}, we can assume that the 
     $\Sigma_k^\kappa$ satisfy \eqref{518b} with $m(\kappa) = \kappa$ and 
     $\Pi_{k, -}$ equal to the orthogonal projection on 
     $\EE_{k, -}$. 
     Therefore,  for $h = \sum V_k  h_k   $,  there holds
     \begin{equation*}
( \Sigma^\kappa h, h ) = \sum (\Sigma^\kappa_k h_k, h_k) 
\ge \kappa \sum \vert \Pi_{k +} h_k \vert^2  -  \sum \vert \Pi_{k, -} h_k \vert^2 . 
\end{equation*}
Let  $\Pi_- = \sum V_k \Pi_{k, -} W_k $. It is a projector on $\EE_-$ and 
\begin{equation*}
\vert \Pi_+ h \vert^2 \le \vert U \vert^2  \sum   \vert \Pi_{k, +}  h_k \vert^2 , 
\quad 
 \sum   \vert \Pi_{k, -}  h_k \vert^2 \le  \vert U^{-1}  \vert^2 \vert \Pi_-  h \vert^2  . 
\end{equation*}
Therefore, $\tilde \Sigma^\kappa = \vert V^{-1} \vert^{-2} \Sigma^\kappa$
satisfies \eqref{518b} with 
$m(\kappa) = \kappa / \vert U \vert^2 \vert U^{-1} \vert^2$. 
   \end{proof}
     
     There are no simple analogue of Lemma~\ref{lem512} fro Kreiss symmetrizers, 
     since the boundary conditions $M$ 
     do not split, in general,  into boundary independent conditions    $M_k$ for 
     each block $G_k$. However, we give in Section 6 an interesting result 
     in this direction.

     We end this section with the following elementary remark. 
     
     \begin{rem}
     \label{rem413}
     \textup{ If $\Sigma $ is a symmetrizer for $G$, then 
     $\Sigma_k = V_k^* \Sigma V_k  $ is a symmetrizer for $G_k$, 
     since $\Sigma_k G_k =  V_k^* \Sigma G V_k$. 
     In particular, if the original system $L$ is hyperbolic symmetric, with symmetrizer
     $S(p)$, then $\Sigma (p) = - S(p) A_d(p)$ is a symmetrizer for 
     $G$ and $\Sigma_k (p, \zeta) = - V_k^*(p, \zeta) S(p) A_d(p) V_k(p, \zeta)$ 
     is a smooth symmetrizer for $G_k(p, \zeta)$. 
     }
     
     \end{rem}

     %
     %

     \section{Proof of the maximal estimates}
     
     In this section we prove Main Theorem~\ref{theo52}. Because of its special interest
     for MHD, we also point out in   Theorem~\ref{2main} below, the particular case   where 
     only geometrically regular and totally nonglancing modes are present. 
     In this case, smooth Kreiss symmetrizers are available, while they are not in the more general
    situation of Theorem~\ref{theo52}. The construction of smooth symmetrizers is studied in more details in the next section.

We fix $\uzeta \in \overline S^d_+$. With 
$\underline \mu_k$  denoting the    distinct eigenvalues of $G(\up, \uzeta)$,
we consider on a neighborhood $\Omega$
of $(\up, \uzeta)$, the  block reduction 
 \eqref{eq517} of $G$  such that the spectrum of $G_k(\up, \uzeta)$ is 
 $\{ \underline \mu_k \}$. 
 In this section, give different  constructions of  K-families of symmetrizers, 
 depending on the properties of the blocks $G_k$. 
 Be denote by $N_k$ the dimension of the block $G_k$, that is the algebraic 
 multiplicity of $\underline \mu_k$. 
 When $\gamma > 0$,  $G_k(p, \zeta) $ has no real eigenvalue,   
 we denote by $\EE_{k,-}(p, \zeta) $ the invariant space of $G_k$ associated to eigenvalues 
 in $\{\im \mu  < 0 \}$  and, in accordance with  \eqref{eq56}, we set 
 \begin{equation}
\label{nn51}
\dim \EE_{k, -} (p, \zeta) = N_{k, +}  \quad   \gamma > 0. 
\end{equation}

   \subsection{Kreiss' Theorem}
     
     We first recall Kreiss' construction with Majda's extension. 
     
      \begin{prop}
 \label{prop513}
   Suppose that  $G_k$ satisfies the block structure condition at 
   $(\up, \uzeta)$. There are  
  spaces $\underline \EE_{k,+}$ and $\underline \EE_{k, -}$, neighborhoods $\Omega^\kappa $
 of $(\up, \uzeta)$  and a smooth 
family of symmetrizers  $\Sigma^\kappa_k (p, \zeta)$ on 
$\overline \Omega^\kappa _+ = \Omega^\kappa  \cap \{ \gamma \ge 0 \}$ such that: 
 \begin{eqnarray}
 \label{41}
 &&\CC^{m'_k}  = \underline \EE_{k,  -}  \oplus \underline \EE_{ k, +},    
 \qquad \dim \underline \EE_{k,-} = N_{k,+}, 
 \\
\label{44}
&& \Sigma^\kappa_k (\up, \uzeta) \ge \kappa \underline \Pi_{k, +}^* \underline \Pi_ {k, +}  - 
\underline \Pi_{k, -}^* \underline \Pi_{k, -}. 
\end{eqnarray}
  where $\underline \Pi_\pm$ is the projector on $\underline \EE_\pm$ 
  in the decomposition $\eqref{41}$. 
 \end{prop}

The block structure condition covers two cases: 

\quad  1) \textit{Elliptic modes}.  If $\mu_k$ is not real, then  
  $G_k$  satisfies the block structure condition. 
 Then $\EE_{k, -} = \EE_k $  [resp. $\EE_{k, -} = \{0 \}$]  if 
 $\im \mu_k < 0$  [resp.  $\im \mu_k > 0$]; 
 There is a symmetrizer $\Sigma $ such that 
 \begin{equation}
 \label{nn55}
 \im \Sigma_k G_k   \ge c > 0 
 \end{equation} 
 and $\Sigma$ is definite negative [resp. positive] when 
 $\im \mu_k < 0$  [resp.  $\im \mu_k > 0$];  
 one chooses $\Sigma^\kappa = \Sigma$
 [resp. $\Sigma^\kappa = \kappa \Sigma$]. 
 
  Elliptic modes are the only possibility  when $\uga \ne  0$. 
  
  \smallbreak
  
  \quad 2) \textit{Geometrically regular modes.}
 If $\uga= 0$ and 
  $(\up, \utau, \ueta, -\underline \mu_k)$  is geometrically regular 
   then, by Theorem~\ref{theo31},
   $G_k$ satisfies the block structure condition. The 
   symmetrizers $\Sigma^\kappa_k $  constructed in \cite{Kr} (see also \cite{Ch-P})
   satisfy   
   \begin{equation}
\label{nn54}
\im (\Sigma^\kappa_k G_k )  = \gamma E^\kappa_k, \quad   E_k \ge c^\kappa  \Id ,
\end{equation}
  which implies \eqref{global43}. This  improvement  is useful for applications
  to variable coefficients:  it allows one
 to use standard Garding's inequalities 
  for the operator  with symbol $E_k$ (see \cite{Ch-P, Mok, MetKochel, MeZ.1}). 
  
  \medbreak
  
  When the block structure condition is satisfied, the spaces 
  $\EE_{k,-}(p, \zeta) $ have limits as 
  $\gamma \to 0$ (see \cite{Kr, Ch-P}). It is trivial for elliptic modes. 
  For geometrically regular modes,  it is proved in \cite{MeZ.2} that the existence of symmetrizers 
  as in Proposition~\ref{prop513} implies that 
  \begin{equation}
\label{nn55one}
\underline \EE_- = 
\lim_{\substack{ (p, \zeta) \to (\up, \uzeta) \\ \gamma > 0} }
   \EE_-(p, \zeta) 
\end{equation}
{}From its definition, geometric regularity is a local property, and 
remains satisfied in a neighborhood of the given point. Therefore, for 
$(p, \zeta) $ in a neighborhood of $(\up, \uzeta)$ with 
$\gamma \ge 0$, the   
$\EE_{k, -}(p', \zeta')$  have limits $\tilde \EE_{k, -} (p, \zeta)$
as  $(p', \zeta') \to (p, \zeta)$   with $\gamma' > 0$. 
Arguing as in \cite{MeZ.2}, this implies that 
the vector bundle $\EE_{k,-}$ has a continous extension to 
$\gamma = 0$ on a neighborhood of $(\up, \uzeta)$. 

 Therefore, with notations as in Proposition~\ref{prop513}, there is a smooth  splitting
 $$
 \CC^{N_k}  = \underline \EE \oplus \EE_{k, -}(p, \zeta)
 $$
for $(p, \zeta)$ close to $(\up, \uzeta)$. With projectors associated to this
decompositions,  by continuity, there holds  locally 
$$
 \Sigma^\kappa_k  (\up, \uzeta) \ge \frac{\kappa }{2}    \Pi_{k,+}^*   \Pi_ {k,+}  -  2
 \Pi_{k, -}^*   \Pi_{k, -}. 
$$
This shows  $ \Sigma^\kappa _k (p, \zeta)$ is a smooth K-family of symmetrizers for $G_k$. 
 Summing up, we have proved: 
 
 \begin{cor}
 \label{cor52}
  Suppose that  $G_k$ satisfies the block structure condition at 
   $(\up, \uzeta)$. Then, 
   
   i) there is  a neighborhood $\Omega$ of $(\up, \uzeta)$ such that the 
   vector bundle $\EE_{k,-}$ has a continous extension to 
$  \Omega \cap \{ \gamma \ge 0 \}$; 

ii)  there are K-families of smooth 
 symmetrizers  $\Sigma^\kappa_k (p, \zeta)$  near $(\up, \uzeta)$. 
 \end{cor}
 
 
 \subsection{Totally nonglancing modes}

 It remains to  consider the case where 
$\uga = 0$ and  the blocks $G_k$ are associated  with 
non-geometrically regular
real roots of the characteristic equations.   We first consider   totally
  nonglancing roots.

\begin{prop}
\label{prop41}
 
ii) If    $( \up, \utau, \ueta, \uxi)$  is totally incoming, then,  for 
$(p, \zeta)$ in a neighborhood of $(\up, \uzeta)$, $\EE_-  (p, \zeta) = \EE(p, \zeta)$

iii) If    $( \up, \utau, \ueta, \uxi)$  is totally outgoing, then,  for 
$(p, \zeta)$ in a neighborhood of $(\up, \uzeta)$, $\EE_- (p, \zeta) = \{0\}$.
\end{prop}

 \begin{proof}
   We omit the parameter $p$ in the notation below. 
   We have to show  that when
 $\gamma > 0$, $\Delta (\tau- i \gamma, \eta, \xi) $ has $m$ roots in $\xi$ close to 
 $\uxi$ in the half space
 $\{ \im \xi > 0 \} $  [resp.,  in   $\{ \im \xi > 0 \} $]  in the totally incoming [resp., outgoing] case. 
 Since $\Delta$   has no real roots when $\gamma > 0$, 
 the number of these roots is constant for 
 $\zeta$ close to 
 $\uzeta$  and $\gamma > 0$, and it is sufficient to count them
 when $\tau = \utau $ and $\eta= \ueta$. 
 
 By \eqref{310b}, there holds
 \begin{equation}
 \label{313b}
\underline  \Delta^{(m)} ( \tau, 0 , \xi)  =   \sum_{j= 0}^m 
 \underline a_{m-j}  \tau^j \xi^{m-j}  = 
\underline a_0  \prod_{j= 1}^m \big( \tau + c_j \xi \big)  \,, 
 \end{equation}
  where $\underline  a_0 \ne 0 $ and the second equality is the definition 
 of the coefficients $c_j$.  The hyperbolicity  of $\underline \Delta^{(m)}$  
 implies that all the $c_j$ are real and the nonglancing condition 
 implies that  the $c_j$ do not vanish. The roots of 
 $\underline \Delta^{(m)}  (- i \gamma, 0, \xi) = 0$ are 
 $\xi = i \gamma / c_j$.  
 Since 
 $$
 \Delta( \utau -i \gamma, 0, \uxi + \gamma  \cx) = e(0, 0, 0) 
 \gamma^m \Big(\underline  \Delta^{(m)} ( - i, 0, \cx) + O(\gamma) \Big),   
 $$
 we see that the roots of $\Delta(\utau - i \gamma, 0,   \xi ) = 0$ close to 
 $\uxi $ are 
 \begin{equation}
 \label{314}
 \xi = \uxi + \frac{i \gamma}{ c_j}  + o( \gamma). 
\end{equation}
 The 
 assumption is that all the $c_j$  have the same sign. Thus, if the $c_j$ are positive [resp., negative], all the roots are in 
  $\{ \im \xi > 0 \} $  [resp.,  in   $\{ \im \xi > 0 \} $]  for $\gamma > 0$ small.
 \end{proof}

 \begin{rem}
 \textup{When $-\utau$ is semi-simple of multiplicity $m$, the geometric multiplicity of 
 $- \uxi$ as en eigenvalue of $G(\up, \uzeta) $ is $m$. 
 The nonglancing hypothesis is that  the algebraic multiplicity is also $m$. 
  In this case, 
 the proof above shows that the condition that all the $c_j$ have the same 
 sign is necessary and sufficient  for having 
 $\EE_- = \EE$ or $\EE_-= \{ 0 \}$.   }
 \end{rem}

  \begin{prop}
  \label{prop515}
  Suppose that $L$ is symmetric hyperbolic and 
  $G_k$ is associated to a totally nonglancing root.  Then, 
   there are smooth K-families of symmetrizers for $G_k$  near 
 $(\up, \uzeta)$.

  \end{prop}

\begin{proof}
  By assumption, there is a definite positive matrix $S(p)$ such that 
the $S A_j$ are symmetric and there is a 
$N \times m$ matrix $V_k$ such that for $(p, \zeta)$ close to
$(\up, \uzeta)$: 
$$
\EE_k  = V_k  \CC^m \,, \quad 
V_k G_k = G V_k 
$$
By Remarks~\ref{lem56} and \ref{rem413}, the symmetric matrices 
\begin{equation*}
\Sigma_k  (p, \zeta) =  - V_k^*(p, \zeta)  S(p)  A_d(p)  V_k(p, \zeta)\,.  
\end{equation*}
are symmetrizers for $G_k$. More precisely, there holds 
\begin{equation}
\label{eq59x}
\im \Sigma_k G_k = \gamma V_k^*  S V_k \ge  c_k \gamma \Id , \quad c_k > 0. 
\end{equation}

 With notations as in \eqref{n211}, 
 denote by $\underline A_{k, d} = \underline W_k A_d(\up) \underline V_k $ 
 the boundary matrix of  
 the tangent system  $\underline L'_k $ at the real root  $(\utau, \ueta, \uxi_k)$ under consideration.
By   \eqref{n28b}, the symmetry implies that $\underline W _k 
= \underline V_k^* S(\up)$, hence  
 \begin{equation}
\label{eq524}
\Sigma_k (\up, \uzeta)   =   -  \underline A_{k, d} . 
\end{equation}
 By Lemma \ref{lem35} $\underline A_{k, d}$ is definite positive [resp. negative]
 when the mode is totally incoming [resp. outgoing].  
 With Proposition~\ref{prop41}, this implies that 
 \begin{equation}
\Sigma^\kappa _k = \left\{\begin{array}{rl}
   \Sigma_k  & \quad \mathrm{in \ the \ incoming \ case, } 
\\
\kappa \Sigma_k  & \quad \mathrm{in \ the \ outgoing \ case. } 
\end{array}
\right. 
\end{equation}
are K-familes of symmetrizers for $G_k$. 
  \end{proof}

With Corollary  \eqref{cor52} and the results of Section~4, this implies 
the next theorem. 
 
   \begin{theo}[Second Main Theorem]
\label{2main}
 
 Suppose that $L$ is symmetric hyperbolic in the sense of Friedrichs and that 
  all the real roots $(\up, \utau, \ueta, \uxi)$ of the characteristic equation 
 are either geometrically regular or  totally nonglancing.  Then, there is a neighborhood
 $\omega$ of $\up $ such that: 
 
 	i) there there are smooth  K-families of symmetrizers $\Sigma^\kappa$ on 
 $\omega  \times S^{d}_+$;
 
 	ii) the vector bundle $\EE_-(p, \zeta)$  has a continuous extension to $\gamma = 0$;  
 
 iii)  if in addition the boundary conditions  satisfy the uniform Lopatinski, 
then there are smooth Kreiss symmetrizers and 
  the maximal estimates $\eqref{eq54}$ are satisfied. 
 \end{theo}

 \begin{exam}
 \label{ex517}
 \textup{The  assumptions are  satisfied by the equations of MHD, under appropriate conditions 
 on the parameters.  Indeed, all points are either geometrically
regular or algebraically regular and totally nonglancing; see Appendix A.
}
 \end{exam}

 \begin{rem}
 \label{rem58} 
 \textup{In Proposition~\ref{prop515} and Theorem~\ref{2main}, the symmetry 
 of  $L$ is used at only one place, to construct  definite (positive or negative) symmetrizers
 for the blocks 
 which are not geometrically regular. Thus, the symmetry assumption on $L$ can be relaxed, 
 as soon as one can construct such symmetrizers. For instance, this can be done for 
 double linearly splitting eigenvalues, see Section 6. 
 } 
  \end{rem}
 
 \begin{rem}
 \label{rem59x}
\textup{ The symmetrizers constructed in Corollary~\ref{cor52} and Theorem~\ref{2main}
are obtained by localization and block reduction. Gluing the 
properties \eqref{nn55} for elliptic modes and 
\eqref{nn54} or \eqref{eq59x} for real roots, we see that 
the symmetrizers satisfy the following condition which implies \eqref{global43}:
\begin{equation}
\label{eq512x}
\im \Sigma G = \sum V_j^* E_j V_j \,,   
\end{equation}
with either $E_j $ definite positive on the support of $V_j$ or $E_j =\gamma E_j^1$ with 
$E_j^1$ definite positive on the support of $V_j$. This improvement is useful in applications 
to variables coefficients and nonlinear problems. }
 \end{rem}

 
  \subsection{Linearly splitting modes}
  
  Next we consider  a block  $G_k$   associated to  
 a  linearly splitting semi simple nonglancing eigenvalue: we are given a smooth manifold
 $\cM$ passing through $(\up, \ueta, \uxi_k)$ as in Definition~\ref{def213}.

  \begin{prop}
  \label{prop518}
  Suppose that $G_k$ is associated to a   linearly splitting eigenvalue 
  transversally to a smooth manifold, and assume that it is  nonglancing. 
  Then,  there are neighborhoods 
  $\Omega^\kappa$ of $(\up, \uzeta)$ and K-families of bounded symmetrizers $\Sigma^\kappa$ on 
  $\Omega^\kappa_+ = \Omega^\kappa \cap \{ \gamma > 0 \}$. 
  \end{prop}
  
  The main difference with the previous cases is that the symmetrizers $\Sigma^\kappa$
  are bounded, but not necessarily smooth, that is, they may have no continuous extensions
  to $\overline \Omega^\kappa_+$. 

 \begin{proof}
 {\bf a) } 
 We first consider the case where   
 \begin{equation}
\label{eq513x}
dx \notin T_{\ueta, \uxi_k} \cM_{\up}. 
\end{equation}
   In local coordinates $(q, \tilde \zeta)$ as in \eqref{n47}, the block $G_k$ has the 
   form 
   \begin{equation}
\label{eq514x}
G_k  = \mu(q) \Id + \tilde G_k(q, \tilde \zeta), \quad \mathrm{with} \  \tilde G(q, 0) = 0. 
\end{equation}
Introduce polar coordinates :
\begin{equation}
\label{eq515xpolar}
\tilde \zeta = \rho \cx , \quad \rho = \vert \tilde \zeta \vert, \  \  \cz \in S^{\nu}_+ ,
\end{equation}
where $S^\nu$ is the unit sphere in $\RR^{\nu + 1}$ where 
$\tilde \zeta $ lives, and $S^\nu_+$ is the half sphere 
$\cg > 0$. With $\ccG$  as in \eqref{nn335}, there holds
 \begin{equation}
\label{532c}
G_k(p, \zeta) = \mu(q) \Id +  \rho  \ccG_k(q, \rho, \cz  ).  
\end{equation}
Introduce next the matrix $\ccA_k(q, \rho, \ch, \cx)$  as in \eqref{nn335}.
By Proposition \ref{prop310} and the strict hyperbolicity of $\ccA$, there holds 
for $(q, \rho)$ in a neigborhood of $(\uq, 0)$ and $\cz \in S^{\nu}$: 
\begin{equation}
\label{eq515x}
\det \big( \xi \Id + \ccG_k(q, \rho, \cz) = e(q, \rho\cz, \rho \cx)  
\prod \big( (\ct - i \cg) + \check \lambda_j( q, \rho, \ch, \cx) \big),
\end{equation}
where $e(\uq, 0, 0) \ne 0$ and the $\check \lambda_j$ are 
smooth, analytic in $\cx$, and real when $\cx$ is real. 
Moreover, there are smooth vectors $e_j(q, \rho, \ch, \cx)$, analytic in $\xi$
and linearly independent, such that 
\begin{equation}
\label{eq526x}
\big(  \xi \Id + \ccG_k(q, \rho, \ct^j, \ch, \cg^j) \big) e_j (q, \rho, \ch, \cx) = 0 , \quad 
\ct^j - i \cg^j = -\check \lambda_j(q, \rho, \ch, \cx). 
\end{equation}
By 
Theorem~\ref{theo31}, 
this shows that $\ccG(q, \rho, \cz)$ satisfies the block structure 
condition near $(\uq, 0, \cz)$, for all $\cz \in \overline S^{\nu}_+$. 
Therefore, the Kreiss construction   applies to $\ccG$; Lemmas~\ref{lem511}, 
\ref{lem512} and Corollary~\ref{cor52} imply the following. 

\begin{prop}
\label{prop512x}
With notations as above, there is a neighborhood $\tilde \omega$ of 
$(\uq, 0)$ such that: 

i) the vector bundle of  invariant spaces $\ccEE_{k, -}(q, \rho, \cz)$ associated 
 to eigenvalues of $\ccG_k$ in $\{\re \check \mu < 0 \}$, 
 defined for $\cg > 0$, has 
 a continuous extension to $\tilde \omega \times \overline S^\nu_+$; its dimension 
 is equal to the number of positive eigenvalues of  $\underline A_{k, d}$;  
 
 ii) there are smooth K-families  of symmetrizers $\ccSi^\kappa_k(q, \rho, \cz)$ 
 for $\ccG_k$ on $\tilde \omega \times \overline S^\nu_+$. 
\end{prop}

 By \eqref{532c}, the negative space of $G_k(p, \zeta)$ is  
 \begin{equation}
\label{eq519x}
\EE_{k, -} (p, \zeta) = \ccEE_{k, -} \Big(q, \vert \tilde \zeta \vert , \frac{\tilde \zeta}
{\vert \tilde \zeta \vert } \Big) . 
\end{equation}
Since $\mu(q)$ is real, it follows that 
\begin{equation}
\label{eq520x}
\Sigma^\kappa_k( p, \zeta) = \ccSi^\kappa_k \big(q, \vert \tilde \zeta \vert, \frac{\tilde \zeta} 
{\vert \tilde \zeta \vert} \big) 
\end{equation}
is a K-family of bounded symmetrizers for $G_k$, for 
$(p, \zeta)$ close to $(\up, \uzeta)$ and $\gamma > $.  

\medbreak
{\bf b) } Next we consider the case where 
\begin{equation}
\label{eq521x}
dx \in T_{\ueta, \uxi_k} \cM_{\up}. 
\end{equation}
By \eqref{nn339}, the normal matrix of the tangent system at
$(\uq, \utau, \ueta, \uxi_k)$ is
\begin{equation}
\label{eq522x}
\underline A_{k, d} = \underline a_k \Id \, , \quad  \underline a_k \ne 0. 
\end{equation}
   In local coordinates $(q, \tau, \tilde \zeta)$ as in \eqref{nn342}, the block $G_k$ for 
   $\gamma = 0$, has the 
   form 
   \begin{equation}
\label{eq523x}
G_k {}_{\vert \gamma = 0}
 = \mu(q, \tau) \Id + \tilde G_k(q, \tau , \tilde \eta), \quad \mathrm{with} \  \tilde G_k(q, \tau, 0) = 0. 
\end{equation}
We switch to polar coordinates $\tilde \eta = \rho \ch$ with $\ch \in S^{\nu}$ and 
introduce $\ccG_k$ as in \eqref{nn343}. 
By Proposition~\ref{prop311}, the eigenvalues of $\ccG$ are real and simple
when $\ch \ne 0$. Therefore, there exists a neigborhood 
$\omega$ of $(\uq, \utau)$, a constant $c> 0$,  and a smooth symmetric matrix
$\ccSi_k (q, \tau, \rho, \ch)$ on $\omega \times S^\nu$, such that 
\begin{equation}
\label{eq524x}
\begin{aligned}
\ccSi_k(q, \tau, \rho, \ch)  \ge c \Id, \quad \im (\ccSi \ccG ) = 0. 
\end{aligned}
\end{equation}
Introduce 
\begin{equation}
\label{eq525x}
\begin{aligned}
\Sigma'_k  ( p, \tau, \eta, \gamma) = 
\left\{ \begin{matrix}\ccSi_k \big(q, \tau,  \vert \tilde\eta \vert,  \frac{\tilde \eta} 
{\vert \tilde \eta \vert} \big)  & \quad \mathrm{when} \ \tilde \eta \ne 0, 
\\ 
\Id & \quad  \mathrm{when} \ \tilde \eta =  0. 
\end{matrix}
\right.
\end{aligned}
\end{equation}
These matrices are uniformly bounded for 
$(p, \zeta)$ in a neighborhood of $(\up, \uzeta)$ and  
\begin{equation}
\label{eq526xtwo}
\begin{aligned}
\Sigma'_k(q, \zeta)  \ge c \Id, \quad \im (\Sigma'_k  G_k  )_{ \vert \gamma = 0}  = 0  . 
\end{aligned}
\end{equation}
Symmetrizer $\Sigma'_k(q, \zeta)$ is independent of $\gamma$ and 
$G_k$ is smooth; using  Proposition~\ref{prop311},    this implies 
\begin{equation}
\label{eq527x}
\im (\Sigma'_k  G_k  ) =   - \frac{\gamma }{\underline a} \Sigma_k  + 
O\big(\gamma \vert (p, \zeta) - (\up, \uzeta)  \vert \big)  + O( \gamma^2).  
\end{equation}
Therefore, 
$\Sigma_k = -   \underline a \Sigma' _k$ is a bounded symmetrizer for 
$G_k$ on a neighborhood of $(\up, \uzeta)$ which is positive definite [resp. negative]
when $\underline a$ is negative [resp. positive].
In the first case, by \eqref{eq522x} and Lemma~\ref{lem35}, the mode is totally outgoing  
and by Proposition~\ref{prop41}, the negative space 
$\EE_{k, -}$ is $\{0 \}$. 
In the second case, the mode is totally incoming and $\EE_{k, -} = \EE_k$. 
Thus, $\Sigma_k^\kappa = \kappa \Sigma_k$ in the outgoing case
and $\Sigma_k^\kappa =  \Sigma_k$ in the incoming case provide us 
with K-families of bounded symmetrizers for $G_k$. 
 \end{proof}

  \begin{rem}
  \textup{ If the mode is totally nonglancing, the space
  $\EE_{k,-}$ is either $ \CC^{N_k}$ or $\{ 0 \}$; it is 
  trivially continous at $(\up, \uzeta)$. 
  If the mode  is not totally nonglancing, then the matrix 
  $\underline A_{k, d}$ has positive and negative eigenvalues. 
 In this case, the negative space is given by \eqref{eq519x}: 
 its dimension is positive and smaller then $N_k$, and it is a smooth function 
 of $\tilde \zeta/ \vert \tilde \zeta \vert$. 
 In general, 
  $\EE(p, \zeta)$ is not continuous at $(\up, \uzeta)$ since 
  the limit   depends on the direction of $\tilde \zeta$. 
  }
  \end{rem}

 The Main Theorem~\ref{theo52} is a consequence 
 of  Corollary~\ref{cor52}, Propositions~\ref{prop515} and  \ref{prop512x}, 
 using  the general results of Section 4.


 \section{Smooth Kreiss symmetrizers}

We next investigate when symmetrizers may be chosen in 
\emph{smooth fashion}; this is important for the treatment
of variable-coefficient or nonlinear problems.
Consider again the boundary value problem
\begin{equation}
\label{sk1}
\D_x u + i G(p, \zeta) u = f, \quad  M(p, \zeta)  u_{\vert x = 0 }  = g,
\end{equation}
for $(p, \zeta)$ in a neighborhood of $(\up, \uzeta)$ with 
$\underline \gamma = 0$ and $G$ associated to a hyperbolic system $L$ as in 
\eqref{nn32}. 

We suppose that we are given two invariant spaces 
$\EE_0(p, \zeta)$ and $\EE_1(p, \zeta)$ for $G(p, \zeta)$, which depend smoothly
on $(p, \zeta)$   and such that 
$\CC^N = \EE_0 \oplus \EE_1$. 
We denote by $G_k$ the restriction of $G$ to $\EE_k$ and by 
$\EE_{k, -}$ the negative space of $G_k$, for $\gamma > 0$. 
The negative space for $G$ is $\EE_- = \EE_{0, -} \oplus \EE_{1,-}$. 
 
\begin{ass}
\label{asssk1}
$\EE_{0, -} (p, \zeta)$ has a continuous extension to $\gamma = 0$ near 
$(\up, \uzeta)$. Moreover,  $\dim \ker M =  N - \dim \EE_-$ and 
\begin{equation}
\label{sk2}
\ker M(\up, \uzeta) \cap \EE_{0, -}(\up, \uzeta) = \{ 0 \}. 
\end{equation}
\end{ass}
The latter two conditions are implied by the uniform Lopatinski condition.  

\subsection{Reduced boundary value problems}

Introduce 
 \begin{equation}
 \label{sk3}
 \FF_0 (p, \zeta) = M(p, \zeta)  \EE_{0, -} (p, \zeta).
\end{equation}
This is a continuous bundle on $\overline \Omega_+$ where 
$\Omega $ is a neighborhood of $(\up, \uzeta)$ and 
$\Omega_+ = \Omega \cap \{ \gamma > 0 \}$. By Assumption \ref{asssk1}, 
$\dim \FF_0 = \dim \EE_{0, -}$. 
Let $\FF_1 $ denote a  continuous bundle such that 
\begin{equation}
\label{sk4}
M \CC^N =  \FF_0 \oplus \FF_1. 
\end{equation}
 
 \begin{defi}
 \label{redbc}
  The reduced boundary conditions for block $\EE_j$ is 
  \begin{equation}
\label{sk5}
M_j   = \pi_j  M _{\vert \EE_j} 
\end{equation} 
where   $\pi_j$ denotes the projection on $\FF_j$ in the splitting $\eqref{sk4}$.
Likewise, we define  reduced lopatinski determinants 
$\Delta^j:= \det|_{\EE^j}(\EE^j_-, \ker M^j )$ 
and a
reduced uniform Lopatinski condition consisting of a uniform lower bound
on the modulus of the reduced Lopatinski determinant. 
\end{defi}

\begin{rem}\label{asym}
\textup{
Note that there is a useful asymmetry in Definition
\ref{redbc}, in that we do not require continuity of $E_{1,-}$
at $\uzeta$; indeed, $\EE_{1, -} $ is mentioned only in the Lopatinski
determinant evaluated for $\gamma >0$.  This allows us to
separate off discontinuous blocks.
The continuous block $\EE_0$ may be further subdivided into
arbitrary many blocks on which the negative subspace is continuous,
each with its own reduced boundary condition.}
\end{rem}

\begin{rem}
\label{remsk4}
\textup{The choice of $\FF_0$ implies that 
\begin{equation}
\label{sk6}
\pi_1 M _{\vert \EE_{0, -} }  = 0. 
\end{equation}
On the other hand, $\pi_0 M $ is not necessarily $0$ on $\EE_{1,-}$.  
This is another asymmetry in the role of $\EE_0$ and $\EE_1$.
It reflects that the boundary conditions do not decouple in the block reduction
of $G$. }
\end{rem}

\begin{prop}\label{redlop}
The uniform Lopatinski condition is satisfied on the full space 
$\EE$ if and only if the reduced uniform Lopatinski conditions 
are satisfied on each $\EE_j$.
\end{prop}

\begin{proof} 
{\bf a)  }  If the reduced problems satisfy the uniform Lopatinski condition on 
$\Omega_+$, there are constants $C_j$ such that for all 
$(p, \zeta) \in \Omega_+$, 
\begin{equation}
\label{sk7}
\vert u_j \vert \le C_j  \vert M_j u_j \vert    \quad \mathrm{for} \ u_j \in \EE_{j,-}. 
\end{equation}
Let $u  = (u_{0}, u_{1}) \in \EE_-$.  Then, using \eqref{sk6}, we have
\begin{equation*}
\begin{aligned}
& \vert u_1  \vert \le C_1 \vert \pi_1 M u  \vert  \le C_1 \vert \pi_1 \vert \, \vert M u \vert \quad 
\\
&
\vert u_0  \vert \le C_0  \vert   M u_0  \vert \le 
C_0 \vert M u \vert + C_0 \vert M \vert  \, \vert u_1 \vert ,
\end{aligned}
\end{equation*}
implying  that 
\begin{equation}
\label{sk8}
\vert u \vert \le C   \vert M  u  \vert    \quad \mathrm{for} \  u \in \EE_{-}. 
\end{equation}

{\bf b) }  Conversely, if the uniform Lopatinski condition is satisfied on $\Omega_+$, 
then, since $M = \pi_0 M$ on 
$\EE_{0,-}$, \eqref{sk8}  restricted to $\EE_{0, -}$ implies that 
\begin{equation*}
  | u_0 |\le  C  |M  u_ 0| =  C  |  M_ 0  u^0_-|  \quad  \mathrm{for} \ u_0 \in \EE_{0, -}. 
\end{equation*}
verifying the reduced Lopatinski condition on $\EE^0$.
Next, by definition of $\FF_0$, there  is an inverse $R$ of $M$ from 
$\FF_0 $ to $\EE_{0,-}$.  
By continuity of $\EE_{0,-}$, this inverse is uniformly bounded on a neighborhood
of $(\up, \uzeta)$. 
For   $  u_1 \in \EE_{1, -} $, let  $u_0 = R \pi_0 M u_1  \in \EE_{0, -} $  
such that $u = (u_0, u_1) \in \EE_-$, $ \pi_0 M  u  =0$ and $ M u = \pi_1 M u = M_1 u_1$.
Then 
\begin{equation}
 \vert u_1 \vert \le  C_2 \vert u \vert \le C_2 C \vert M u \vert  = C_2 C \vert M_1 u_1 \vert. 
\end{equation}
where $C_2$ is the norm of the projection 
from $\CC^N$ on $\EE_0$. This  verifies the reduced Lopatinski condition on $\EE^1$.
\end{proof}

 
\subsection{Symmetrizers}

We shall seek  symmetrizers commuting with the projectors
onto $\EE_j$; we call such a symmetrizer {\it consistent}
with the decomposition $\EE=\EE_0 \oplus \EE_1$. 
Note that if $\Sigma$ is a symmetrizer for $G$, then 
the  $\Sigma_j$ defined by restricting the hermitian form 
$(\Sigma u, v)$ to $\EE_j $ is a symmetrizer for $G_j$ (see Remark~\ref{rem413}). 
Thus, $\widetilde \Sigma  = \Sigma_0 \oplus \Sigma_1$ is a symmetrizer for 
$G$, consistent with the decomposisiton.

We have shown in Section 5 that the notion of  K-families of symmetrizers 
is well adapted to a block decomposition: if there are such families 
$\Sigma^\kappa_k$ for each block,  
then $\Sigma^\kappa = \mathrm{diag } ( \Sigma^\kappa_k) $ is a 
k-family for the full system. Such an analogue is not true for Kreiss symmetrizers: 
a K-family depends only on $G$, while the definition of a Kreiss symmetrizer also 
depends on the boundary conditions, which do no   decouple, in general, in the block
reduction of $G$. However, we have the following  useful partial result.

\begin{prop} 
 \label{propsk7}
Consider  $G$ with boundary condition $M$,
and a specified decomposition $\EE_0$, $\EE_1$ in a neighborhood
of $(\up, \uzeta) $ satisfying Assumption $\ref{asssk1}$.  
If there exists a local smooth Kreiss symmetrizer  $\Sigma$ respecting this decomposition 
then  

\quad i) the reduced boundary value problem on $\EE^0$
satisfies the (reduced) uniform Lopatinski
condition, 

\quad ii) the restriction $\Sigma_1$ of $\Sigma$ to $\EE_1$ 
is a smooth Kreiss symmetrizer    for  the
  reduced boundary value problem $(G_1, M_1)$. 
  
  Conversely, assume ii) and 
  
  \quad iii) the reduced boundary value problem on $\EE^0$
satisfies the  uniform Lopatinski condition and   there exists a smooth K-family of symmetrizers for the reduced problem
  $(G_0, M_0)$. 
  
  Then, there exists a local smooth Kreiss symmetrizer for $(G,M)$. 
\end{prop}

\begin{proof}
 {\bf  a) }  
Suppose that $\Sigma =\bdiag\{\Sigma_0,\Sigma_1\}$ is a local smooth Kreiss symmetrizer 
for $(G, M)$. Then there are $C >0$ and $\eps > 0$ such that 
\begin{equation}
\label{sk10}
\langle u, \Sigma u\rangle= 
\langle u_0, \Sigma_0u_0\rangle+ \langle u_1, \Sigma_1u_1\rangle \ge \epsilon |u|^2
-C|Mu|^2. 
\end{equation}
We have recalled in Section 4 that the existence of a Kreiss symmetrizer 
implies that the uniform Lopatinski condition must hold, and thus, by Proposition \ref{redlop}, 
the uniform reduced Lopatinski condition must hold on $\EE^0$. 

Further, for all $u \in \EE_1$, there  is  $\tilde u_{0,-} \in \EE_{0, -} $ 
such  that $\pi_0 M(u^1+\tilde u_{0, -})  =0$. 
By  Lemma \ref{sign} and continuity of 
$\EE_{0, -}$, $\langle \Sigma_0 u_{0,-}, u_{0,-} \rangle \le 0$. 
Moreover,   by \eqref{sk6}, 
$M ( u_1 + u_{0, -}) = \pi_1 M u_1 = M_1 u_1$. Thus,  
\begin{equation}
\begin{aligned}
\langle u_1, \Sigma_1 u_1\rangle &\ge 
\langle u_1, \Sigma_1 u_1\rangle+
\langle \tilde u_{0, -} , \Sigma_0 \tilde u_{0,-} \rangle \\
&=
\langle (u_1+\tilde u_{0,-}), \Sigma (u_1+\tilde u_{0,-})\rangle \\
&\ge \epsilon (|u_1|^2+|\tilde u_{0,-}|^2) -C|M(u_1+\tilde u_{0, -} )|^2 
\\
&= \epsilon |u^1|^2 -C|M^1u^1|^2,
\end{aligned}
\end{equation}
verifying ii).

 {\bf  b)  }  Conversely,  let $\Sigma^\kappa_0$ denote a K-family 
 of symmetrizers for $(G_0, M_0)$. Then,  by Proposition \ref{prop58}, 
 for $(p, \zeta) \in \Omega^\kappa_+$, the uniform Lopatinski condition implies that  
 \begin{equation}\label{s0rel}
\langle u_0, \Sigma^\kappa_0 u_ 0\rangle  + C_0 \vert M_0 u_0 \vert^2 
\ge  \vert u_0 \vert^2 + 
  m(\kappa) |u_{0, +} |^2  
\end{equation}
where $m(\kappa) \to +\infty$ as $\kappa \to +\infty$. 
 By assumption, there is also   $S_1$ satisying    
\begin{equation}
\langle u_1, \Sigma_1 u_1\rangle \ge |u_1|^2 - C_1|M_1u_1|^2 . 
\end{equation}
Therefore, 
\begin{equation*}
\begin{aligned}
\delta  \la u_0 , \Sigma_0 u_0 \ra +   \la u_1 , \Sigma_1 u_1 \ra     
+ C_0 \delta \vert M_0 u_0 \vert^2 & + C_1 \vert M_1 u_1 \vert^2  \ge  
 \\
   \vert u_1 \vert^2  & + \delta m(\kappa) |u_{0, +} |^2 +  \delta |u_{0} |^2. 
 \end{aligned} 
\end{equation*}
Using  \eqref{sk6}, we have for $u = (u_0, u_1)$: 
\begin{equation*}
\begin{aligned}
2       \vert \pi_1 M u \vert ^2    
 &      
\ge  \vert M_1 u_1 \vert^2 -  C_2    \vert   u_{0, +} \vert^2 
\\ 
  2  \vert \pi_0 M u \vert^2  & \ge    \vert M_0 u_0  \vert^2  - C_3 \vert u_1 \vert^2. 
\end{aligned}
\end{equation*}
Therefore, 
\begin{equation*}
\begin{aligned}
\delta  \la u_0 , \Sigma_0 u_0 \ra  & +   \la u_1 , \Sigma_1 u_1 \ra     
+    2 C_0 \delta     \vert \pi_0 M u   \vert^2   + 2 C_1 \vert \pi_1 M  u  \vert^2  \ge  
 \\
   \vert u_1 \vert^2  & +   \delta |u_{0} |^2 + \delta m (\kappa) |u_{0, +} |^2  - C_3 
 \delta   \vert u_1 \vert^2  - C_2 \vert u_{0, +} \vert^2 . 
 \end{aligned} 
\end{equation*}
Hence, if $C_3 \delta \le 1/2$,  $ \delta m(\kappa) \ge C_2$, 
$C \ge 2 C_0 \delta \vert \pi_0 \vert^2 +  2 C_1 \vert \pi_1 \vert^2$
and 
$\Sigma = \mathrm{diag} (\delta \Sigma_0, \Sigma_1) $: 
\begin{equation*}
 \la u , \Sigma u \ra + C \vert M u \vert^2 \ge \delta \vert u \vert^2 . 
\end{equation*}
 \end{proof}

\begin{rem}
\textup{
 If $G_0$ satisfies the block structure condition, then 
 $i)$ and $iii)$ are equivalent. }
 \end{rem}


\subsection{Smooth symmetrizers for $\EE_1$}
 
 
 We investigate the existence of  smooth Kreiss symmetrizers for blocks
$\EE_1$ associated to nonregular eigenvalues. 
We have in mind the case of linearly splitting eigenvalues. 
We first give necessary and next sufficient conditions. 
We assume, as we may, that $G_1$ has a block diagonal form 
\begin{equation}
\label{eq614x}
G_1(p, \zeta) = \bdiag ( G_{1,k}(p, \zeta) ) 
\end{equation}
where the $G_{1,k}$ are some  of the invariant blocks of $G(p, \zeta)$, for 
$(p, \zeta)$ close to $(\up, \uzeta)$ associated to pairwise distinct real eigenvalues 
$\underline \mu_k = - \uxi_k$ of $G(\up, \uzeta)$. We further assume that
the modes $(\up, \utau, \ueta, \uxi_k)$ are semi-simple and nonglancing. 
In this case,  by Proposition~\ref{prop34}, the $ \underline \mu_k$  are semi-simple eigenvalues of $G$ and 
\begin{equation}
\label{eq615x}
G_{1,k} (up, \uzeta + \zeta) = \underline \mu_k \Id + 
\underline G'_{1,k} (\zeta) + O( \vert \cz \vert^2). 
\end{equation}
Moreover, denoting by 
$$
\underline L'_{1, k}  = \tau \Id + \sum_{j=1}^{d-1} \eta_j  \underline A_{k, j}  + \xi  \underline A_{k, d}  
$$
 the tangent  system at  $(\up, \utau, \ueta, \uxi_k)$, there holds
\begin{equation}
\label{eq616x}
\underline G_{1, k} =  \underline A_{k, d} ^{-1} \Big( 
(\tau - i \gamma ) \Id + \sum_{j=1}^{d-1} \eta_j  \underline A_{k, j} \Big). 
\end{equation}
The Taylor expansion of $G_1$ at the base point reads  
\begin{equation}
\label{sk17}
G_1 (\up, \uzeta + \cz ) = G_1 (\up, \uzeta) + \underline G'_1 (\cz) + O( \vert \cz \vert^2). 
\end{equation}
$ G_1 (\up, \uzeta) $ and $\underline G'_1 $ are block diagonal; we further introduce 
\begin{equation}
\label{eq618x}
\underline L'_1 = \bdiag (\underline L'_{1, k}) 
\end{equation}

\begin{prop}
\label{prop68x}
Suppose that $\Sigma (p, \zeta)$ is a smooth  symmetrizer for 
$G$.  Let $\Sigma_1$ denote the restriction of $\Sigma$ to 
$\EE_1$. Then  $\Sigma_1 $ is a smooth  symmetrizer for 
$G_1 $ and 

\quad i)  $\underline \Sigma_1 := \Sigma_1( \up, \uzeta)$ is block diagonal and 
 satisfies $\im (\underline \Sigma _1 \underline G_1) = 0$; 
 
\quad ii)    $\underline L'_1$ is hyperbolic symmetric in the sense of Friedichs. 

\end{prop}
 
 \begin{proof}
 Since $\Sigma$ is a symmetrizer, $\im (\Sigma G ) \ge c\gamma $, 
 which by restriction   to $\EE_1$, implies that for $\gamma > 0$: 
 \begin{equation}
 \label{sk17b}
  \im \Sigma _1 G_1    \ge c \gamma \Id_{\EE_1}. 
\end{equation}
Thus, by continuity, 
\begin{equation*}
\im (\underline \Sigma_1 \underline G_1)  \ge 0. 
\end{equation*}
Since $\underline G_1$ is block diagonal, with real entries $\mu_j$, this means that 
for all pair $(k, l)$ with $k \ne l $, the nondiagonal block entries $\underline \Sigma_{1,k,l}$ 
of $\underline \Sigma_1$ are such that 
$$
\begin{pmatrix}
  0     &   P    \\
    P^*      &  0 
\end{pmatrix} , \quad   P = \frac{1}{2i} ( \mu_k - \mu_j)  \underline \Sigma_{j,k} 
$$ 
is nonnegative. This implies that $P=0$, hence  $\underline \Sigma_{1,k,l} = 0$ when $k \ne l$. 
thus 
\begin{equation}
\label{sk18}
\im (\underline \Sigma_1 \underline G_1)  =  0. 
\end{equation}

The diagonal block  $  \Sigma_{1, k, k}$  is a symmetrizer  for $G_{1, k}$ by 
Remark~\ref{rem413}. 
Taking the first order term in $ \zeta$  in
$\im (\Sigma_{1,k,k}  G_{1,k})$,  implies that 
 \begin{equation}
 \label{eq620x}
  \im \underline \Sigma_{1,k, k}    \underline G'_{1,k} (\zeta ) \ge c \gamma     . 
\end{equation}
We have used that the $\underline G_{1,k} =\underline  \mu_k \Id $ and that 
the derivatives of $\Sigma_{1, k, k}$ are self-adjoint. Because $\underline G'_1$ is a linear 
in $\zeta$, changing $\zeta$ to $- \zeta$ implies that 
$ \im \underline \Sigma_{1,k, k}    \underline G'_{1,k} (\zeta ) = 0 $  when $\gamma = 0$. 
Thus, the matrices 
$\underline S_{1, k} :=  - \underline \Sigma_{1,k,k} \underline A_{k, d}^{-1} $
and $ \underline S_{1, k} \underline A_{k, j} $ are self-adjoint. 
Moreover, \eqref{eq620x} implies that $\underline S_{1,k} $ is definite positive. 
This shows that   $ \underline S_{1, k} $ is a Friedrichs symmetrizer for 
$\underline L'_{1,k}$, therefore that 
\begin{equation}
\label{sk622y}
\underline S_1 = \bdiag (\underline S_{1,k})
\end{equation}
  is a Friedrichs symmetrizer for 
$\underline L'_{1}$. 
\end{proof}

\begin{prop}
\label{prop69y}
Suppose that $\Sigma_1$ is a smooth symmetrizer for 
$G_1$. Then, it is a smooth Kreiss symmetrizer on a neighborhood of
$(\up, \uzeta)$ if and only if 
the reduced boundary condition $M_1(\up, \uzeta)$ is maximal 
dissipative for $\underline L'_1$ with symmetrizer \eqref{sk622y}. 

\end{prop}

\begin{proof}
If  
$\Sigma_1$ is a Kreiss symmetrizer for $(G_1, M_1)$, then  
\begin{equation}
\label{sk623y}
\la \Sigma_1 u_1, u_1 \ra  + C \vert M_1 u_1 \vert^2  \ge c \vert u_1 \vert^2 . 
\end{equation}
Specializing at $(\up, \uzeta)$  implies that $\underline S_1 \underline A_d = - \underline \Sigma_1$ is maximal negative 
on $\ker   M_1(\up, \uzeta) $, thus that the boundary condition $M_1(\up, \uzeta)$   is maximal 
dissipative for $\underline L'_1$ with symmetrizer $\underline S_1$. 

Conversely, if $M_1(\up, \uzeta)$   is maximal 
dissipative for $\underline L'_1$ with symmetrizer $\underline S_1$, then 
\eqref{sk623y} holds at $(\up, \uzeta)$, therefore on a neighborhood of that point. 
 \end{proof}

 \begin{theo}
[Third  Main Theorem]
\label{3main}
 Suppose that $L$ is a   hyperbolic system  and $M$ are boundary conditions 
 which satisfy the uniform Lopatinski condition.  For 
 all $(p, \zeta)$ in a neighborhood 
 of $(\up, \uzeta)$, denote by 
 $\EE_0(p, \zeta) $ the invariant space of $G$  corresponding to eigenvalues 
 which are either nonreal, or real and associated to  geometrically regular or totally nonglancing
 modes of $L$. Denote by $\EE_1(p, \zeta)$ the invariant space associated to the remaining 
 eigenvalues of $G$. 
Then, Assumption~$\ref{asssk1}$ is satisfied and, with notations as in Proposition~$\ref{prop69y}$, there exists a smooth Kreiss symmetrizer for 
$(G, M)$ near $(\up, \uzeta)$ if and only if there is a smooth symmetrizer 
$\Sigma_1$ for $G_1$ such that the associated tangent system 
$\underline L'_1$ \eqref{eq618x} is hyperbolic symmetric in the sense of Friedrichs
with symmetrizer $\underline S_1$ \eqref{sk622y} and 
 the reduced boundary condition is maximally dissipative  for 
 $\underline L'_1$. 
  \end{theo}  

\begin{proof}
By Corollary \ref{cor52},  Propositions \ref{prop41}   \ref{prop515} and Lemma \ref{lem510},  there are smooth K-families of symmetrizers 
for $G_0$ and Assumption \ref{asssk1} is satisfied. 

If $\Sigma$ is a smooth Kreiss symmetrizer for $(G, M)$, then by Proposition \ref{propsk7}
$\Sigma_1$ is a smooth Kreiss symmetrizer for the reduced problem  $(G_1, M_1)$ and 
Proposition \ref{prop69y} implies that  the reduced boundary condition is maximally dissipative  for 
 $\underline L'_1$. 

Conversely, if $\Sigma_1$ is a smooth symmetrizer for $G_1$ such that 
the reduced boundary condition is maximally dissipative  for 
 $\underline L'_1$,  then $\Sigma_1$ 
is a 
Kreiss symmetrizer for the reduced problem  $(G_1, M_1)$ and 
Proposition \ref{propsk7} implies that there is a smooth Kreiss symmetrizer for 
the full system $(G,M)$. \end{proof}

This reduces the analysis to the construction of smooth Kreiss symmetrizers for $G_1$.  
 When  $G $ is associated to a \textit{symmetric hyperbolic system $L$}, with 
 symmetrizer $S$, then $\Sigma = - S A_d$ is a smooth symmetrizer for $G$. 
 and the restriction $\Sigma_1$ of $\Sigma $ to $\EE_1$ is a  smooth  symmetrizer
on a neighborhood of $(\up, \uzeta)$. For this symmetrizer, the maximal dissipativity condition 
of the reduced boundary condition  holds if and only if 
 $\Sigma_1(\up, \uzeta)$ is definite positive on 
 $\ker M_1(\up, \uzeta)$, that is if  and only if there is $c > 0$ such that 
 \begin{equation}
\label{sk19}
\quad \forall 
u \in \underline \EE_1  , \quad  \underline M  u \in \underline M  \underline \EE_{0,-}   
\ \Rightarrow \ 
\la  \underline S \underline A_d u , u \ra \ge c \vert u \vert^2 ,  \ 
\end{equation}
  where $\, \underline {\ } \,$ means evaluation 
at the base point $(\up, \uzeta)$. 
Therefore, 

\begin{cor}
\label{cor611y}
With assumptions and notations as in Theorem~$\ref{3main}$, 
if the condition \eqref{sk19} holds at $(\up, \uzeta)$, then there is a 
smooth Kreiss symmetrizer  for $(G, M)$ on a neighborhood of 
$(\up, \uzeta)$. 
\end{cor} 

\begin{rem}
\label{rem612y}
\textup{ Suppose that $L$ is symmetric and split $G_1$ in diagonal blocks    $G_{1,k}$ 
as in \eqref{eq614x}. Then, by Remark \ref{rem413}, there are smooth symmetrizers 
$\Sigma_{1, k} $ for $G_{1,k}$, of the form 
$\Sigma_{1,k} = - V_{1,k}^* ( S A_d) V_{1,k}$. Therefore, there are smooth symmetrizers for 
$G_1$ of the  form
\begin{equation}
\label{sk625y}
\widetilde \Sigma_1 = \bdiag ( \alpha_k \Sigma_{1,k}) , \qquad \alpha_k > 0\,. 
\end{equation}
For these symmetrizers, one can state the analogues of condition \eqref{sk19} and Corollary~\ref{cor611y}.  }
\end{rem}

 \begin{rem}
\label{rem613y}
\textup{ We have seen that a necessary condition for the existence of  a smooth symmetrizer
$\widetilde \Sigma_1$ 
is  that $\underline L'_1$ must be hyperbolic symmetric in the sense of Friedrichs, 
 with symmetrizer 
 $\underline {\widetilde S}_1 = \bdiag  \underline {\widetilde S}_{1,k}$, where 
 $  \underline {\widetilde S}_{1,k} = - \underline {\widetilde \Sigma}_{1,k}  
 \underline A^{-1}_{k, d} $ and  
  $\underline {\widetilde \Sigma}_1 = \bdiag  \underline {\widetilde \Sigma}_{1,k}$. 
  Therefore, in the construction of a smooth Kreiss symmetrizer $\widetilde \Sigma_1$, the first 
  step is construct a  Friedrichs symmetrizer 
  $\underline {\widetilde S}_1 = \bdiag  \underline {\widetilde S}_{1,k}$ for $\underline L'_1$ 
  such that $\underline M_1$ is maximal dissipative. This defines $\underline {\widetilde \Sigma}_1$, 
  and the second  is to determine a smooth symmetrizer for $G_1$ such that 
  $\widetilde \Sigma_1 (\up, \uzeta) = \underline {\widetilde \Sigma}_1$. 
  Since $\underline {\widetilde \Sigma}_1$ is block diagonal, 
  it is sufficient to construct smooth symmetrizers for $G_{1,k}$ such that 
  $\widetilde \Sigma_{1,k} (\up, \uzeta) = \underline {\widetilde \Sigma}_{1,k}$. }

\textup{  
When $L$ is symmetric hyperbolic, the $\underline L'_{1, k}$ have
 Friedrichs symmetrizers $ \underline S_{1, k} $, as in Remark~\ref{rem612y}. 
 In  generic cases,  these symmetrizers are unique up to a constant, implying that 
 necessarily
 \begin{equation}
 \label{626y}
 \underline {\widetilde S}_{1, k} = \alpha_k \underline S_{1, k},  \qquad  \alpha_k > 0. 
 \end{equation} 
 Therefore, 
 \begin{equation}
\label{sk627y}
\underline{\widetilde \Sigma}_1 = \bdiag ( \alpha_k \underline \Sigma_{1,k}).  
\end{equation}
In this case,  if there is a smooth Kreiss symmetrizer, there is one of the 
form \eqref{sk625y}, indicating that this form contains almost all the possible
choices of smooth symmetrizers. 
 }
 \end{rem}

 \begin{rem}
\label{remreducible}
\textup{More precisely,
symmetrizers $\underline {\widetilde S}_{1, k}$ are unique up
to constants if and only if the reduced tangent systems 
$\underline L'_{1, k}$ are {\it irreducible} in the sense
that there do not exist nontrivial constant invariant subspaces of
$\underline{A}_{1,k}(\eta,\xi)$.
For, without loss of generality taking $\underline {A}_{1,k}^d$ diagonal,
so that $\underline{\widetilde S}_{1, k}$ must be diagonal as well,
and expressing $\underline{\widetilde S}_{1, k}$ as a block-diagonal 
matrix with blocks consisting of distinct scalar multiples of the identity,
we find that symmetry of 
$\underline{\widetilde S}_{1, k}\underline{A}_{1,k}^j$ implies
holds if and only if each $\underline{A}_{1,k}^j$ shares the same
block-diagonal structure.
}

\textup{
In the linearly splitting case, we find by spectral separation that
these invariant subspaces must correspond to the limit as
$(\eta,\xi)\to 0$ of group eigenprojections
associated to fixed subsets of the associated eigenvalues 
${\lambda}_{1,k}^j(\eta,\xi)$ of the nearby perturbed system.
In particular, this implies that these group eigenprojections are
continuous along rays across the origin, which is sufficient to give
analyticity separately in coordinates $\eta_j$, $\xi$.  By
Hartog's Theorem, therefore, they are jointly analytic in $(\eta,\xi)$,
and thus there is an analytic decomposition of the perturbed system
into two invariant blocks.  
By an argument similar to that for the
geometrically regular case, we find that, at least for the 
easier nonglancing case, the associated block $G_{1,k}$ 
decomposes analytically as well.
Thus, in the linearly splitting, nonglancing case, we may always
arrange that the tangent system have unique (up to constants) symmetrizers.
}
\end{rem}
 
\subsection{The case of double roots}
 
In this section, we assume that $\EE_1$ is the direct sum of two dimensional 
subspaces $\EE_{1,k} $  associated to  double linearly splitting characteristic 
roots: 
  
\begin{ass}
\label{ass41}
 Suppose that $L$ is symmetric hyperbolic and noncharacteristic
with respect to the boundary $x =0$, and that its characteristic roots 
are either geometrically regular, totally nonglancing, or second order, 
nonglancing and linearly splitting transversally to a smooth manifold 
 of codimension two. 
 \end{ass} 
 
 The following result is proved in Appendix D. 
 
  \begin{prop}
  \label{prop615y} Suppose that $G_{1,k}$ is a $2 \times 2$ block associated 
  to a non\-glancing linearly splitting eigenvalue, near 
  $(\up, \uzeta)$ with $\uga = 0$. Then there exist 
  smooth symmetrizers $\Sigma_{1,k}$ near $(\up, \uzeta)$ for $G_{1,k}$ 
  such that 
  $$\underline S_{1, k} := -   \Sigma_{1, k} (\up, \uzeta) A_{k, d}^{-1}(\up, \uzeta)$$
  is a Friedrichs symmetrizer for $\underline L'_{1,k}$. 
  
  Moreover, the symmetrizer for $\underline L'_{1,k}$ is unique, up to multiplication by a positive constant. 
  \end{prop}
  
  In this case, the only block diagonal symmetrizers of 
  $\underline L'_1$ are 
   \begin{equation}
\label{sk628y}
\underline{\widetilde S}_1 = \bdiag ( \alpha_k \underline  S_{1,k}),   
\end{equation}
with  positive constants $\alpha_k$. Moreover, given the $\alpha_k$, 
 \eqref{sk625y} defines a smooth symmetrizer of $G_1$ such that 
  \begin{equation}
\label{sk628ytwo}
 {\widetilde \Sigma}_1(\up, \uzeta) = \bdiag ( \alpha_k \underline \Sigma_{1,k}).  
\end{equation}
   Therefore, Theorem~\ref{3main} takes the special form: 
   
   \begin{theo}
   \label{theo616y}
   Suppose that $L$ satisfies Assumption~$\ref{ass41}$ and $M$ are boundary conditions 
   which satisfy the uniform Lopatinski conditions. With notations as in Theorem~$\ref{3main}$, 
   there exists a smooth Kreiss symmetrizer for $(G,M)$ if and only if for all 
   $(\up, \uzeta)$ with $\uzeta \in \overline S^d_+$ and $\uga = 0$, the 
   reduced boundary conditions $M_1(\up, \uzeta)$ are maximal dissipative 
   for some Friedrichs symmetrizer of $\underline L'_1$. 
  
  \end{theo}

Moreover, we have by
direct calculation the following useful fact, generalizing
a corresponding observation of Majda-Osher in the case of the
$2\times 2$ wave equation.  For the proof, see Appendix D
and Proposition \ref{gencone} below.
 
 \begin{prop}
  If $G_1$ itself is a $2 \times 2$ matrix, there is a smooth Kreiss symmetrizer  for the reduced problem
  $(G_1, M_1)$, if and only if it satisfies the uniform Lopatinski condition. 
 \end{prop}
 
 This immediately implies the following

 \begin{theo}
[Fourth Main Theorem]
\label{4main}
Suppose that Assumption~$\ref{ass41}$ is satisfied
and  that for all $\uzeta \in \overline S^d_+ $ with $\uga  = 0$, 
there is at most one real eigenvalue associated to a mode which is neither 
geometrically regular nor totally nonglancing. 
 
Then, for boundary conditions $M$ there is a smooth Kreiss symmetrizer 
for $(G, M)$ if and only if $M$ satisfies the uniform Lopatinski condition.
 \end{theo}

We conclude by briefly considering in the simple double-root case,
under what conditions our structural hypotheses hold.

\begin{prop}\label{genreg}
For a symmetrizable system in dimension $d\le 2$,
geometric regularity holds at any point $\xi_0\ne 0$.
\end{prop}

\begin{proof}
Restricting by homogeneity to the unit sphere, we obtain
a symmetric matrix perturbation problem in a single parameter,
for which eigenprojections necessarily vary analytically
\cite{Kat}.
\end{proof}

\begin{prop}\label{gencone}
In dimension $d\ge 3$,
a variable multiplicity point of multiplicity two
is generically nonglancing, nonregular, and linearly splitting
of codimension two.
Moreover, a linearly splitting point of multiplicity two
is either of codimension two
or else geometrically regular of codimension one.
\end{prop}

\begin{proof}
The reduced, symmetrizable 
system in dimension $d-1=2$ 
generically has noncharacteristic
$\tilde A_1$, hence is nonglancing.  
Algebraic regularity can hold only if the reduced,
frozen coefficient system has linearly varying
eigenvalues, in which case $\tilde A_1$ and $\tilde A_2$ commute, a 
degenerate case.
Moreover, if $\sum_{j+1}^{d-1} \xi_jA_j$ has equal eigenvalues for
$\xi\ne (0,0)$, then by
symmetry it is a multiple $\alpha I$ of the identity,
which is three conditions (again taking into account symmetry to
eliminate one entry) on the unknowns $\xi$, $\alpha$, generically
determining a codimension two manifold $M$ in $\xi$.

The latter argument shows also that a linearly splitting point can be
at most codimension two.
Suppose that it is codimension one.  Then, by an analytic change of
coordinates taking the manifold $M$ on which eigenvalues agree to
the $(\xi_1,\dots, \xi_{d-1})$ plane, and subtracting off the 
(analytic) average of the two eigenvalues, we obtain a matrix perturbation
problem of form $\tilde A(\xi)= \xi_d D(\xi)$, where $D(0)$ is diagonalizable
with distinct eigenvalues.  The eigenvalues and associated
eigenprojections of $D$ are therefore analytic, as are the 
eigenvalues and eigenprojections for the original problem.
\end{proof}

\begin{prop}\label{genalg}
In any dimension,
an algebraically regular variable multiplicity point 
of multiplicity two
is either geometrically regular or else has 
characteristics tangent to first order in $\xi$,
in which case nonglancing and totally nonglancing
are equivalent.
\end{prop}

\begin{proof}
Denote the two characteristics as $\tau_1(\xi)$, $\tau_2(\xi)$.
If $\partial (\tau_1-\tau_2)/\partial \xi_d)=0$,
then one of the last two possibilities must occur.
If $\partial (\tau_1-\tau_2)/\partial \xi_d)\ne0$,
on the other hand, then by the Implicit Function Theorem,
there is a hypersurface $\xi_d=\phi(\xi')$ on which
$\tau_1-\tau_2\equiv 0$.
But symmetrizability implies semisimplicity, so that the
reduced matrix $\tilde A(\xi)$ obtained by projection onto the associated
two-dimensional total eigenspace of $\tau_1$, $\tau_2$
must be a multiple of the identity; likewise, $\partial \tilde A/\partial \xi_d$
must be diagonalizable, with distinct eigenvalues.
The matrix $\tilde B(\xi):= (\tilde A-\tau_1I)(\xi)$
is thus of form $(\xi_d-\phi(\xi'))D(\xi')+ \CalO(
|\xi_d-\phi(\xi')|^2)$, analytic in $\xi$, where
$D$ is diagonalizable with distinct eigenvalues,
and it follows that the associated projectors are analytic.
The final assertion follows by the same argument applied to
each direction $\xi_j$.
\end{proof}

Collecting results, we have the following conclusion.

\begin{cor}\label{verification}
Assumption \ref{ass41} is satisfied always in dimension $d=2$.
For multiplicity two crossings in dimension $d=3$,
assumption \ref{ass41} is generically satisfied,
and at algebraically regular points is always satisfied.
\end{cor}

\section{Application to MHD} 

Finally, we discuss our main application, to the shock wave problem in MHD.
As is well known, the Euler equations of MHD are symmetrizable hyperbolic
for a thermodynamically stable equation of state, 
in particular for an ideal gas law; 
see \cite{G, Kaw, KSh, MZ.4, Zhandbook}.
However, they do not satisfy the condition of constant multiplicity, and
so the standard Kreiss--Majda theory cannot be applied; see calculations, 
Appendix A.
Moreover, the boundary condition that arises through linearization of the
shock transmission problem is not dissipative, and so the older theory
of dissipative boundary conditions cannot either be applied; 
see \cite{BT.1}.

For this reason, and also the inherent complexity of the system, progress
on shock stability for MHD has been so far somewhat limited.
In particular, the only nonlinear stability result up to now has been
obtained by Blokhin and Trakhinin \cite{BT.1, BT.4} 
using the ``direct method''
introduced by Blokhin \cite{Bl} for the study of gas dynamics,
in the very special case that magnetic field $H$ goes to zero: essentially,
the fluid-dynamical limit.
This consists of constructing a Lyapunov function in the form of certain
``dissipative integrals'', from which linearized and nonlinear stability
follow directly.
The construction, which is rather intricate, depends crucially on symmetrizability
of the system; as the authors describe it, it consists of augmenting the
original system with various derivatives, chosen in such a way that the
(initially nondissipative) boundary condition becomes dissipative.
At least as presented, it does not appear easily generalizable to other cases.

On the other hand, with additional symmetry, either in the ``transverse'' case
that the magnetic field be orthogonal to fluid velocity, or
the``parallel'' case that they are parallel, the Lopatinski condition has
been explicitly evaluated even for large $H$, to yield physical
stability conditions analogous to those of Majda \cite{Maj} in the gas-dynamical
case;  see \cite{GK, BT.1, BT.2, BT.3, BT.4}.
More generally, the Lopatinski condition has been extensively 
studied numerically for arbitrary parameters, 
with interesting results \cite{BTM.1, BTM.2}.
For a detailed description of the current status of the theory, we refer to
the excellent survey \cite{BT.1}.
Up to now, however, the precise relation between the Lopatinski condition
and linearized and nonlinear stability was not known.

Using the framework developed in this paper,
we can now resolve this issue, at least for the simpler, isentropic case.
We expect that similar considerations hold in the general case.

\subsection{Equations}
The equations of isentropic magnetohydrodynamics (MHD) appear in 
basic form as
\begin{equation}
\label{mhdeq}
\left\{ \begin{aligned}
 & \D_t \rho +  \div (\rho u) = 0
 \\
 &\D_t(\rho  u) + \div(\rho u^tu)+ \na p + H \times \curl H = 0
 \\
 &  
 \D_t H + \curl (H \times u) = 0
 \end{aligned}\right.
\end{equation}
\begin{equation}\label{divfree}
\div H=0,
\end{equation}
where $\rho\in \RR$ represents density, $u\in \RR^3$ fluid velocity,
$p=p(\rho)\in \RR$ pressure, and $H\in \RR^3$ magnetic field.
With $H\equiv 0$, \eqref{mhdeq} reduces to the equations of isentropic
fluid dynamics.

Equations \eqref{mhdeq} may be put in conservative form using identity
\begin{equation}
H\times \curl H= (1/2)\div (|H|^2I-2H^tH)^\trans + H\div H
\end{equation}
together with constraint \eqref{divfree} to express the second equation
as
\begin{equation}\label{cons2}
 \D_t(\rho  u) + \div(\rho u^tu)+ \na p + (1/2)\div (|H|^2I-2H^tH)^\trans=0.
\end{equation}
They may be put in symmetrizable (but no longer conservative) 
form by a further change, using identity
\begin{equation}
\curl (H \times u) = 
 (\div u) H+ (u\cdot \na)H -(\div H) u- (H\cdot \na)u
\end{equation}
together with constraint \eqref{divfree} to express the third equation as
\begin{equation}\label{symm3}
 \D_t H + 
 (\div u) H+ (u\cdot \na)H - (H\cdot \na)u=0.
\end{equation}

\subsection{Initial boundary value problem with constraint}
Constraint \eqref{divfree} is preserved by the flow, hence may be 
regarded as a constraint on the initial data.
Thus, a sufficient condition for well-posedness
under constraint \eqref{divfree} is that for some
form of the equations, rewritten modulo \eqref{divfree},
there hold well-posedness with respect to general initial data.
Under the assumptions of the previous section, i.e. noncharacteristicity,
symmetrizability, and the property that characteristics everywhere be
geometrically regular, totally nonglancing, or linearly splitting of
order two,
this is equivalent to satisfaction of the uniform Lopatinski
condition for the chosen representation of the equations.
This is sufficient for our present needs, and seems to
be the general approach followed in the literature.

To obtain sharp (i.e., necessary) conditions, one should
rather use the constraint together with the linearized equations
to make a pseudodifferential change of coordinates
eliminating one variable, arriving at an unconstrained
system in one fewer variables. 
Equivalently, as described in \cite{Ch-P}, one might
by a pseudodifferential choice of augmented variables
reduce the unconstrained second-order hyperbolic 
equation from which the system is derived 
to an unconstrained first-order system with modified boundary conditions.
However, the resulting analysis would be special to the system
considered, and not of the general form considered by Majda.

In any case, our analysis here concerns sufficient conditions 
for stability, and so we shall follow the first, simpler
approach and consider the equations: 
\begin{equation}
\label{mhdex}
\left\{ \begin{aligned}
 & \D_t \rho +  u \cdot \nabla \rho + \rho  \div u = 0, 
 \\
 &\rho ( \D_t u + u \cdot \nabla u)  + \na p + H \times \curl H = 0, 
 \\
 &  
 \D_t H + u \nabla H +  (\div u ) H  -  H \cdot \nabla u  = 0. 
 \end{aligned}\right.
\end{equation}
The Rankine--Hugoniot conditions are deduced form the 
  the conservative form
of the equations. Using notations $(y, x) \in \RR^2 \times \RR$ for spatial variables
as above, consider a
  shock front  
\begin{equation}
\label{shockf} 
 x  = \varphi (t, y)  
\end{equation} 
with space time normal  
\begin{equation}
\label{planef}
(- \sigma, n)  , \quad  \sigma = \D_t \varphi, \  n = ( - \D_{y_1} \varphi, \ldots, - \D_{y_{d-1}} \varphi, 1).   
\end{equation}
The jump conditions  read, with $u_n := n \cdot u$ and  $H_n = n \cdot H $, 
\begin{equation}\label{RH}
\left\{ \begin{aligned}
 & [\rho (u_n- \sigma) ] = 0, 
 \\
 &[\rho u (u_n - \sigma) ] +  n \left[  p + \frac{1}{2} |H|^2\right] - [ H_n H ]=0, 
 \\
 &  [(u_n - \sigma ) H ] - [ H_n u ] = 0, 
 \\
 & [ H_n ] = 0.    
\end{aligned}\right. 
\end{equation}
The last jump condition comes from the constraint equation \eqref{divfree}. Apparently
this system of 8 scalar equations is too large. However, projecting the third equation 
in the normal direction yields $\sigma [H_n ] = 0$ which is implied by the last equation. 
This shows that \eqref{RH} is made of 7 independent equations, as expected. 
Denoting by $u_{\rm tg}$  and $H_{\rm tg}$ the tangential part of $u$ and $H$, 
that is their orthogonal projection on $n^\perp$, \eqref{RH} is equivalent to 
\begin{equation}\label{RHm}
\left\{ \begin{aligned}
 & [\rho (u_n- \sigma) ] = 0, 
 \\
 &[\rho u (u_n - \sigma) ] +  n \left[  p + \frac{1}{2} |H|^2\right] - [ H_n H ]=0, 
 \\
 &  [(u_n - \sigma ) H_{\rm tg}  ] - [ H_n u_{\rm tg}  ] = 0,  \quad  [ H_n ] = 0.    
\end{aligned}\right. 
\end{equation}
The constraint \eqref{divfree} is preserved by the equations in the following sense. 

\begin{prop}
\label{prop71}
Suppose that $ U= (\rho, u, H)$ is  a smooth solution of \eqref{mhdex}  for 
$t \in [0, T]$ on both side 
of   the front  $x = \varphi(t, y)$  and satisfy the jump conditions 
\eqref{RHm}. Assume in addition that the traces of $u_n - \sigma$ on both side never vanish on the front and at least one among  the traces  $u_n^- - \sigma$ and $\sigma - u_n^+$ is positive. 
Then, if the initial  value of  $H$ satisfies 
$\div H_{ \vert t = 0}$, then $\div H = 0$ for all $t \in [0, T]$. 

\end{prop}
\begin{proof}  
The third equation in \eqref{mhdex} implies that on both side of the front: 
$$
\D_t H + \curl (H \times u) +  (\div H) u = 0. 
$$
Given a smooth function $v$ on both side 
of the front, denote by $\tilde v  $  the distribution equal to 
$v$  on both side of the front. The jump condition $[H_n ] = 0$ implies that 
$\div \tilde H = \widetilde { \div H}$. Moreover, the  jump conditions imply that 
$$
\D_t \tilde H + \curl (\tilde H \times \tilde u) =  (\D_t H + \curl (H \times u))\widetilde{} 
= - \widetilde {(\div H)} \tilde  u  = - (\div \tilde H) \tilde u. 
$$
 Therefore $ \tilde w := \div \tilde H$ is a bounded and piecewise continuous solution of 
$$
 \D_t  \tilde w +  \div (\tilde u \tilde w  )   = 0 , \quad   \tilde w_{\vert t = 0} = 0. 
$$ 
In particular, the jumps of $w = \div H$ satisfy 
$ [(u_n - \sigma) w ] = 0$, which could have been derived directly from the 
equations and the jumps conditions.  Note that there are no general uniqueness 
theorem for equations $\D_t w+ \div (a w) = 0$ with 
$a \in L^\infty$. But here $a$ is piecewise smooth ant the front is never characteristic
for this equation since we assumed that $u_n - \sigma \ne 0$. Moreover, the characteristics
are impinging the front at least on one side. 
Therefore, integrating along characteristics, this implies that $w = 0$ and thus that 
$\div \tilde H = 0$. 
\end{proof}

\subsection{The shock wave problem for MHD}

Thanks to Proposition \ref{prop71}, we consider the equations \eqref{mhdex} 
on both side of a front \eqref{shockf} , plus the transmission conditions \eqref{RHm}. 
 
Following the standard approach for linearized stability analysis of
fluid interfaces, (see \cite{Maj})
we carry the front location as an additional dependent variable,
and make the change of (independent) variables
$$
\tilde x :=x -\varphi (x',t),
$$
to fix the shock location in the new variables at $z\equiv 0$.  In the new coordinates, the problem becomes
\begin{equation}
\label{7.12}
\D_t U +
\sum^2_{j=1} A_j (U)   \D_{y_j} U      
  + \widetilde A_3 (U, d\varphi) \D_{\tilde x } U =0,  
 \qquad  \mathrm{on} \ \ \pm \tilde x > 0, 
\end{equation}
with 
$$
 \widetilde A_3 (U, d\varphi) =  A_3 (U) - \sum_{j=1}^2 (\D_{y_j} \varphi) \, A_j(U) 
 - ( \D_t \varphi )\,    \Id. 
 $$
The explicit form of the  matrices $A_j$ is  given in Appendix A, where it is also
verified that the system is hyperbolic symmetric in the sense of Friedrichs as soon as 
$c^2 = dp/ d\rho > 0$. 
The unknowns are $U = (\rho, u, H)$. The equations are supplemented with 
the transmission conditions \eqref{RHm} which read 
\begin{equation}
\label{7.13}
\D_t \varphi [F_0(U) ] + \sum_{j=1}^2 \D_{y_j} \varphi [F_j(U)] = [F_3(U)]
 \qquad  \mathrm{on} \ \  \tilde x = 0. 
\end{equation}
The stability analysis yields to consider the linearized equations
\begin{equation}
\label{7.14}
\D_t \dot U +
\sum^2_{j=1} A_j (U)   \D_{y_j} \dot U      
  + \widetilde A_3 (U, d\varphi) \D_{\tilde x } \dot U = \dot f ,  
 \qquad  \mathrm{on} \ \ \pm \tilde x > 0, 
\end{equation}
\begin{equation}
\label{7.15}
\D_t \dot \varphi [F_0(U) ] + \sum_{j=1}^2 \D_{y_j}  \dot \varphi [F_j(U)] = 
[\widetilde F' _3(U, d\varphi) \dot U ] + \dot g   \qquad  \mathrm{on} \ \  \tilde x = 0, 
\end{equation}
with 
$$
\widetilde F' _3(U, d\varphi) = F'_3 (U) - \sum_{j=1}^2 (\D_{y_j} \varphi) \, F'_j(U) 
 - ( \D_t \varphi )\,    F'_0(U). 
 $$

The Majda-Lopatinski determinant is computed as follows. One consider the linearized equations
\eqref{7.14} \eqref{7.15} for $U^-$ and $U^+$ constant and $d \varphi $ constant. 
This yields constant coefficients system, depending on parameters 
$p = (U^-, U^+, d \varphi)$  to which we apply the analysis developed in the previous 
sections. 
Equations \eqref{7.14} \eqref{7.15} comprise a slightly nonstandard, 
constant--coefficient initial-boundary-value problem on the half-open
space  $\tilde x >0$, in the ``doubled'' variables
$$
\tilde  U^+(z):= \dot U^+ (z), \quad \tilde  U^-(z):=\dot U^- (-z),
$$
the new feature being the appearance of front variable $\dot \varphi$.
As pointed out by Majda \cite{Maj}, this may be 
handled within the standard framework, thanks to the cascaded form
of the equations.

For, assuming that the shock front is noncharacteristic, 
which follows from Lax' shock conditions, 
taking the Fourier transform in the transverse coordinates $y $
and the Laplace transform in $t$, we obtain 
a system of ODE in variables $\hat U^+,\hat U^-$,
\begin{equation}
\label{7.16}
 \D_{\tilde x} \hat U^\pm  + G^\pm (p, \zeta) \hat U^\pm  = \hat f^\pm  , 
\quad \mathrm{on} \ \ \tilde x  \gtrless 0,
\end{equation}
with $G^\pm $ as in \eqref{nn32},  and  boundary condition  , 
\begin{equation}
\label{7.17}
\hat \varphi  X(p, \zeta)  =  [\tilde F'_3 \, \hat U] +  \hat g,
\quad \mathrm{on} \ \ \tilde x  =  0,
\end{equation}
 with 
 $$
 X(p, \zeta) = (i \tau - \gamma) [F_0(U) ] + 
 \sum_{j=1}^2 \eta_j [F_j(U)]. 
 $$
 Noting that $\hat \varphi $ does not appear in the first equation,
we may convert it to a separate, standard boundary-value problem by
introducing a new, decoupled boundary condition of rank $ N-1= 6 $,
obtained by taking the inner product of \eqref{7.17} with  a basis of unit vectors
in  $ X(p, \zeta)^\perp$. This reduces the boundary conditions to a system
of the form 
\begin{equation}
\label{7.18}
M(p, \zeta) \hat U =  \hat g' , \qquad  \hat \varphi = \Phi(p, \zeta) \hat U + \hat g_1, 
\end{equation}
indicating how the boundary value problem analysis developed in the previous
sections applies to the present  framework. 

More precisely, one introduces for $\gamma > 0$, the spaces 
$\EE^\pm (p, \zeta)$ of initial conditions 
   corresponding
to decaying solutions of \eqref{7.16} at $ \pm \infty$.  
They are generated by the generalized eigenspaces of 
$G^\pm (p, \zeta)$ located in $\{ \pm \im \mu <  0 \}$. 
Lax' condition  implies that 
\begin{equation}
\label{Lax}
\dim \EE^- (p, \zeta) + \dim \EE^+(p, \zeta) = N-1 = 6. 
\end{equation}
The stability condition  for the boundary condition \eqref{7.17} states that the homogeneous
equation  with $\hat g = 0$ has no solution 
$(\hat \varphi, \hat U^-, \hat U^+) $ in $\CC \times \EE^- \times \EE^+$. 
In accordance with the general Definition \ref{lop}, 
we obtain Majda's Lopatinski
determinant for the full shock problem, 
\begin{equation}\label{7.20}
 D(p, \zeta)  =  \det
\begin{pmatrix}  \widetilde F'_3{}^- \EE^- ,  &    \widetilde F'_3{}^+ \EE^+ , 
&  \CC X 
\end{pmatrix},
\end{equation}
or 
\begin{equation}\label{shocklop}
 D(p, \zeta)  =  \det
\begin{pmatrix}  r^-_1, \ldots, r^-_{p-1}  ,  &    r^-_{p+1}, \ldots, r^-_{7}  , 
&    X 
\end{pmatrix},
\end{equation}
where $r_j^-$ and $r_j^+$ are bases for $\EE^-$ and $\EE^+$. 
  
  \subsection{Applications of Theorem \ref{2main} } 
  
  Following the analysis presented in the  previous subsection, 
  we consider the linearized equations \eqref{7.14} \eqref{7.15} around 
  a planar shock. Without loss of generality, we can assume that 
  the spatial normal to the front is 
  $n = (0, 0, 1)$ so that the normal matrix 
  is $\widetilde A_3 = A_3   - \sigma \Id $. 
  The seven  eigenvalues of 
  $A(\eta, \xi) = \eta_1 A_1 + \eta_2 A_2 + \xi \widetilde A_3$  are 
  \begin{equation}
\label{7.21}
\left\{\begin{aligned} 
  \lambda_0  & = \eta_1 u_1 + \eta_2 u_2 + ( u_3 - \sigma)\xi 
\\
 \lambda_{\pm 1} & =  \lambda_0 \pm c_s  (\vert \eta \vert^2 + \xi^2)^{1/2} 
\\
 \lambda_{\pm 2} & = \lambda_0  \pm  (\eta_1 H_1  + \eta_2 H_2 + \xi H_3 )/ \sqrt \rho
\\
  \lambda_{\pm 3} & = \lambda_0 \pm c_f  (\vert \eta \vert^2 + \xi^2)^{1/2} 
\end{aligned}\right.
\end{equation}
  with 
$$ 
  \begin{aligned}
    c_f^2  := \frac{1}{2} 
  \Big( c^2 + h^2  ) + \sqrt{ (c^2 - h^2)^2 + 4 b^2 c^2  } \Big)
  \\
    c_s^2  := \frac{1}{2} 
  \Big( c^2 + h^2  ) -  \sqrt{ (c^2 - h^2)^2 +  4 b^2  c^2   } \Big)
\end{aligned}
$$
with $c^2 = p'(\rho)$ the sound speed, $h^2   = \vert H \vert^2 / \rho$, 
$b^2 = \vert (\hat \eta, \hat \xi) \times H \vert^2 / \rho$ and 
$ (\hat \eta, \hat \xi) = (\eta, \xi) / \vert  (\eta,\xi) \vert$ 
(see Appendix A for the computations). 
  The first eigenvalue corresponds to the transport of the 
  constraint. It can be decoupled from the system :  there is a smooth one dimensional subspace, 
  $\EE_0$ such that $A(\eta, \xi) = \lambda_0$  on this space 
  and $\EE_0^\perp$ is stable for $A(\eta, \xi)$. 
  The other eigenvalues are in general simple.  
  
  The boundary $\{ \tilde x = 0 \} $ is noncharacteristic if and only if 
  \begin{equation}
\label{7.23}
u_3 - \sigma \notin \big\{ 0, \  \pm H_3 / \sqrt \rho, \ \pm c_s(n), \ \pm c_f(n)
\big\}  
\end{equation}
  where $c_s(n)$ and $c_f(s)$ are the slow and fast speed computed in the 
  normal direction $n = (0, 0, 1)$. 
  The next lemma is proved in Appendix A.
  
 \begin{lem}
 \label{lem72}
  Assume that $ 0 < \vert H \vert^2  \ne \rho c^2$. Consider 
  $(\eta, \xi)$ in the unit sphere $S^2$. 
  
 \quad i)  When $  (\eta, \xi) \cdot H \ne 0$ and 
 $  (\eta, \xi) \times H \ne 0$, the eigenvalues are simple. 
 
 \quad ii) On the manifold $  (\eta, \xi) \cdot H =  0$  the eigenvalues $\lambda_{\pm 3}$ 
 are simple and the multiple eigenvalue $\lambda_0 = \lambda_{\pm 1} = \lambda_{\pm 2}$
 is geometrically regular.  
 
 \quad iii) On the manifold $  (\eta, \xi) \times  H =  0$, $\lambda_0$ is simple. 
 When $\vert H \vert^2 < \rho c^2$  [resp. $\vert H \vert^2 > \rho c^2$ ], 
 $\lambda_{\pm 3}$ [resp. $\lambda_{\pm 1}$ are simple, the other eigenvalues  
 $\lambda_{+2}$ and $\lambda_{-2}$ are double, algebraically regular, 
 and totally nonglancing provided that  $u_3 - \sigma \ne \pm H_3/ \sqrt \rho$. 
 \end{lem}
  
   This shows that the assumptions of Theorem \ref{2main} 
 are satisfied as soon as 
 the hyperbolicity condition $c^2 := dp/d\rho>0$  holds, 
  and  $H\ne 0$ with length $\vert H \vert \ne  c \sqrt \rho $
  and condition \eqref{7.23} is satisfied. 
   
\begin{cor}
\label{mhdlop}  
 Under the generically satisfied conditions $0 <\vert H^\pm \vert^2 \ne \rho^\pm (c^\pm)^2$,
  the uniform Lopatinski condition implies
linearized and nonlinear stability of noncharacteristic Lax-type MHD shocks.
\end{cor}

\begin{proof}
Lemma \ref{lem72} and Theorem \ref{2main} imply the existence of smooth
Kreiss  symmetrizers. The  stability follows through the standard
shock stability framework of \cite{Maj, MetKochel} (see the useful Remark \ref{rem59x}).
\end{proof}

\subsection{The $H \to 0$ limit}

Using Corollary \ref{mhdlop}, we easily recover the result of Blokhin
and Trakhinin in the isentropic case, by a simple perturbation argument.
When $H = 0$, the system \eqref{mhdex} reduces to isentropic Euler's equations
and \eqref{RHm} to the corresponding Rankine Hugoniot condition. 
Denote by $D_0( p_0, \zeta)$ the Majda-Lopatinksi determinant 
of this problem, where the parameters are now 
$p_0 = (U_0^-, U_0^+, d\varphi) $ with $U_0 = (\rho, u)$. 
Note that  because we have everywhere
either geometric regularity or total nonglancing, the stable subspace $\EE_-$,
and thus the Lopatinski 
determinants have   continuous extensions to $\zeta \in \overline S^3_+$
(see Theorem \ref{2main}).

When $H = 0$, the eigenvalues are 
\begin{equation}
\begin{aligned}
\label{7.24}
&\lambda_0 = \lambda_{\pm1} = \lambda_{\pm 2} = \eta_1 u_1 + \eta_2 u_2 + \xi (u_3- \sigma), 
\\
& \lambda_{\pm 3} = \lambda_0 \pm c \vert (\eta, \xi) \vert. 
\end{aligned}
\end{equation} 
The boundary is noncharacteristic when
  \begin{equation}
  \label{7.25}
  u_3 - \sigma \notin \{  - c,  0, + c \}. 
  \end{equation}
  
  \begin{lem}
  \label{lem74}
 Suppose that $\underline U = (\underline \rho, \underline u, 0)$
  and $\sigma$ satisfy  \eqref{A26} holds. Then the multiple eigenvalue $\lambda_0$ is totally 
  nonglancing at $(\underline U, \eta, \xi)$ for all 
  $(\eta, \xi) \ne 0$ and the eigenvalues $\lambda_{\pm 3}$ are simple. 
  \end{lem} 
 
 The proof is given in the Appendix A. 
 By Theorem \ref{2main}, this implies that the negative spaces depend continuously 
 of $H$ for $H$ close to $0$. Therefore, the Lopatinski determinant 
$D(p, \zeta)$ is a continuous function of 
$\zeta \in \overline S^3_+ $ and $p= (U^-, U^+, d \varphi) $, with 
$U^\pm = (\rho^\pm, u^\pm, H^\pm)$, when  
$U^\pm$ satisfy \eqref{7.25} and $H^\pm$ are small. 
Therefore:

\begin{cor}
\label{cor75}
Suppose that $p_0$ is a Lax shock for 
the insentropic Euler's equation. Then, as $(H^-, H^+) \to 0$, there holds 
$D( p, \zeta) \to D_0(p_0, \zeta)$, uniformly for 
$\zeta \in S^3$, with $\gamma \ge 0$. 
\end{cor}

\begin{cor}\label{pert}
In the $H\to 0$ limit, Lax-type MHD shocks approaching a noncharacteristic
limiting fluid-dynamical shock satisfy
Assumption \ref{ass41}, so that the Lopatinski condition reduces
to \eqref{eq410x}; moreover, they
satisfy the Lopatinski condition if and only if it is satisfied
by the limiting fluid-dynamical shock (i.e., satisfies the physical stability
conditions of Erpenbeck--Majda \cite{Er, Maj}, in which case
the MHD shocks are linearly and nonlinearly stable.
In particular, for an ideal gas equation of state, they are
always stable in the $H\to 0$ limit.
\end{cor}

Likewise, using Corollary \ref{mhdlop}, we may immediately convert the
results of \cite{GK}, \cite{BT.1}--\cite{BT.4}, etc. on satisfaction of the
Lopatinski condition to full results of linearized and nonlinear stability.
More generally, Corollary \ref{mhdlop}
gives a useful framework for the systematic study of MHD shock
stability and other variable-multiplicity PDE problems through the
investigation of the linear-algebraic Lopatinski condition.

\begin{rem}
\label{uniform}
\textup{
Corollary \ref{pert} guarantees stability for sufficiently small
but nonzero values $H^\pm$.  But since the system is symmetric 
 hyperbolic for all $H$, a direct application of Theorem \ref{2main} 
provides us with Kreiss symmetrizers which are smooth in 
$ p = (U^-, U^+, d \varphi)$ for  $|H^\pm| $ small. 
Therefore, we have   maximal    stability estimates 
for the linearized equations  which are uniform in $H^\pm$ 
for $H^\pm$ small.  
 This implies  uniform nonlinear stability, since the time of existence 
 depends only on the estimates (see \cite{Maj}, \cite{MetKochel}).
 This is likely to imply the continuity of the shock solution as $H \to 0$. 
}
\end{rem}

\begin{rem}
\label{smallamp}
\textup{
Calculations as in \cite{MetKochel, Zhandbook} show
that Lax-type MHD shocks
are stable in the small-amplitude limit $[U]\to 0$,
again, for $H_\pm \ne 0$.
}
\end{rem}
 

%
\appendix

\section{Appendix A. The symbolic structure of MHD}

 \subsection{Multiple eigenvalues}
 
The first order  term  
of the linearized equations of \eqref{mhdex} about $(u, H)$ is
\begin{equation}\label{linearized}
\left\{ \begin{aligned}
 & D_t  \dot \rho + \rho \div \dot u    
 \\
 &  D_t \dot u +   \rho^{-1} c^2 \na \dot \rho  + \rho^{-1} H \times \curl \dot H  
 \\
 &  
 D_t \dot H + (\div \dot u) H   - H \cdot \na  \dot u  
 \end{aligned}\right.
\end{equation}
with $D_t = \D_t + u \cdot \na$ and $c^2  = dp/d\rho$. 
This system is hyperbolic symmetric
with symmetrizer 
$S = \bdiag( c^2, \rho \Id, \Id ) $.  The associated symbol is
\begin{equation}
\label{Asymbol}
\left\{ \begin{aligned}
 & \tilde \tau   \dot \rho + \rho (\xi \cdot \dot u)   
 \\
 &  \tilde \tau  \dot u +  \rho^{-1} c^2  \dot \rho \xi  + \rho^{-1} H \times (\xi \times \dot  H  )
 \\
 &  
\tilde  \tau \dot H +  ( \xi \cdot \dot u) H  - (H \cdot \xi)  \dot u   
 \end{aligned}\right.
\end{equation}
with $\tilde \tau  = \tau + u \cdot \xi$.  We use here the notation
$\xi = (\xi_1, \xi_2, \xi_3)$ for the spatial  frequencies and 
$$
\xi = \vert \xi \vert\, \hat\xi \,, \quad  u_\parallel = \hat \xi \cdot u \,, \quad 
u_\perp = u - u_\pr \hat \xi = - \hat \xi \times (\hat \xi \times u)\,. 
$$
 
 We write \eqref{Asymbol} in the general form 
 $\tau \Id + A (U, \xi)$ with parameters $U = (\rho, u, H)$. 
 The eigenvalue equation $A(U, \xi) \dot U = \lambda \dot U$ 
reads 
\begin{equation}
\label{Achar}
\left\{ \begin{aligned}
 & \tilde\lambda \dot \rho =  \rho   \dot u_\pr   , 
 \\
 &\rho \tilde\lambda \dot u_\pr =   c^2  \dot \rho   + H_\perp \cdot \dot H_\perp  , 
  \\
  &\rho \tilde\lambda \dot u_\perp =  - H_\pr \dot H_\perp, 
  \\
   &  
 \tilde\lambda \dot H_\perp  =   \dot u_\pr  H_\perp   - H_\pr  \dot u_\perp,
\\
 & \tilde \lambda \dot H_\pr = 0 ,    
 \end{aligned}\right.
\end{equation}
with $\tilde \lambda = \lambda - (u \dot \xi)$. 
The last condition decouples. 
On the space 
\begin{equation}
\label{Ap2}
\EE_0 (\xi) = \big\{ \dot \rho = 0, \, \dot u = 0, \ \dot H_{\perp} = 0 \big \},
\end{equation}
$A$  is equal to 
$\lambda_0  :  = u\cdot \xi$.  
{}From now on we work on $\EE_0 ^\perp = \{ \dot H_\pr = 0 \}$
which is invariant by $A(p, \xi)$. 
 
Consider $ v = H/ \sqrt \rho $,  $\dot v =  \dot H/ \sqrt \rho $ and $\dot \sigma = \dot \rho / \rho $. 
The characteristic system reads: 
\begin{equation}
\label{Alambdachar}
\left\{ \begin{aligned}
 & \tilde\lambda \dot \sigma  =     \dot u_\pr   
 \\
 &  \tilde\lambda \dot u_\pr =   c^2  \dot \sigma   + v_\perp \cdot \dot v_\perp  
  \\
  &  \tilde\lambda \dot u_\perp =  - v_\pr \dot v_\perp
  \\
 &  
 \tilde\lambda \dot v_\perp  =      \dot u_\pr  v_\perp   - v_\pr  \dot u_\perp
 \end{aligned}\right.
\end{equation}
Take a basis of $\xi^\perp$ such that 
$v_\perp = (b, 0)$ and let $a = v_\pr$. 
In such a basis, the matrix of the system reads
\begin{equation}
\label{calA}
   \tilde\lambda  - \tilde\cA 
      := \left[ \begin{matrix}
 \tilde\lambda  & -1 & 0 &0 & 0 &0 
 \\  - c^2  &   \tilde\lambda & 0  & 0 & -b  &0
\\ 0 & 0 & \tilde\lambda &  0 &  a &  0
\\0 & 0 & 0 & \tilde\lambda    &  0 & a 
 \\
 0 & -b  & a &0 & \tilde\lambda &0 
 \\
 0 & 0 & 0 & a & 0 & \tilde\lambda  
\end{matrix} \right]
\end{equation}
The characteristic roots satisfy 
\begin{equation}
\label{Aroots}
(  \tilde\lambda^2 -  a^2  ) \big(
  (  \tilde\lambda^2 -  a^2  ) (\tilde\lambda^2 - c^2)  - \tilde\lambda^2  b^2 \big) = 0 \,. 
\end{equation}
Thus, either 
\begin{eqnarray}
\label{A8} &  & \tilde\lambda^2 = a^2 
\\ \label{A9}
  & & 
 \tilde\lambda^2 = c_f^2  := \frac{1}{2} 
  \Big( c^2 + h^2  ) + \sqrt{ (c^2 - h^2)^2 + 4 b^2 c^2  } \Big)
  \\ \label{A10}
   & &  \tilde\lambda^2 = c_s^2  := \frac{1}{2} 
  \Big( c^2 + h^2  ) -  \sqrt{ (c^2 - h^2)^2 +  4 b^2  c^2   } \Big)
\end{eqnarray}
with $h^2 = a^2 + b^2 = \vert H \vert^2 / \rho$. 

With  $P(X) = (X - a^2) (X-c^2) - b^2 X  $, $\{ P \le 0\} = [ c_s^2, c_f^2]$ and 
$P (X) \le 0$ for $X \in [\min (a^2, c^2), \max(a^2, c^2) ]$. 
Thus, 
\begin{eqnarray}
 & & c_f^2 \ge \max (a^2, c^2)  \ge a^2  \\
&& c_s^2 \le \min (a^2, c^2)  \le a^2  
\end{eqnarray}
\bigbreak

\noindent{\bf 1. }  The case $v_\perp \ne 0$  i.e. $ w = \hat \xi \times v \ne 0$.
Thus, the basis such that  \eqref{calA}  holds is smooth in $\xi$. 
In this basis, $ w = (0, b)$, $b = \vert v_\perp \vert > 0$. 

{\bf  1.1 }  The spaces 
$$
\EE_{\pm } (\hat \xi) = \{ \dot \sigma = 0,  \dot u_\pr = 0,\  \dot v_\perp \in
 \CC ( \hat \xi \times v) ,  \ 
\dot u_\perp =  \mp  \dot v_\perp \} 
$$
are  invariant for $\tilde\cA$ and 
 \begin{equation}
 \tilde\cA   =  \pm a  \quad \mathrm{on} \quad \EE_{\pm }\,. 
 \end{equation} 

\bigbreak
{\bf 1.2 }   In $(\EE_+ \oplus \EE_-)^\perp$, which is invariant,  the matrix of $\tilde\cA$ is 
\begin{equation}
 \tilde\cA_0 
      := \left[ \begin{matrix}
0  & 1 & 0 &0 
 \\   c^2  &   0 & 0  & -b  
\\ 0 & 0 & 0  &  -  a 
 \\
 0 & -b  & a &0 & 
\end{matrix} \right]
\end{equation}

Since $P(c^2) = - b^2 c^2 < 0$, there holds
$c_s^2 < c^2 < c_f^2$. 

\qquad {\bf 1.2.1 ) } Suppose that  $a \ne 0$. Then, $P(a^2 ) = - a^2 c^2 < 0$ and 
$c_s^2 < a^2 < c_f^2$.  Thus, all the eigenvalues are simple. 
Moreover, $c_s^2 c_f^2 = a^2 c^2$ and $c_s^2 > 0$. 
The space 
$$
\FF_{\tilde\lambda}  = \Big\{  \dot \sigma =  \tilde\lambda \dot u_\pr , \ 
\dot u_\pr = \frac{\tilde\lambda b  v_1 } { \tilde\lambda^2 - c^2 } , \ 
\dot u_1 =  \frac{ - a v  } { \tilde\lambda  } , \ v_1 \in    \CC  \Big\}
$$
is  an eigenspace associated to the eigenvalue $\tilde\lambda$ 
when $\tilde\lambda = \pm c_f $ and $\tilde\lambda = \pm c_s$. 
 Here $\dot u_1$ and $\dot v_1$ denote the first component of $\dot u_\perp$ 
 and $\dot v_\perp$ respectively in the basis $(v_\perp, w)$. 
 
\medbreak
\qquad {\bf 1.2.1 ) } Suppose that  $a$ is close to $0$. 
Since $c_f^2 > c^2 > 0$, the spaces $\FF_{\pm c_f}$ are still eigenspaces 
associated to the eigenvalues  $\tilde\lambda = \pm c_f $.

By direct computations:
\begin{equation*}
c_s^2 = \frac{c^2 a^2}{ c^2 + h^2} + O( a^4)\,. 
\end{equation*}
Therefore, 
\begin{equation*}
\frac {c^2_s}{a^2}  \to  \frac{c^2 }{ c^2 + h^2}  > 0  \quad  \mathrm{as } \ a \to 0\,. ,. 
\end{equation*}
 Therefore, consider $\tilde c_s =  \frac{a}{\vert a \vert} c_s$ is an analytic 
 function of $a$ (and $b \ne 0$)  near $a= 0$  and 
 $$
 \widetilde \FF_{\pm , s} (a, b) =  \FF_{\pm   \tilde  c_s } 
 $$
 are analytic determinations of Eigenspaces, associated to the eigenvalues 
 $\pm \tilde c_s$. Moreover, the values at $a=0$ are 
 $$
 \widetilde \FF_{\pm , s} (0, b) = \big\{  \dot \sigma =  0  , \ 
\dot u_\pr =  0  , \ 
\dot u_\perp =  \frac{ \mp   c  v  } { \sqrt{c^2 + b^2}  }  \big\}
 $$
 and $\widetilde \FF_{+, s} \cap \widetilde \FF_{-, s} = \{0 \}$, 
 thus we still have an analytic diagonalization of $\tilde\cA_0$.

\bigbreak 

\bigbreak 
\noindent{\bf 2. }  Suppose now that $b$ is close to zero. 
At  $ b = 0$, the eigenvalues of $\tilde\cA$ are $\pm c $ (simple) and $\pm h$  (double). 
{\bf Assume that $c^2 \ne h^2$}. 
Note that when $b = 0$, then $\vert a \vert = h$ and 
$$
\begin{matrix}
\mathrm{when} \ c^2 > h^2  \ : &   c_f = c , &  c_s = h ,
\\
\mathrm{when} \ c^2 <  h^2\ :  &   c_f = h ,  &  c_s = c . 
\end{matrix}
$$

{\bf 2.1 }  The eigenvalues close to $ \pm c$  remain simple. 

\medbreak

{\bf 2.2 }  We look for the eigenvalues close to $h$. 
The characteristic equation implies that 
\begin{equation}
\label{Arootdescription}
\left\{ \begin{aligned}
&   c^2  \dot \sigma    = \tilde\lambda \dot u_\pr - v_\perp \cdot \dot v_\perp  
\\
 &  (\tilde\lambda^2 - c^2)  \dot u_\pr =   \tilde\lambda  v_\perp \cdot \dot v_\perp   
\end{aligned}\right. 
\end{equation} 
Eliminating $\dot u_\pr$, we are left with the $4 \times 4$ system
in $\xi^\perp \times \xi^\perp$:
\begin{equation}
\label{A4x4}
\left\{ \begin{aligned}
  & \tilde\lambda \dot u_\perp =  -a  \dot v_\perp
  \\
 &  
 \tilde\lambda \dot v_\perp  =      -  a   \dot u_\perp +
 \frac{\tilde\lambda}{\tilde\lambda^2 - c^2}   (v_\perp \otimes v_\perp) \dot v_\perp. 
 \end{aligned}\right.
\end{equation}
  Thus, 
  \begin{equation*}
(  \tilde\lambda^2 - a^2)  \dot v_\perp  =       
 \frac{\tilde\lambda^2}{\tilde\lambda^2 - c^2}   (v_\perp \otimes v_\perp) \dot v_\perp. 
\end{equation*}
Recall  that $\vert v_\perp \vert = b$ is small. 
We recover  4  {\sl  smooth  eigenvalues } 
\begin{equation}
\label{A17}
 \pm a  \,, \qquad \pm \sqrt  {a^2 + O (b^2)} = \pm ( a + O(b^2))  \,. 
\end{equation}
(remember that $a = \pm h + O(b^2)$. 
However, the eigenspaces  are not smooth in $v$, since they are 
$\RR v_\perp$ and $ \RR \hat \xi \times v_\perp$ and have no limit at 
$v_\perp \to 0$.

\bigbreak
Summing up, we have proved the following. 

\begin{lem}
\label{lemA1} Assume that $ c^2 = dp/d\rho < 0$.  The eigenvalues of $A(U, \xi)$ are 
 \begin{equation}
\label{A20}
\left\{\begin{aligned} 
  \lambda_0  & =  \xi \cdot u 
\\
 \lambda_{\pm 1} & =  \lambda_0 \pm c_s (\hat \xi)  \vert  \xi  \vert 
\\
 \lambda_{\pm 2} & = \lambda_0  \pm  (\xi \cdot H)/ \sqrt \rho
\\
  \lambda_{\pm 3} & = \lambda_0 \pm c_f  (\hat \xi  )  \vert \xi \vert 
\end{aligned}\right.
\end{equation}
with  $\hat \xi = \xi / \vert \xi \vert$ and 
\begin{eqnarray}
  &&c_f^2(\hat \xi)  := \frac{1}{2} 
  \Big( c^2 + h^2  ) + \sqrt{ (c^2 - h^2)^2 + 4 b^2 c^2  } \Big)
  \\
   & &   c_s^2 (\hat \xi)   := \frac{1}{2} 
  \Big( c^2 + h^2  ) -  \sqrt{ (c^2 - h^2)^2 +  4 b^2  c^2   } \Big)
\end{eqnarray}
where  $h^2 =   \vert H \vert^2 / \rho $, $b^2 = \vert \hat \xi \times H \vert^2/ \rho $.

\end{lem}

\begin{lem}
\label{lemA2}
  Assume that $ 0 < \vert H \vert^2  \ne \rho c^2$ 
where $c^2=dp/d\rho > 0 $.

\quad i)  When $ \xi  \cdot H \ne 0$ and 
 $   \xi \times H \ne 0$, the eigenvalues of $A(U, \xi)$  are simple.

 \quad ii) On the manifold $   \xi  \cdot H =  0$, $\xi \ne 0$,   the eigenvalues $\lambda_{\pm 3}$ 
 are simple and the multiple eigenvalue $\lambda_0 = \lambda_{\pm 1} = \lambda_{\pm 2}$
 is geometrically regular.  
 
 \quad iii) On the manifold $  \xi  \times  H =  0$, $\xi \ne 0$, $\lambda_0$ is simple. 
 When $\vert H \vert^2 < \rho c^2$  [resp. $\vert H \vert^2 > \rho c^2$ ], 
 $\lambda_{\pm 3}$ [resp. $\lambda_{\pm 1}$]  are simple;  
 $\lambda_{+2} \ne \lambda_{-2}$ are double, equal to $\lambda_{\pm 1}$
 [resp. $\lambda_{\pm 3}$ ] depending on the sign of $ \xi \cdot H$ , algebraically regular, 
  but not geometrically regular.  
  
\end{lem}

 
 \subsection{Glancing conditions for MHD}
 
 Write the space variables $(y, x)$, $  y \in \RR^2$, $x \in \RR$ and 
 consider the planar front  $\{x  = \sigma t \}$.   As indicated in Section 7, 
 we perform the change of variables $\tilde x = x - \sigma t$, obtaining 
 the linearized system \eqref{7.14} with reads
 \begin{equation}
 \label{} 
 \D_t - \sigma \D_{\tilde x} +  A (U, \D_y, \D_x) = \D_t + A(\tilde U, \D_y, \D_x) 
 \end{equation}
 with 
 $\tilde U = (\rho, \tilde u, H)$, $\tilde u = (u_1, u_2, u_3 - \sigma)$. 
 This shows that the eigenvalues of the symbol are 
 \begin{equation}
 \tilde \lambda(U, \eta, \xi) = \lambda(\tilde U, \eta, \xi)
 \end{equation}
 where $(\eta, \xi)$ denote the tangential Fourier frequencies dual to 
 $(y, tilde x)$. Therefore, Lemma \ref{lemA1} implies that the eigenvalues are given
 by \eqref{7.21}. 

In particular, the eigenvalues of the boundary matrix are
\begin{equation}
u_3 - \sigma, \quad
u_3 - \sigma  \pm H_3 / \sqrt \rho  \,, \quad u_3 - \sigma 
  \pm c_f(n ) \,, \quad u_3 - \sigma  \pm c_s(n) \,. 
\end{equation}
with $n = (0,0, 1)$.  The boundary is noncharacteristic if 
  they are different from zero, that we assume from now on. 
  The next result finishes the proof of Lemma \ref{lem72}. 

\begin{lem}
Suppose that $ 0 < \vert H \vert^2  \ne \rho c^2$ and the front in noncharacteristic.  
On the manifold $ (\eta, \xi)  \times  H =  0$, $(\eta, \xi) \ne 0$,   the double eigenvalues  
 $\lambda_{+2}$ and $ \lambda_{-2}$  are totally nonglancing. 
\end{lem}

\begin{proof}  Consider 
\begin{equation}
\label{A24}
\lambda_2  =  \eta_1 (u_1 + v_1) + \eta_2 (u_2 + v_2) + \xi (u_3 + v_3 - \sigma),  
\end{equation}
 with $v = H/ \sqrt \rho$. We have shown  at \eqref{A17} that the eigenvalue is algebraically  regular, and the nearby eigenvalue 
 is 
$$
\lambda_2 + O( \vert \xi \times H \vert^2). 
$$
 It is tangent to $\lambda_2$ to second order.
   Therefore, the nonglancing condition reduces to 
$$
\frac{\D \lambda_2}{\D \xi } = u_3 + v_3 - \sigma   \ne 0\,,
$$
 that  is automatically satisfied when the front is noncharacteristic. 
  \end{proof}
 
 We now pass to the proof of Lemma~\ref{lem74}.   
When $H = 0$, the eigenvalues are 
\begin{equation}
\begin{aligned}
\label{A25}
&\lambda_0 = \lambda_{\pm1} = \lambda_{\pm 2} = \eta_1 u_1 + \eta_2 u_2 + \xi (u_3- \sigma), 
\\
& \lambda_{\pm 3} = \lambda_0 \pm c \vert (\eta, \xi) \vert. 
\end{aligned}
\end{equation} 
The boundary is noncharacteristic when
  \begin{equation}
  \label{A26}
  u_3 - \sigma \notin \{  - c,  0, + c \}. 
  \end{equation}
  
  \begin{lem}
  Suppose that $\underline U = (\underline \rho, \underline u, 0)$
   and $\sigma$ satisfy  \eqref{A26} holds. Then the multiple eigenvalue $\lambda_0$ is totally 
  nonglancing at $(\underline U, \eta, \xi)$ for all 
  $(\eta, \xi) \ne 0$. 
  \end{lem} 
 
\begin{proof}
Note that $\lambda_0$ has constant multiplicity  5  in $(\eta, \xi)$ for 
$H = 0$, and is linear in $(\eta, \xi)$. 
Therefore, the tangent system at $(\underline U, \eta, \xi)$ is $\lambda_0 \Id$  in dimension 5. 
In particular, with notations as in Section 3,   the boundary matrix of this system is $\underline A_3 = (u_3 - \sigma) \Id$, implying that the eigenvalue is totally nonglancing.
\end{proof}


\section{Appendix B. Maxwell's equations in a bi-axial crystal}

We give here an example  where  geometric regularity fails.
The consequence is that the tangent operator
  is no longer a transport operator
along a single direction, but rather a system of hyperbolic equations
of wave equation type.
This example permits the description of
the phenomenon of \textit{conical refraction} 
(cf \cite{Lud, Taylor, BW}) 

Maxwell's equations in the absence of exterior charge may be written as
\begin{equation} 
\label{Maxbiax}
\left\{\begin{aligned}
{}& \partial _t B + \rot E = 0\, , \quad \div B = 0\, , 
\\
& \partial _t D - \rot H   = 0\, ,  \quad  \div D  = 0\, , 
\end{aligned}\right. 
\end{equation} 
with
\begin{equation} 
D = \cE E \, , \quad B = \mu H\, 
\end{equation} 
where  $\cE$ and $\mu$ are $3 \times 3$ positive definite matrices.
In the case of a bi-axial crystal, $\mu$ is scalar and 
$\cE$ has three distinct eigenvalues.
We change variables so that $\mu = 1$ and choose
coordinate axes so that $\cE$ is diagonal:
\begin{equation} 
\cE^{-1} = \left( \begin{array}{ccc}
\alpha^{}_1 & 0 & 0 \\0 & \alpha_2^{} & 0 \\0&0&\alpha_3^{}
\end{array}\right)
\end{equation}  
with $\alpha_1 > \alpha_2 > \alpha_3$. 
Ignoring the divergence conditions, the characteristic
equation and the polarization conditions are obtained
as solutions of system
\begin{equation} 
\label{symbMaxbiax}
L(\tau, \xi)  \left(\begin{array}{c} B \\E 
\end{array}\right) := 
\left(\begin{array}{lll}
  \tau  B &+ &\xi  \times  E     
\\
 \tau      E  &-  &\cE^{ -1} (\xi \times   B)    
\end{array}\right)  \, =\, 0\, . 
\end{equation} 
For $\eta \ne 0$, 
$\tau = 0$ is a double eigenvalue, with eigenspace generated by
$(\eta, 0)$ and $(0, \eta)$.  These modes are 
incompatible with the divergence conditions.
The nonzero eigenvalues are given as solutions of
$$ 
E = \cE^{-1} (\frac{\xi }{\tau}\times  B) \, , \quad 
(\tau^2 + \Omega(\xi ) \cE^{-1} \Omega(\xi ) ) B = 0
$$
where $i \Omega(\xi )$ is the symbol of the operator $\rot$.
We have
$$
A(\xi ) :=  \Omega(\xi ) \cE^{-1} \Omega(\xi )  = 
\left(\begin{array}{ccc}
- \alpha_2 \xi_3^2 & \alpha_3 \xi_1 \xi_2 & \alpha_2 \xi_1
\xi_3 \\
\alpha_3 \xi_1 \xi_2 & - \alpha_1 \xi_3^2 - 
\alpha _3 \xi_1^2 & \alpha_1 \xi_2 \xi_3 \\
\alpha_2 \xi_1 \xi_3 & \alpha_1 \xi_2 \xi_3 &
- \alpha_1 \xi_2^2 - \alpha_2 \xi_1^2 
\end{array}\right)\, , 
$$
$$ 
\det (\tau^2 + A(\xi) ) \, =\, 
 \tau^2 \big( \tau^4 - \Psi(\xi) \tau^2 + 
\vert \xi \vert^2 \, \Phi(\xi)\big)
$$ 
with 
$$ 
\left\{\begin{aligned}
\Psi(\xi) 
& = (\alpha_1 + \alpha_2) \xi_3^2 +
(\alpha_2 + \alpha_3) \xi_1^2 + (\alpha_3 + \alpha_1)\xi_2^2
\\
\Phi(\xi) 
& = \alpha_1 \alpha_2 \xi_3^2 + \alpha_2 \alpha_3 \xi_1^2 +
\alpha_3 \alpha_1 \xi_2^2\, . 
\end{aligned} \right. 
$$ 
The equation $\det (\tau^2 + A(\xi) = 0$ yields
$\tau = 0$ and an equation of second order in $\tau^2$, of which
the discriminant is
$$
\Psi^2(\xi) - 4 \vert \xi \vert^2 \Phi(\xi)  = P^2 + Q 
$$
with
$$
\begin{aligned}
P &= (\alpha_1- \alpha_2) \xi_3^2  + (\alpha_3 - \alpha_2) \xi_1^2
+ (\alpha_3 - \alpha_1) \xi_2^2
\\
Q  &= 4 (\alpha_1-\alpha_2) (\alpha_1- \alpha_3) \xi_3^2 \xi_2^2
\ge 0\, . 
\end{aligned}
$$ 
The root is double if and only if
\begin{equation}
\label{racdble} 
\xi_2 = 0, \, \quad 
 \alpha_1 \xi_3^2 +   \alpha_3  \xi_1^2  =   \alpha_2
(\xi_1^2 + \xi_3^2)\, = \tau^2\, . 
\end{equation} 

Let $(\utau , \uxi) $ be a solution of \eqref{racdble}. 
For $\theta$ such that
$\alpha_3  \cos^2 \theta + \alpha_1 \sin^2 \theta = \alpha_2   $, 
we have
\begin{equation} 
\label{freqrefracco}
\uxi = \lambda b_1 \, , \quad \utau = \pm
\lambda \sqrt{\alpha_2}\, . 
\end{equation} 
where we have introduced the basis of $\RR^3$ 
$$ 
b_1 = \left( \begin{array}{c}
\cos \theta \\0 \\ \sin \theta  \end{array}\right)\, ,
\quad  b_2 = \left( \begin{array}{c}
0 \\1 \\ 0 \end{array}\right)
\, , \quad 
b_3 = \left( \begin{array}{c}
- \sin \theta \\0 \\ \cos \theta  \end{array}\right)
$$ 
In \eqref{freqrefracco}, suppose now that $\utau =   \lambda \sqrt{\alpha_2}\,
$, the other case being similar.
The kernel of $ L(\utau,  \uxi)$ is of dimension two, generated by
\begin{equation} 
u_2 := \left( \begin{array}{c} b_2 \\e_2
\end{array}\right)  , \quad 
u_3 := \left( \begin{array}{c} b_3 \\e_3
\end{array}\right)\quad \mathrm{with} \quad 
e_2 = \left( \begin{array}{c}
- \delta_1 \sin \theta  \\0 \\ \delta_3 \cos \theta 
\end{array}\right)\, ,
  \  
e_3 = \left( \begin{array}{c}
 0 \\ - \delta_2 \\ 0 \end{array}\right) 
\end{equation} 
and $\delta_j := \alpha_j / \sqrt{\alpha_2}$. 
As $L(\tau, \xi)$ is self-adjoint, $(u_2, u_3)$
also form a basis of $\ker L^*(\utau, \uxi)$. 
The tangent system  $\underline L'$  at 
$(\utau, \uxi)$  is obtained  by multiplying on the left by a basis
for the left kernel and on the right by a basis for
the right kernel of $L(\tau, \xi)$,  see  Remark \ref{groupeq}. 
We obtain 
a $2 \times 2$ system with symbol 
  \begin{equation} 
 \left(\begin{array}{cc}
\tau - ( \xi_1 \delta_3 \cos \theta + \xi_3 \delta_1 \sin\theta )
&   \xi_2 (\delta_1- \delta_3) \cos\theta \sin\theta \\
 \xi_2 (\delta_1- \delta_3) \cos\theta \sin\theta & 
\tau - ( \xi_1 \delta_2 \cos \theta + \xi_3 \delta_2 \sin\theta )
\end{array}
\right).
\end{equation} 
Compare to equations \eqref{defv} and \eqref{212}  
obtained in the semisimple and geometrically regular cases.
By a change of variables, we arrive at the form
\begin{equation}
\label{B9}
 \left(\begin{array}{cc}
\tilde \tau  - \tilde \xi_1  & \tilde \xi_2 \\
\tilde \xi_2   & \tilde \tau  +  \tilde \xi_1 
\end{array} \right)\, . 
\end{equation}
This shows that the eigenvalue is linearly splitting.


 \bigbreak
 \section{Appendix C. The block structure condition}

We consider the system
\begin{equation}
\label{bl1}
\D_x + i G (p, \zeta) 
\end{equation}
With $\zeta = (\tau - i \gamma, \eta)$.  Assume that 
\begin{ass}
\label{assbl1}
The characteristic polynomial in $\xi$, 
$\Delta (p, \zeta, \xi) = \det \big(\xi \Id + G(p, \zeta )\big) $,  has real coefficients
when $\gamma = 0$ and  has no real roots when $\gamma > 0$. 
\end{ass}

\begin{defi}
\label{Ddefbl2}
$G$ has the block structure property near $(\up, \uzeta)$ if there exists a smooth matrix 
$V$ on a neighborhood of that point such that 
$V^{-1} G V = \mathrm{diag} (G_k) $ is block diagonal, with blocks $G_k$
of size $\nu_k \times \nu_k$,  holomorphic in $\tau - i \gamma$,   having one 
the following properties:

 i) 	  the spectrum of $G_k(p, \zeta) $ is contained in $\{ \im \mu \ne 0 \}$. 

ii)  $\nu_k = 1$,  $G_k(p, \zeta)$  is real when when 
$\gamma = 0$, and  
$ \partial_\gamma G_k(\up, \uzeta)  \ne 0$, 

iv)   $\nu_k > 1$, $G_k(p, \zeta)$  has real coefficients   
  when $\gamma = 0$,   there is
 $ \underline \mu_k  \in \RR$ such that
 \begin{equation}
\label{Dbl2}
 G_k(\up, \uzeta)    = \underline \mu_k \Id +  
\left[\begin{array}{cccc}
0  & 1 & 0&   
\\
0  &0  & \ddots  &  0
\\
  & \ddots &      \ddots & 1 \\
 &  &  \cdots & 0
\end{array}\right]\, , 
\end{equation}
 and  the lower left hand corner of 
${\partial_\gamma  Q_k} (\up, \underline \zeta)$ 
does not vanish. 

\end{defi}

\begin{theo}
\label{theoC3}
The block structure condition is always satisfied when $\uga > 0$.  
When $\uga = 0$, it holds if and only if 
for all real root $\uxi$ of $\Delta (\up, \uzeta, \uxi) = 0$,

i)  there  are smooth functions $\lambda_j(p, \eta, \xi)$, analytic in $\xi$,   
real when $\xi $ is real, near $(p, \ueta, \uxi)$ such that
\begin{equation}
\label{bl3one}
\Delta(p, \zeta, \xi) =  e(p, \zeta, \xi)
\prod_{j=1}^{\underline j}  \big( \tau - i \gamma + \lambda_j  (p, \eta, \xi) \big)
\end{equation}
with $e(\up, \uzeta, \uxi) \ne 0$ and 
$\lambda_j(\up, \ueta, \uxi) + \utau = 0$, 

ii) there are smooth vectors $e_j(p, \eta, \xi)$  near $(p, \ueta, \uxi)$, analytic in $\xi$, 
linearly independent, such that
\begin{equation}
\label{bl3}
\big(  \xi  \Id +    G^{(j)} (p, \eta, \xi)   \big) e_j (p, \eta, \xi)  = 0
\end{equation}
where $ G^{(j)}  $ is the matrix $G$ evaluated at 
$ \tau - i \gamma = -  \lambda_j (p, \eta, \xi) $. 
\end{theo}

\begin{proof}
The first statement is clear from Assumption \ref{assbl1}. Moreover, one can always
make a first block diagonal reduction with blocks $G_k$ associated to 
distinct eigenvalues of $G(\up, \uzeta)$. 
Non real roots always satisfy condition $i)$ of  Definition \ref{Ddefbl2}.
So is is sufficient to consider separately each block, associated 
to one real eigenvalue.  Thus, it is sufficient to prove the theorem when 
\begin{equation}
\label{bl4}
\uxi \Id  +  G(\up, \uzeta)  = 0, \quad  \uxi \in \RR. 
\end{equation}

\medbreak
{\bf a) }  
Suppose that the block structure condition holds. Then 
\begin{equation}
\label{bl5}
\Delta (p, \zeta, \xi) = \prod \Delta_k (p, \zeta, \xi) , 
\quad  \Delta_k (p, \zeta, \xi) := \det \big(\xi \Id + G_k(p, \zeta)\big). 
\end{equation}
The form \eqref{Dbl2} of $G_k$, implies that 
$$
\det  \big (\uxi + G \up, \utau - i \gamma, \ueta)\big)   = (-1)^{\nu_k -1} c
\gamma + O(\gamma^2)
$$
where $c$ is the lower left hand corner of $\D_\gamma Q_k (\up,  \uzeta)$. 
Since $c \ne 0$ and $Q_k$ is holomorphic in 
$\tau - i \gamma$, this implies that 
\begin{equation*}
\D_ \tau \Delta_k(\up, \uzeta, \uxi) \ne 0. 
\end{equation*}
Since $\Delta_k$ is real when $\xi$ and $\tau - i \gamma $ are real, 
the implicit function theorem  implies that there is a smooth function $\lambda_k(p, \eta, \xi)$, 
analytic in $\xi$ such that 
\begin{equation}
\label{bl6}
\Delta_k(p,  \zeta, \xi) =  e(p, \zeta, \xi) \big( \tau - i \gamma + 
\lambda_k(p, \eta, \xi)
\end{equation}
on a neighborhood of $(\up, \uzeta, \uxi)$, with $e(\up, \uzeta, \uxi) \ne 0$.

Decompose the matrix $V$ into   $(V_1, \ldots, V_k)$ 
where $V_k$ is of size $N  \times \nu_k$
and denote by $V^{(k)}_k(p, \eta, \xi)$ the matrix
$V_k$ evaluated at $\tau - i \gamma = - \lambda_k(p, \eta, \xi)$. 
At the base point 
\begin{equation}
\label{bl7}
 V^{(k)}_k (\up, \ueta, \uxi) = V_k (\up, \uzeta). 
\end{equation}

If $\nu_k = 1$, then the coefficient 
$\xi + Q_k(p, \zeta) $ vanishes at when $\tau - i \gamma = - \lambda_k$
and   $ r_k =   V^{(k)}_k$ is a smooth element in the kernel of 
$\xi +   G^{(k)} (p, \eta, \xi)$.

If $\nu_k  > 1$, 
consider next the first $\nu_k - 1 $ rows of the eigenvector equation
$(\xi \Id + G_k ) u = 0$. Denoting by $v_1$ the first component and by $v'$ the remaining 
$\nu_k - 1$, these equations read 
\begin{equation}
\label{bl8}
Q (p, \zeta, \xi ) v'  =  v_1 h(p, \zeta, \xi)
\end{equation}
with $Q(\up, \uzeta, \uxi) = \Id$ and $h(\up, \uzeta, \uxi) = 0$. 
Therefore, one can solve $v' =  v_1 Q^{-1} h $. 
When $\tau - i \gamma = - \lambda_k(p, \eta, \xi)$, the determinant 
of $ \xi \Id + G^{(k)}_k (p, \eta, \xi)$  vanishes. 
This shows that the kernel of 
$ \xi \Id + G^{(k)}_k (p, \eta, \xi)$ has dimension one and is spanned 
by the smooth vector $e_k = {}^t (1, Q^{-1} h )$ where the functions are evaluated 
at $\tau - i \gamma = - \lambda_k(p, \eta, \xi)$.  
This implies that 
$ r_k := V^{(k)} e_k$ is a smooth vector in the kernel of 
$\xi \Id + G^{(k)} $.  
By \eqref{bl7}, the vectors $r_k$ are linearly independent 
at $(\up, \ueta, \uxi)$, and thus remain independent in the vicinity of that point.

\medbreak
{\bf b) }   Conversely, suppose that \eqref{bl3} and \eqref{bl4} are satisfied. 
By \eqref{bl4}, there is only one real eigenvalue $- \uxi $
at $(\up, \uzeta)$. 
For all $j$, there is an integer $\nu_j \ge 1$ such that 
\begin{equation}
\label{bl9}
\D_\xi \lambda_j = \ldots = \D_\xi^{\nu_j-1} \lambda_j = 0 , \quad 
\D_\xi^{\nu_j} \lambda_j \ne 0 \quad
\mathrm{at}  \quad (\up, \ueta, \uxi). 
\end{equation}
Repeating the analysis in \cite{MetLMS}, this implies that 
\begin{equation}
\label{bl10}
\tau - i \gamma + \lambda_j( p, \ueta, \xi) = 
e_j(p, \zeta, \xi) \Delta_j (p, \zeta, \xi)
\end{equation}
with $e_j(\up, \uzeta, \uxi) \ne 0$ and $ \Delta_j (p, \zeta, \xi)$ a monic polynomial
in $\xi$, of degree $\nu_j$, with real coefficients when 
$\gamma = 0$,  such that 
\begin{equation}
\label{bl11}
\Delta_j( \up, \uzeta, \xi ) = (\xi - \uxi)^{\nu_j}. 
\end{equation}

Suppose that 
$\nu_j = 1$.  Then, $\Delta_j = \xi  - \mu_j (p, \zeta)$ with 
$\mu_j$ smooth and 
$\tilde e_j (p, \zeta) = e_j (p, \eta, \mu_j)$ is a smooth eigenvector
of $ G(p, \zeta)$ with eigenvalue $\mu_j$. 
Comparing with \eqref{bl10}, implies that 
$\D_\gamma \mu_j (\up, \uzeta) \ne 0$. 
Thus, the one dimensional space 
$\CC \tilde e_j  $ provides us with a one dimensional block 
which satisfies property $ii)$ of Definition \ref{Ddefbl2}.

Suppose next that $\nu_j > 1$. 
Differentiate the equation \eqref{bl3} with respect to $\xi$, at 
$(\up, \ueta, \uxi)$.  Since 
$G^{(j)} (\up, \ueta, \xi) = G(\up, \uzeta) + O\big( (x - \uxi)^{\nu_j}\big) $, 
this implies that 
\begin{equation}
\label{bl12}
\underline e_{j, k} = \frac{(-1)^k} {k!} (\D_\xi^k e_j ) (\up, \ueta, \uxi)
\end{equation}
satisfy 
\begin{equation}
\label{bl13}
\begin{aligned}
&  (\uxi + G(\up, \uzeta)) \underline e_{j, 0}  = 0,  
\\
&  (\uxi + G(\up, \uzeta)) \underline e_{j, k} = -  \underline e_{j, k-1} ,  
\quad  k = 1, \ldots, \nu_j- 1. 
\end{aligned}
\end{equation}
The definition of the $e_{j,k}$ is extended to 
$(p, \zeta)$ in a neighborhood  of $(\up, \uzeta)$, as in 
\cite{MetLMS}:
\begin{equation}
\label{bl14}
e_{j, k} (p, \zeta) = \frac{(-1)^k (\nu_j - k -1)!}{2 i \pi  \nu_j!} 
\int_{ \vert \xi - \uxi \vert = r}  
\frac{ \D_\xi^{k+1} \Delta_j(p, \zeta, \xi)}{\Delta_j(p, \zeta, \xi)} e_j(p, \eta, \xi) d\xi,  
\end{equation}
where $r > 0$ is such that $\Delta_j(p, \zeta, \xi) $ has no root on the circle 
$\{  \vert \xi - \uxi \vert = r \}$. 
By \eqref{bl11}, the two definitions agree  at $(\up, \uzeta)$, so that 
$e_{j, k} (p, \zeta)  = \underline e_{j,k}$. 
Repeating the analysis of \cite{MetLMS},  one shows that the space 
$\EE_j(p, \zeta) $ spanned by the $e_{j, k} (p, \zeta)$ is invariant by 
$G(p, \zeta)$. It is smooth and of dimension $\nu_j$. 

Moreover,   using  \eqref{bl13} and the independence of the 
$e_j(\up, \uzeta)$, one shows that 
$ e_{j, k}$, for $k \in \{0, \ldots, \nu_j -1\}$ and $j  \in \{ 1, \ldots, \underline j\}$
are linearly independent. Therefore, this provides a smooth  decomposition 
of $\CC^N$ into invariant subspaces of $G$, thus a block decomposition 
of $G$. Moreover
the characteristic polynomial of 
$Q_j$ is $\Delta_j$. The property 
\eqref{Dbl2} 
is stated in \eqref{bl13}. 

The proof goes on as in \cite{MetLMS}. For each block, one can perform an additional change of  bases 
such that 
$$
 Q_k   =  Q_k(\up, \uzeta)  + 
\left[\begin{array}{cccc}
  q_1(z)  &  0  & \cdots &0
\\
 \vdots & \vdots &  \ddots & \vdots 
\\
q_\nu(z) & 0 &\ldots &   0  
\end{array}\right]\, . 
$$
(see \cite{Ral}). Thus,   
$$
\Delta_j (p, \zeta,  \xi)  
(\xi- \underline \xi)^\nu + \sum_{l = 1}^\nu  (-1)^l
q_{l } (z)\, (\xi- \underline \xi)^{\nu - l} \,   
$$
Since $\Delta_j $ has real coefficients when $\gamma = 0$, 
this implies that
the $q_l(z)$   and therefore   
$  Q_k(p, \zeta)$ is  real when $\gamma = 0$. 
In addition,  by \eqref{bl10}, 
$$
\frac{\partial \Delta_j}{\partial \gamma}  
  (\up, \uzeta,  \underline \xi) \, = \, 
\frac{\partial q_\nu}{\partial \gamma} (\up, \uzeta)\,
\ne 0 . 
$$
 Therefore, $\widetilde Q(z)$ satisfies the property
$iv)$ of Definition~\ref{Ddefbl2}   and the proof  of 
 Theorem~\ref{theoC3} is   complete.     
\end{proof} 


 \section{Appendix D.  $2 \times 2$ sytems}

Consider    the $2 \times 2$ complex system in space dimension 2: 
\begin{equation}
\label{ts1}
\tau \Id  + \xi A(p) + \eta B(p) 
\end{equation}

\begin{ass}
\label{E1}
The system is strictly hyperbolic in the direction $dt$ and 
noncharacteristic and nonhyperbolic  in the direction $dx$. 
\end{ass}

\subsection{Reductions}

In particular, the eigenvalues of $A$  are real and distinct. Thus, there is a smooth change
of basis such that $A$ is diagonal, with  nonvanishing real  diagonal coefficients $a_1 \ne a_2$. 
One can perform changes of variables which preserve 
the directions $dt$ and $dx$:
 \begin{equation}
\label{ts2}
\tau' = \tau + \alpha \eta, \quad \xi' = \xi + \beta \eta , \quad \eta' = \eta. 
\end{equation}
This transform \eqref{ts1} to 
\begin{equation}
\label{ts3}
\tau' \Id  + \xi' A(p) + \eta' B'(p) . 
\end{equation}
With appropriate choices of $\alpha(p)$ and $\beta(p)$, one can cancel the 
the diagonal terms in $B'$: 
\begin{equation}
\label{ts4}
B' = \begin{pmatrix}
    0  &   b \\
   c   &  0
\end{pmatrix}
\end{equation}
Strict hyperbolicity implies that $bc$ is real and positive. By a change of basis, one 
can assume that $b = c $ is real, so that the system reads
\begin{equation}
\label{ts5}
\tau' \Id  + \xi' \begin{pmatrix}
    a_1  & 0   \\
    0  &  a_2  
\end{pmatrix} + \eta   \begin{pmatrix}
     0 &  b  \\
     b  &   0
\end{pmatrix}
\end{equation}
The assumption implies that $a_1 a_2 < 0$. 
With  $a = ( - a_1/a_2)^{1/2}  > 0$ and 
\begin{equation}
\label{ts6}
\tau'' = \tau' , \quad \xi'' =  \frac{a}{  a_1 }  \xi'  , \quad \eta'' =   b \eta,  
\end{equation}
we obtain the reduced form 
\begin{equation}
\label{ts7}
\tau'' \Id  + \xi''\begin{pmatrix}
    a   & 0   \\
    0  &  - 1/a   
\end{pmatrix} + \eta''   \begin{pmatrix}
     0 &  1  \\
     1  &   0
\end{pmatrix}. 
\end{equation}
It remains only one parameter $a = a(p) > 0$.

\subsection{Negative spaces} 

Dropping the $\, ''\, $ we consider \eqref{ts7}, and more precisely 
\begin{equation}
\label{ts8}
\xi \Id + G(p, \tau, \eta) := \xi \Id + \begin{pmatrix}
     \tau / a & \eta/a     \\
      - a \eta  &   - a \tau
\end{pmatrix}
\end{equation}
The eigenvalue equation for $G$ is 
\begin{equation}
\label{ts9}
\mu^2 + a \tau \mu - \frac{1}{a} \tau \mu = \tau^2 - \eta^2. 
\end{equation}
If $\mu + a \tau \ne 0$, the corresponding eigenspace is  
\begin{equation}
\label{ts10}
\EE  = \CC r  , \quad  r = \begin{pmatrix}
  1  \\
      - a \eta  / ( \mu + a \tau   )
\end{pmatrix}. 
\end{equation}
For $\im \tau < 0$, there is one eigenvalue $\mu $ in $\im \mu < 0$. 
We denote by $\EE_-(\tau, \eta)$ and $r(\tau, \eta)$ the associated 
eigenspaces and eigenvectors. 

Let $z = \eta/( \mu + a \tau)$. Then, for $z \ne 0$,  \eqref{ts9} is equivalent to 
\begin{equation}
\label{ts11}
 \big( a + \frac{1}{a} \big)  \tau  =  \eta \big(  z + \frac{1}{z} \big). 
\end{equation}
In particular, 
\begin{equation}
\label{ts12}
 ( a + \frac{1}{a} )  \im \tau  =  \eta \im  z \Big( 1  -  \frac{1}{\vert z\vert ^2} \Big). 
\end{equation}
When $ \im \tau < 0$ and $\eta \ne 0$, $\im (\mu + a \tau)  < 0 $ and 
$\im \eta z  > 0$, and \eqref{ts12} implies that $\vert z \vert <  1$.    

Conversely, for $\vert z \vert < 1$, $\im z \ne 0$, we can choose  $\eta = \pm1 $
such that $\eta \im z > 0$, implying that  $\tau$ defined by \eqref{ts11} satisfies $\im \tau < 0$. 
Then $\mu = - a \tau + \eta/ z$ is a root of \eqref{ts9} and 
\begin{equation*}
\im \mu = - a \im\tau - \frac{\eta \im z}{ \vert z \vert^2}  =  - 
\frac{\eta \im z}{\vert z \vert^2}   \frac{1 + a^2 \vert z \vert^2}{1+ a^2}  < 0. 
\end{equation*}

This shows   that the mapping $(\tau, \eta ) \mapsto z $ 
sends $\{ \im \tau < 0 ; \eta \ne 0 \} $ onto 
$\{ \vert z \vert < 1 ; \im z \ne 0 \}$. By continuity argument, or by direct computation, this implies
the following.

\begin{lem}
The union of the complex lines $\EE_-(\tau, \eta)$ for 
$ \im \tau \le 0$ and $(\tau, \eta) \ne 0$ is the cone $\Gamma  \subset \CC^2$
of vectors $ u = {}^t(u_1, u_2)$ such that 
$ a \vert u_2 \vert \le \vert u_1 \vert $. 
\end{lem}

\begin{cor}
The boundary condition $ Mu = u_1 - c u_2 = 0$ satisfies the Lopatinski condition, 
if and only if $ a \vert c \vert   < 1$.
\end{cor}

\subsection{Symmetrizers}

In the form \eqref{ts7}, or \eqref{ts5} there is a unique symmetrizer (up to a scalar factor)
 which is the identity matrix. 
Tracing back to the original system \eqref{ts1}, 
yields a unique (up to scalars) symmetrizer $S(p) = V^*(p)  V (p) \underline S $ where 
$V$ is the change of basis. The symmetrizer for $G $
is $\Sigma = S A $. 

\begin{prop}
\label{propts14}
The boundary condition $ Mu =  0$ satisfies the Lopatinski condition
if and only if it is maximal dissipative. 
\end{prop}

\begin{proof}
This is invariant by  the change of basis. For systems \eqref{ts7}, maximal dissipative and 
uniform Lopatinski conditions are such that $\ker M = \{ u_1 = c u_2 \}$. 
Since 
$$
( \Sigma u, u) = a \vert u_1\vert^2 - \frac{1}{a} \vert u_2 \vert^2, 
$$
the boundary condition is maximal dissipative if and only if this form is definite negative on $\ker M$,  that is if and only if $  a^2 \vert c \vert ^2 < 1$. 
\end{proof}


\subsection{  $2 \times 2$ linearly splitting blocks}

We are given a $2 \times 2$ matrix $G(q, \tau, \eta, \gamma)$, 
$\zeta = (\tau, \eta, \gamma)$ in a neighborhood of $0$ in $\RR^3$, and 
$q$ a set of parameters in a neighborhood of $\uq$. 
with the property that  
\begin{equation}
\label{ts21}
G(q, \zeta ) =  \phi(q) \Id + \tilde G(q, \zeta),  \quad  \tilde G(q, 0) = 0.  
\end{equation}

\begin{ass}
\label{assts21}
i)    The characteristic polynomial in $\xi $,   $\Delta (q, \zeta, \xi) = 
\det \big( \xi \Id + G(q, \zeta)$,  has real coefficients when $\gamma = 0$. 

ii) The first order expansion in $\zeta$ of $ \tilde G$ at $\uq$ is 
\begin{equation}
\label{ts22}
\underline G' (\dot \zeta) = (\dot \tau - i \dot \gamma) \underline G_0  +  \dot \eta \underline G_1
\end{equation}
where $\det \underline G_0   <  0$ and 
$\tau \Id + \underline G_0 ^{-1} (\xi + \eta \underline G_1)$ is strictly hyperbolic. 
\end{ass}

\begin{prop}
\label{propts22}
There is a smooth self adjoint matrix $\Sigma(p, \zeta)$ on a neighborhood
of $(\uq, 0)$ such that 

\quad i)  $ - \im (\Sigma G ) =  \gamma  E $, with $E(\uq, 0)$ positive definite,

\quad ii) with $\underline \Sigma = \Sigma (\uq, 0)$, 
$\underline \Sigma  \underline G_0  $ and  $\underline \Sigma  \underline G_1  $
are self adjoint.  
\end{prop}

\begin{proof}
{\bf a) } We can assume that $\underline G_0 $ is diagonal. Introduce the notations
\begin{equation}
\label{ts23}
\tilde G =    \begin{pmatrix}
     a &  b  \\
    c  &  d
\end{pmatrix}, \quad 
 \underline G_0 = \begin{pmatrix}
   \underline a_0    & 0    \\
   0   &   \underline d_0
\end{pmatrix}, \quad 
 \underline G_1 = \begin{pmatrix}
   \underline a_1    & \underline b_1    \\
   \underline c_1    &   \underline d_1
\end{pmatrix}. 
\end{equation}
The assumption implies that 
\begin{equation}
\label{ts24}
\begin{aligned}
&\underline a_0 = \D_ \tau a(\uq, 0)
\  \mathrm{and} \  \underline d_0 = \D_\tau b(\uq, 0) 
\  \mathrm{are \ real \ and } \ \ \underline a_0 \underline d_0 < 0 , 
\\
& 
\underline b_1 = \D_\eta b(\uq, 0), \ \underline c_1 = \D_\eta c(\uq, 0) , \ 
\underline b_1 \underline c_1 \  \mathrm{is  \ real \ and } \   \underline b_1 \underline c_1  <  0 . 
\end{aligned}
\end{equation}
By assumption, there holds
\begin{equation}
\label{ts24b}
a + d \in \RR \,,  \quad  ad - bc \in \RR \quad \mathrm{when} \ \gamma = 0. 
\end{equation}
  We look for a symmetrizer   
   \begin{equation}
\label{ts25}
\Sigma =    \begin{pmatrix}
     \alpha  &  \beta  \\
    \overline \beta   &  \delta
\end{pmatrix},  \qquad \alpha , \delta \in \RR, \ \beta \in \CC. 
\end{equation}

 {\bf b) } We show that one can choose  such that
 \begin{equation}
\label{ts26}
\im (\Sigma G ) = 0  \quad  \mathrm{when} \quad  \gamma = 0. 
\end{equation}
In this part, we assume that $\gamma = 0$.  \eqref{ts26}  is equivalent to the conditions
 \begin{eqnarray}
 \label{ts27}
\im (\alpha a + \beta c )  & =  & 0  , 
\\
\label{ts28}
\im ( \overline \beta b +  \delta d ) & = & 0, 
\\
\label{ts29}
 (\overline a -   d ) \beta & =  & \alpha b - \delta \overline c.      
\end{eqnarray}
 Because of  \eqref{ts24b}, $\sigma : =  \overline a - d = \re (a - d) $ is real and 
 by \eqref{ts24}, $\D_\tau \sigma (\uq, 0) \ne 0$. 
 Therefore, we can take $(q, \sigma, \eta)$ as local coordinates  near 
 $(\uq, 0, 0)$, and we first solve the equation 
 \begin{equation}
\label{ts211}
\alpha b - \delta \overline c = 0 , \quad \mathrm{on } \ \sigma = \gamma = 0. 
\end{equation}
 By \eqref{ts24}, there holds
 \begin{equation}
\label{ts212}
b_{ \vert \sigma = 0 , \gamma = 0} = \eta  
\underline b_1 + O(\eta^2), 
\quad 
c_{ \vert \sigma = 0 , \gamma = 0} = \eta \underline c_1 + O(\eta^2)
\end{equation}
 Moreover, when $\sigma = 0$, $ a d =  \vert a \vert^2 $ and   
  by \eqref{ts24b}, $bc $ is real. Therefore, factoring out $\eta $ in \eqref{ts211}, 
  we see that there is a smooth local solution 
  $\alpha$,  $ \delta $, such that 
  \begin{equation}
\label{ts213}
\alpha (q, \zeta) \in \RR, \quad  \delta(q, \zeta) \in \RR, \quad 
\alpha(\uq, 0) = \vert \underline c_1\vert^2,\quad  \delta(\uq, 0) = \underline b_1 \underline c_1. 
\end{equation}
 By \eqref{ts211}, one can factor out $\sigma $  in 
 $\alpha b - \delta \overline c $ restricted to  $\gamma= 0$ and there is a smooth 
 function $\beta$ such that  \eqref{ts29} holds when  $\gamma = 0$. 
 In particular, since $b$ and $c$ 
 \begin{equation}
\label{ts214}
\beta(\uq, 0) = 0. 
\end{equation}
For $\gamma = 0$ and $\sigma \ne 0$, 
 $$
 \im (\alpha a + \beta c ) = \frac{\alpha}{\sigma}  \im  ( a \sigma +  {bc} ) = 
  \frac{\alpha}{\sigma}  \im  (    {bc} - ad )
 $$
 vanishes by \eqref{ts24b}. Thus, \eqref{ts27} and similarly, \eqref{ts28} are satisfied on 
 $\gamma = 0$.  
  Therefore, we have constructed a smooth solution 
  $\Sigma(q, \zeta) $ of \eqref{ts26}.

  Morever, differentiating this equation with respect to $\tau$ and $\eta$, implies that \begin{equation}
\label{ts215}
\im (\underline \Sigma \underline G_0 ) = 0, \quad 
\im (\underline \Sigma \underline G_1) = 0. 
\end{equation}
The first equality,    \eqref{ts23} and  \eqref{ts24}, imply  that 
  $\beta (\uq, 0) = 0$. Thus,  with \eqref{ts213}, we see that 
  \begin{equation}
\label{ts216}
\underline \Sigma = \begin{pmatrix}
   \underline \alpha    &  0  \\
   0   &  \underline \delta 
\end{pmatrix}, \quad  \mathrm{with} \quad \underline \alpha \underline \delta < 0. 
\end{equation}
    
 { \bf c) }   Changing $\Sigma$ to $- \Sigma$ if necessary, we can assume that 
 $\underline \alpha \underline a_0 >  0$, implying that 
 $\underline \delta \underline d_0  >  0$, hence   that 
 $\underline \Sigma \underline G_0 $ is definite positive.  
 
 The identity \eqref{ts26} implies that $\im (\Sigma G) =  -  \gamma E$
 with $E$ depending smoothly on $(q, \zeta)$ in a neighborhood of $(\uq, 0)$. 
 Moreover,  at $q = \uq, \tau = \eta = 0$, there holds
 \begin{equation*}
\im (\Sigma G)  =  \im ( - i \gamma \underline \Sigma \underline G_0)  + O( \gamma^2)
\end{equation*}
showing that 
$E(\uq, 0) =     \underline \Sigma \underline G_0$ is definite positive. 
   \end{proof}

   \begin{theo}
   \label{theots23}
   Suppose that  Assumption $\ref{assts21}$ is satisfied and that 
   $ M (q, \zeta)$ is a $1 \times 2 $ matrix such that the 
   equation 
   \begin{equation}
\label{ts217}
\D_x u + i G(q, \zeta) u = f , \quad  M (q, \zeta) u_{\vert x= 0 }  = g
\end{equation}
   satisfies the uniform Lopatinski condition near $(\uq, 0)$. 
   Then, there is a smooth symmetrizer $\Sigma (q, \zeta)$ near $(\uq, 0)$ 
   and there are constants $c> 0$ and $C$, such that 
   \begin{eqnarray}
   \label{ts218}
 &   & -  \im (\Sigma G ) \ge c \gamma \Id , \quad  when \  \gamma \ge 0, 
  \\
  \label{ts219}
 &   &  - \Sigma + C M^* M \ge c \Id.  
\end{eqnarray}
   
   \end{theo}

   \begin{proof}
   The symmetrizer $\Sigma$  has been constructed in Proposition \ref{propts22}. 
   To prove that \eqref{ts219} holds on a neighborhood of $(\uq, 0)$, it is sufficient 
   to prove that it holds at $(\uq, 0)$, thus it is sufficient 
   to check that $- \underline \Sigma  $ is definite positive 
   (or $\underline \Sigma$ definite negative) on $\ker \underline M$, 
   where $\underline M = M(\uq, 0)$. 
   
   In addition, since $\underline \Sigma \underline G_0$  is self adjoint and definite positive, 
   and $\underline \Sigma \underline G_1$ is self adjoint, 
   $\underline \Sigma $  is a symmetrizer for 
  $\underline G' = (\tau - i \gamma) \underline G_0 + \eta \underline G_1$. 
  Thus, by Proposition \ref{propts14}, it is sufficient to prove that 
  the boundary condition $\underline M$ satisfies the uniform Lopatinski condition
  for the operator $\D_x + i \underline G'$. 
  
  This follows from the analysis in Section 4. We have shown that the negative space 
  for $G(q, \zeta)$ with $\gamma > 0$, is the negative space for  a matrix
  $\ccG (q, \rho, \cz)$ with $\rho = \vert \zeta \vert $ and 
  $\cz = \zeta / \vert \zeta \vert$. Moreover, 
 $  \ccG(\uq, 0, \cz)  = \underline G' (\cz)$. 
 The Lopatinski condition implies that 
 \begin{equation}
\label{ts220}
\vert u \vert \le C \vert M(\uq, \rho \cz)  u \vert 
\end{equation}
   for all $u \in  \EE_-(\uq, \rho, \cz)$. As $\rho$ tends to zero, 
   this implies that 
    \begin{equation}
\label{ts221}
\vert u \vert \le C \vert M(\uq, 0)  u \vert 
\end{equation}
   for all $u \in \EE_-(\uq, 0, \cz)$ and $\vert \cz \vert = 1$ with $\check \gamma > 0$. 
   This means that $\underline M$ does satisfy the uniform Lopatinski condition
  for  $\D_x + i \underline G'$.
   \end{proof}

\bigskip
\section{Appendix E. The viscous case}

In conclusion, we briefly discuss the case of a hyperbolic
system with viscous regularization, that is,
an $N \times N$ second order system with symbol
\begin{equation}
\label{vdefL} 
\begin{aligned}
L_\nu(p, \tau , \xi) &= \tau \Id + A(p, \xi)+\nu B(p,\xi) \\
&= 
\tau \Id + \sum_{j=1}^d \xi_j A_j (p)s
 -i\nu \sum_{j,k=1}^d \xi_j\xi_k B_{j,k} (p)s\\
&=L(p, \tau , \xi) 
 -i\nu \sum_{j,k=1}^d \xi_j\xi_k B_{j,k} (p)s.
\end{aligned}
\end{equation}
We make the assumptions of {\it symmetrizability} of the
first-order part, i.e., existence of positive definite
symmetric $S$ such that $SA_j$ is symmetric for all $j$,\
and {\it dissipativity} of the second-order part, i.e.
$B(p,\xi)=\sum_{j,k}\xi_j\xi_k SB_{j,k}\ge 0$, with no eigenvector
of $A(p,\xi)$ lying in the kernel of $B(p,\xi)$
(the ``genuine coupling'' condition of \cite{Kaw}). 

\medbreak
{\bf The small-viscosity limit.}
In \cite{MeZ.1, GMWZ.2, GMWZ.3, GMWZ.4}, there
was considered, under the additional assumption of constant
multiplicity of the first-order part, the problem of obtaining
maximal estimates
 \begin{equation}
\label{veq54}
\gamma \Vert u \Vert^2_{L^2(\RR_+) } + 
\vert u(0) \vert^2 \ls \frac{1}{\gamma + \nu |\zeta|} \Vert f \Vert^2_{L^2(\RR_+) } + 
\vert g  \vert^2 
\end{equation}
generalizing those of \eqref{eq54}, where $\nu$ is the coefficient
of viscosity (``Reynolds number'') in \eqref{vdefL} and
$\zeta=(\tau, \eta)$ as in previous sections denotes Laplace--Fourier
transform frequencies, for $\nu |\zeta|$ sufficiently small.
Such small-frequency estimates, together with analogous intermediate-
and high-frequency estimates {\it not} depending on the constant-multiplicity
assumption, 
were used to verify, respectively, existence of viscous boundary
and shock layers in the $\nu\to 0$ limit, 
with rigorous convergence to associated formal asymptotic series.

A consequence of the matrix perturbation analysis of \cite{MeZ.1} is 
that the reduced equations $G^1_{viscous}$ for the viscous system
for a hyperbolic mode associated with basis $V$ satisfies the
intuitively appealing relation
\begin{equation}
V G^1_{viscous}= GV=\big(G_{inviscid} +\nu B(p, \xi)\big) V,
\end{equation}
where $V=R(p,\xi)$ is a matrix of right eigenvectors 
of $A(p, \xi)$ associated with the single
eigenvalue $\tau$, and thus $V^*S=L(p,\xi)$ is a matrix of  
left eigenvectors, where $S$ is the symmetrizer of the system.

A well-known consequence of the dissipativity/genuine coupling
assumption (see especially \cite{Kaw, KSh}) is
that $LBR(p,\xi)=V^*SBV$ is uniformly positive definite.
{}From this fact, it is straightforward to see that, under the 
structural hypotheses of Assumption \ref{ass41}, 
the smooth hyperbolic symmetrizer $\Sigma=V^*S$ constructed in Section
6 near points of variable multiplicity serves also as a symmetrizer
for the viscous system, yielding the desired maximal estimate \eqref{veq54}.
{\it This extends the results of the above papers to the variable-multiplicity
case under Assumption \ref{ass41}, 
in particular yielding existence of boundary and shock layers for MHD
under the uniform Lopatinski (Evans function) condition.}

We remark that a continuity argument like that in Section 7
shows that the uniform Evans function condition reduces in
the vanishing-magnetic field ($H\to 0$) limit to the corresponding
condition on the limiting fluid-dynamical shock.
The uniform Evans function condition is always satisfied
for sufficiently small-amplitude shocks in either gas dynamics
or MHD \cite{PZ, FS}.

\medbreak
{\bf Long-time stability of viscous shock waves.}
In \cite{ZKochel, Zhandbook, ZCime, GMWZ.2},
there was considered under the constant-multiplicity assumption
together with an additional technical hypothesis (satisfied for
gas dynamics and MHD) that the 
glancing set have a certain ``foliated structure'',
the related problem of time-asymptotic $L^1\cap H^s\to L^2\cap H^s$
stability, for $s$ sufficiently large, of planar viscous shock
profiles, for which the relevant symbol is \eqref{vdefL} with
constant coefficients and $\nu=1$.
Again, the constant-multiplicity hypothesis was used only in
the small-frequency regime, this time to obtain maximal
$L^1\to L^p$ stability estimates, $p\ge 2$. 
Away from the glancing set, the only way in which the constant-multiplicity
assumption was used was to show that the real part of ``slow'',
or ``hyperbolic'' eigenvalues was bounded 
below by multiples of $\gamma + |\zeta|$.
But, it is a straightforward exercise to show that this is also implied 
by the existence of a symmetrizer $\Sigma G^1 \ge C^{-1}(\gamma +|\zeta|)$.
{\it This extends the results of the above papers to the variable-multiplicity
case under Assumption \ref{ass41}, 
in particular yielding long-time stability of Lax- and overcompressive
type shock waves for MHD
under the refined Lopatinski (Evans function) condition defined
therein.}


\end{document}